\documentclass{guthesis}
\usepackage[all]{xy}
\usepackage{makeidx}
\makeindex
\xyoption{2cell}
\CompileMatrices  
\CompilePrefix{matrices/matrix}
\UseTwocells
\UseHalfTwocells
\usepackage{amsmath}
\usepackage{amsthm}
\usepackage{epsfig}
\DeclareSymbolFont{AMSb}{U}{msb}{m}{n}
\DeclareMathSymbol{\natural}{\mathbin}{AMSb}{"4E}
\DeclareMathSymbol{\integer}{\mathbin}{AMSb}{"5A}
\DeclareMathSymbol{\real}{\mathbin}{AMSb}{"52}
\DeclareMathSymbol{\rational}{\mathbin}{AMSb}{"51}
\DeclareMathSymbol{\I}{\mathbin}{AMSb}{"49}
\DeclareMathSymbol{\complex}{\mathbin}{AMSb}{"43}
\DeclareMathSymbol{\bbF}{\mathbin}{AMSb}{"46}
\DeclareMathSymbol{\bbI}{\mathbin}{AMSb}{"49}
\DeclareMathSymbol{\bbB}{\mathbin}{AMSb}{"42}
\DeclareMathSymbol{\Acat}{\mathbin}{AMSb}{"41}
\DeclareMathSymbol{\Bcat}{\mathbin}{AMSb}{"42}
\DeclareMathSymbol{\Csmall}{\mathbin}{AMSb}{"43}
\DeclareMathSymbol{\Dsmall}{\mathbin}{AMSb}{"44}
\newcommand{\C}{\ensuremath{\mathcal{C}}}
\newcommand{\D}{\ensuremath{\mathcal{D}}}
\newcommand{\V}{\ensuremath{\mathcal{V}}}
\newcommand{\calS}{\ensuremath{\mathcal{S}}}
\newcommand{\B}{\ensuremath{\mathcal{B}}}
\newcommand{\F}{\ensuremath{\bbF}}
\def\T{\mathcal{T}}
\def\S{\mathcal{S}}
\newcommand{\defterm}[1]{\textbf{#1}}
\newcommand{\cat}[1]{\ensuremath{\mbox{\upshape{\textbf{#1}}}}}
\newcommand{\Cat}{\cat{Cat}}
\newcommand{\Triv}{\cat{Triv}}
\newcommand{\CAT}{\cat{CAT}}
\newcommand{\Set}{\cat{Set}}
\newcommand{\Top}{\cat{Top}}

\newcommand{\Clone}{\cat{Clone}}
\newcommand{\Law}{\cat{Law}}
\newcommand{\FP}{\cat{FP}}
\newcommand{\Mnd}{\cat{Mnd}}
\newcommand{\Multicat}{\cat{Multicat}}
\newcommand{\SymmMulticat}{\cat{$\Sigma$-Multicat}}
\newcommand{\Multigraph}{\cat{Multigraph}}
\newcommand{\CatMulticat}{\cat{Cat-Multicat}}
\newcommand{\CatSymmMulticat}{\cat{Cat-$\Sigma$-Multicat}}
\newcommand{\CatMultigraph}{\cat{Cat-Multigraph}}

\newcommand{\VMultigraph}{\cat{\V-Multigraph}}
\newcommand{\DigraphMultigraph}{\cat{Digraph-Multigraph}}
\newcommand{\Operad}{\cat{Operad}}
\newcommand{\SymmOperad}{\cat{$\Sigma$-Operad}}
\newcommand{\CatOperad}{\cat{Cat-Operad}}
\newcommand{\CatSymmOperad}{\cat{Cat-$\Sigma$-Operad}}

\newcommand{\strictness}[1]{_{\rm #1}}
\newcommand{\SMCCat}{\cat{SMC}}
\newcommand{\wk}{\strictness{wk}}

\newcommand{\str}{\strictness{str}}
\newcommand{\Digraph}{\cat{Digraph}}

\newcommand{\Alg}[1]{\cat{Alg$(#1)$}}
\newcommand{\Algfctr}{\cat{Alg}}
\newcommand{\Algwk}[1]{\cat{Alg$\wk(#1)$}}

\newcommand{\PsAlg}[1]{\cat{Ps-Alg$(#1)$}}
\newcommand{\op}{^{\rm{op}}}
\newcommand{\WkPCat}{\cat{Wk-$P$-Cat}}
\newcommand{\StrPCat}{\cat{Str-$P$-Cat}}

\newcommand{\End}{\mathop{\rm{End}}}
\newcommand{\term}{\mathop{\rm{term}}}
\newcommand{\var}{\mathop{\rm{var}}}
\newcommand{\supp}{\mathop{\rm{supp}}}
\newcommand{\munge}{\mathop{\rm{label}}}
\newcommand{\tree}{\mathop{\rm{tree}}}
\newcommand{\Lan}{\mathop{\rm{Lan}}}
\newcommand{\union}{\cup}
\newcommand{\conc}{\mathrel{++}}
\newcommand{\nary}{$n$-ary}

\newcommand{\st}[1]{\cat{st}(\ensuremath{#1})}
\newcommand{\stfunc}{{\rm\textbf{st} }}
\newcommand{\seq}{_\bullet}
\newcommand{\udot}{^\bullet}
\newcommand{\dseq}{\seq\udot}
\newcommand{\ldseq}{_{\bullet\bullet}}
\newcommand{\tseq}{\udot\ldseq}
\newcommand{\adot}{a\seq}
\newcommand{\cmp}[1]{\circ(#1)} 
\def\hatmap{\hat{\phantom{\alpha}}}
\def\barmap{\bar{\phantom{\alpha}}}
\def\checkmap{\check{\phantom{\alpha}}}
\newcommand{\ijseq}{_{i\bullet}^j}
\newcommand{\xylabel}[1]{\ensuremath{\xymatrix{*+[o][F-]{#1}}}}
\newcommand{\Wk}[1]{{\ensuremath{\mbox{\upshape Wk($#1$)}}}}
\newcommand{\Wkwrt}[2]{{\ensuremath{\mbox{\upshape Wk$_{#2}(#1)$}}}}
\newcommand{\WkP}{\Wk{P}}
\newcommand\oneton[2]{{\def\tmp##1{#1} \tmp{1} #2 \dots #2 \tmp{n}}}
\newcommand\listn[1]{\oneton{#1}{,}}

\newcommand{\prodkn}[1]{#1_{k_1} \times \dots \times #1_{k_n}}
\newcommand{\toletter}[1]{
  \mathop{\stackrel{\scriptstyle{#1}}{\longrightarrow}}
}
\newcommand{\adjunction}[4]{  
	\xymatrix{
		#1 \ar@<1ex>[r]^-{#3}_-{\bot}
		& #2 \ar@<1.1ex>[l]^-{#4}
	}
}
\newcommand{\midlabel}[1]{
	        \xymatrixrowsep{4pc}
		\xymatrixcolsep{4pc}
		\xymatrix{ {} \ar @{}[d]^{#1} \\ {} }
}
\newcommand{\midequals}{\midlabel{=}}

\newcommand{\parallelarsdir}[3]{\ar@<0.5ex>[#3]^-{#1} \ar@<-0.5ex>[#3]_-{#2}}
\newcommand{\uparallelarsdir}[3]{\unar@<0.5ex>[#3]^-{#1}
	\unar@<-0.5ex>[#3]_-{#2}}
\newcommand{\parallelars}[2]{\parallelarsdir{#1}{#2}{r}}
\newcommand{\fork}[6]{
	\xymatrix{#1 \parallelars{#2}{#3} & #4 \ar[r]^-{#5} & #6 }
}
\newcommand{\parallelpair}[4]{ \xymatrix{#1 \parallelars{#2}{#3} & #4} }

\newcommand{\commsquare}[8]{
	\xymatrix{
		#1 \ar[r]^-{#2} \ar[d]_-{#4}
		& #3 \ar[d]^-{#5} \\
		#6 \ar[r]^-{#7}
		& #8
	}
}
\newcommand{\algmap}[6]{
	\commsquare{#2 #3}{#2 #1}{#2 #5}{#4}{#6}{#3}{#1}{#5}
}
\newcommand{\valg}[3]{
	\mbox{$\left(\begin{matrix}
		{#1}{#2} \cr \phantom{#3} \Big\downarrow {#3} \cr {#2}
	\end{matrix}\right)$}
}
\newcommand{\halg}[3]{
	\xymatrix{{#1}{#2} \ar[r]^{#3} & {#2}}
}
\newdir{(>}{{}*!/-10pt/\dir{>}}
\newcommand{\ear}{\ar@{->>}}     
\newcommand{\mar}{\ar@{(>->}}    
\let\booar=\ear                  
\let\lffar=\mar                  
\newcommand{\unar}{\ar@{..>}}  
\newdir^{((}{{}*!/-5pt/\dir^{(}}
\newcommand{\incar}{\ar@{^{((}->}} 
\newcommand{\pullback}{\ar@{}[dr]|-<{\textstyle{\lrcorner}}}
\newcommand{\pushout}{\ar@{}[dr]|->{\textstyle{\ulcorner}}}
\newcommand{\trans}{\overline}
\newcommand{\E}{\ensuremath{\mathcal{E}}}
\newcommand{\M}{\ensuremath{\mathcal{M}}}
\newcommand{\Ebar}{\overline \E}
\newcommand{\Mbar}{\overline \M}
\newcommand{\arin}{\mathop{\mbox{in}}}
\newcommand{\orth}{\mathop{\bot}}
\newcommand{\orthset}{^{\bot}}
\newcommand{\smc}{symmetric monoidal category}

\newcommand{\smcs}{symmetric monoidal categories}

\newcommand{\Fp}{{\ensuremath{F_{\rm pl}}}}
\newcommand{\Up}{{\ensuremath{U^{\rm pl}}}}
\newcommand{\Tp}{{\ensuremath{T_{\rm pl}}}}
\newcommand{\Fsig}{{\ensuremath{F_{\Sigma}}}}
\newcommand{\Usig}{{\ensuremath{U^{\Sigma}}}}
\newcommand{\Tsig}{{\ensuremath{T_{\Sigma}}}}
\newcommand{\Ffpp}{{\ensuremath{F_{\rm fp}^{\rm pl}}}}
\newcommand{\Ufpp}{{\ensuremath{U^{\rm fp}_{\rm pl}}}}

\newcommand{\Fsigp}{{\ensuremath{F_{\Sigma}^{\rm pl}}}}
\newcommand{\Usigp}{{\ensuremath{U^{\Sigma}_{\rm pl}}}}

\newcommand{\Ffpsig}{{\ensuremath{F_{\rm fp}^{\Sigma}}}}
\newcommand{\Ufpsig}{{\ensuremath{U^{\rm fp}_{\Sigma}}}}

\newcommand{\Ffp}{{\ensuremath{F_{\rm fp}}}}
\newcommand{\Ufp}{{\ensuremath{U^{\rm fp}}}}
\newcommand{\Tfp}{{\ensuremath{T_{\rm fp}}}}
\newcommand{\tr}{\mathop{\mbox{tr}}}
\def\nterms#1{(\Up \Fp #1)_n}
\def\allterms#1{\sum_n(\Up \Fp #1)_n}
\def\alleqns#1{\sum_n ((\Up \Fp #1)_n)^2}
\def\fs#1{{\underline{#1}}}	
\newcommand{\id}{{\rm id}} 
\def\n{\fs{n}}
\def\m{\fs{m}}
\let\act=\cdot
\let\kel=\circ
\let\ccomp=\bullet
\let\ocomp=\circ
\let\fcomp=\circ
\def\tcomp{}

\theoremstyle{plain}
\newtheorem{theorem}{Theorem}[section]
\newtheorem{lemma}[theorem]{Lemma}
\newtheorem{corollary}[theorem]{Corollary}
\theoremstyle{definition}
\newtheorem{defn}[theorem]{Definition}
\newtheorem{example}[theorem]{Example}
\newtheorem{remark}[theorem]{Remark}

\newcommand{\dop}{finite product operad}

\newcommand{\dm}{finite product multicategory}

\newcommand{\dops}{finite product operads}
\newcommand{\Dops}{Finite product operads}

\newcommand{\dms}{finite product multicategories}
\newcommand{\Dms}{Finite product multicategories}

\newcommand{\dco}{finite product \Cat-operad}

\newcommand{\Dcos}{Finite product \Cat-operads}

\newcommand{\dvm}{finite product \V-multicategory}

\newcommand{\DOpCat}{{\cat{FP-Operad}}}
\newcommand{\DCOCat}{{\cat{\Cat-FP-Operad}}}
\newcommand{\DCMCat}{{\cat{\Cat-FP-Multicat}}}

\newcommand{\DMCat}{{\cat{FP-Multicat}}}

\submissionmonth{September}
\submissionyear{2008}
\title{Coherence for Categorified Operadic Theories}
\author{Miles Gould}
\begin{document}
\allowdisplaybreaks
\maketitle
\begin{dedication}
To my father 
\[1 + 2 + \dots + n = {n \over 2} (n + 1)\]
\end{dedication}
\chapter*{Abstract}
\addcontentsline{toc}{chapter}{Abstract}

Given an algebraic theory which can be described by a (possibly symmetric)
operad $P$, we propose a definition of the \emph{weakening} (or
\emph{categorification}) of the theory, in which equations that hold strictly
for $P$-algebras hold only up to coherent isomorphism.
This generalizes the theories of monoidal categories and symmetric monoidal
categories, and several related notions defined in the literature.
Using this definition, we generalize the result that every monoidal category is
monoidally equivalent to a strict monoidal category, and show that the
``strictification'' functor has an interesting universal property, being
left adjoint to the forgetful functor from the category of strict
$P$-categories to the category of weak $P$-categories.
We further show that the categorification obtained is independent of our choice
of presentation for $P$, and extend some of our results to many-sorted
theories, using multicategories.

\tableofcontents
\listoffigures
\addcontentsline{toc}{chapter}{List of Figures}

\chapter*{Acknowledgements}
\addcontentsline{toc}{chapter}{Acknowledgements}
First and foremost, I must thank my supervisor, Dr Tom Leinster, for his help
and encouragement with all aspects of this thesis and my research.
I could not have wished for a better supervisor.
I would like to thank Steve Lack for invaluable help with pseudo-algebras: the
argument of Theorem \ref{thm:ps-algs are weak algs} is due to him (any errors,
of course, are mine).
I would like to thank Michael Batanin for suggesting I consider the
construction of Definition \ref{def:weakening_pres} in the context of symmetric
operads.
I would like to thank Jeff Egger and Colin Wright for invaluable motivational
advice, which in both cases came exactly when it was most needed.
I would like to thank Jon Cohen for suggesting Examples \ref{ex:ptd_sets} and
\ref{ex:ptd_sets_ABCD}, and Hitesh Jasani for helping me to see the benefits of
isolating the concept of labelling functions.
I would like to thank the night staff at Schiphol airport, for providing quite
the best environment for doing mathematics I've ever encountered.
I would like to thank Wilson Sutherland, both for his excellent teaching of
undergraduate mathematics and for encouraging me to apply for this PhD, and
Samson Abramsky and Bob Coecke, for first showing me the beauty of category
theory.
Thanks to all those who commented on drafts of this thesis:
Rami Chowdhury, Malcolm Currie, Susannah Fleming, Cath Howdle, John Kirk, Avril
Korman and Michael Prior-Jones.
Thanks to Hannah Johnson for Sumerological ratification.
I am grateful to EPSRC, for funding this research.  

Much of the challenge of this PhD has been retaining some semblance of sanity
throughout, so the people below are those who have provided welcome distraction
(as opposed to the unwelcome kind).
Thanks must go to Ruth Elliot, my co-organizer for the Scottish Juggling
Convention 2008, who did far more than her share despite being in recovery from
a serious motorbike accident; to the rest of Glasgow Juggling Club (and to Alia
Sheikh for first teaching me to juggle); to the hillwalking and rock climbing
crowd, namely Katie Edwards, Martin Goodman, Michael Jenkins, Andy Miller,
Elsie Riley, Jo Stewart, Richard Vale, Bart Vlaar, Dan and Becca Winterstein,
Stuart White, and especially Philipp Reinhard; to my office-mates, James
Ferguson, Martin Hamilton, Gareth Vaughan and John Walker, for the good times;
to the members of NO2ID Scotland, particularly Geraint Bevan, Richard Clay,
James Hammerton, Alex Heavens, Bob Howden, Jaq Maitland, Roddy McLachlan,
Charlotte Morgan, and John Welford; to the students and instructors at Glasgow
Capoeira; and to all at the theatre group Two Shades of Blue.
I would like to extend sincere and heartfelt thanks to the Chinese emperor Shen
Nung, the Ethiopian goat herder Kaldi, the Sumerian goddess Ninkasi, and the
unnamed Irish monk who, according to legend, discovered or invented tea,
coffee, beer and whisky respectively.

I would like to thank my parents, Dick and Jackie Gould, for their patience and
support, and particularly my father for showing me that mathematics could be
beautiful in the first place.
Finally, I would like to thank my wonderful girlfriend Ciorstaidh MacGlone, who
(when not contending for the computer) has been an endless source of love,
sympathy, support and tea.

\chapter*{Declaration}
\addcontentsline{toc}{chapter}{Declaration}
I declare that this thesis is my own original work, except where credited to
others.
This thesis does not include work forming part of a thesis presented
for another degree.

\setcounter{chapter}{-1}
\chapter{Introduction}

Many definitions exist of categories with some kind of ``weakened'' algebraic
structure, in which the defining equations hold only up to coherent
isomorphism.
The paradigmatic example is the theory of weak monoidal
categories, as presented in \cite{catwork}, but there are also definitions of
categories with weakened versions of the structure of groups \cite{baez+lauda},
Lie algebras \cite{baez+crans}, crossed monoids \cite{crossedmonoid},
sets acted on by a monoid \cite{modulecat}, rigs \cite{laplaza}, vector spaces
\cite{k+v} and others.
A general definition of such categories-with-weakened-structure is obviously
desirable, but hard in the general case.
In this thesis, we restrict our attention to the case of theories that can be
described by (possibly symmetric) operads, and present possible definitions of
weak $P$-category and weak $P$-functor for any symmetric operad $P$.
We show that this definition is independent (up to equivalence) of our choice
of presentation for $P$; this generalizes the equivalence of classical and
unbiased monoidal categories.
In support of our definition, we present a generalization of Joyal and Street's
result from \cite{j+s} that every weak monoidal category is monoidally
equivalent to a strict monoidal category: this holds straightforwardly when $P$
is a plain operad.
This generalization includes the classical theorem that every symmetric
monoidal category is equivalent via symmetric monoidal functors and
transformations to a symmetric monoidal category whose associators and unit
maps are identities.

The idea is to consider the strict models of our theory as algebras for an
operad, then to obtain the weak models as (strict) algebras for a weakened
version of that operad (which will be a \Cat-operad).
In particular, we do not make use of the pseudo-algebras of Blackwell, Kelly
and Power, for which see \cite{bkp}.
Their definition is related to ours in the non-symmetric case, however: we
explore the connections in Chapter \ref{ch:others}.
We weaken the operad using a similar approach to that used in Penon's
definition of $n$-category: see \cite{penon}, or \cite{cheng+lauda} for a
non-rigorous summary.

In Chapters \ref{ch:theories}, \ref{ch:operads} and \ref{ch:factsys}, we review
some essential background material on theories, operads and factorization
systems.
Most of this is well-known, and only one result (in Section \ref{sec:synclass})
is new.
In Chapter \ref{ch:categorification}, we present our definitions of weak
$P$-category, weak $P$-functor and $P$-transformation.
We start with a na\"ive, syntactic definition that is only effective for
strongly regular (plain-operadic) theories.
We then re-state this definition using the theory of factorization systems,
which allows us to apply it to the more general symmetric operads.
Section \ref{sec:symmmoncats} uses this definition to explicitly calculate the
categorification of the theory of commutative monoids with their standard
signature, and shows that this is exactly the classical theory of symmetric
monoidal categories.
In Chapter \ref{ch:coherence}, we treat the problem of different presentations
of a given operad: we use this to prove that the weakening of a given theory is
independent of the choice of presentation.
We also prove some theorems about strictification of weak $P$-categories.
In Chapter \ref{ch:others}, we compare our approach to other approaches to
categorification which have been proposed in the literature.

Material in this thesis has appeared in two previous papers: the
material on strictification for strongly regular theories was in my preprint
\cite{WkPeqStP}, and the material on signature-independence was in my paper
\cite{presindep}, which was presented at the 85th Peripatetic Seminar
on Sheaves and Logic in Nice in March 2007, and at CT 2007 in Carvoeiro,
Portugal.

\section{Remarks on notation}
Throughout this thesis, the set $\natural$ of natural numbers is taken to
include 0.
We shall occasionally adopt the ${}\seq$ notation from chain complexes and
write, for instance, $p\seq$ for a finite sequence $p_1, \dots, p_n$ and
$p\seq\udot$ for a double sequence.
We shall use the notation $\fs n$ to refer to the set $\{1, \dots, n\}$ for all $n \in \natural$: the set $\fs 0$ is the empty set.
We shall use the symbol 1 to refer to terminal objects of categories and
identity arrows, as well as to the first nonzero natural number; it is my hope
that no confusion results.

\chapter{Theories}
\label{ch:theories}

The first step will be to obtain a mathematical description of the notion of an
algebraic theory, of which the familiar theories of groups, rings, modules etc.
are examples.
In this chapter, we present some standard ways of doing this, and prove that
they are equivalent.
The most convenient description for our purposes will be the notion of
\emph{clone}, which appears to have been introduced by Philip Hall in
unpublished lecture notes in the 1960s, and may be found on \cite{cohn} page
132, under the name ``abstract clone''.
The treatment here follows \cite{johnstone}.
The remainder of the material in this chapter is all well-known, and may be
found in e.g. \cite{borceux} chapters 3 and 4, or \cite{a+r} chapter 3.

In the next chapter, we shall describe \emph{operads},
which allow us to capture certain algebraic theories in an especially simple
way, suitable for categorification, and we shall show how operads relate to the
clones described in this chapter.

\section{Syntactic approach}

The most traditional way of formalizing algebraic theories is syntactic.
In this approach, we abstract from the standard ``operations plus equations''
description (used to describe e.g. the theory of groups) to create
\defterm{presentations of algebraic theories}, and define a notion of an
\defterm{algebra} for a presentation.

\begin{defn}
\index{signature}
A \defterm{signature} $\Phi$ is an object of $\Set^\natural$.
\end{defn}
In other words, a signature is a sequence of sets $\Phi_0, \Phi_1, \Phi_2,
\dots$. 

Fix a countably infinite set $X = \{x_1, x_2, \dots, \}$, whose elements we
call \defterm{variables}.
\index{variable}
Throughout, let $\Phi$ be a signature.

\begin{defn}
\label{def:term}
\index{term}
Let $n \in \natural$.
An \defterm{\nary\ $\Phi$-term} is defined by the following inductive clauses:
\begin{itemize}
\item $x_1, x_2, \dots, x_n$ are \nary\ terms.
\item If $\phi \in \Phi_m$ and $t_1, \dots, t_m$ are \nary\ terms, then $\phi
\tcomp (t_1, \dots, t_m)$ is an \nary\ term.
\end{itemize}
A \defterm{$\Phi$-term} is an \nary\ $\Phi$-term for some $n \in \natural$.
\end{defn}

\begin{defn}
\index{var}
Let $t$ be an \nary\ $\Phi$-term.
Then $\var(t)$ is the sequence of elements of $\{x_1, \dots, x_n\}$ given as
follows:
\begin{itemize}
\item $\var(x_i) = (x_i)$,
\item $\var(\phi \tcomp (t_1, \dots, t_n)) =
\var(t_1) \conc \var(t_2) \conc \dots \conc \var(t_n)$,
\end{itemize}
where $\conc$ is concatenation.
\end{defn}

\begin{defn}
\index{supp}
Let $t$ be an \nary\ $\Phi$-term.
Then $\supp(t)$, the \defterm{support of $t$}, is the subset of $\{x_1, \dots,
x_n\}$ given as follows:
\begin{itemize}
\item $\supp(x_i) = \{x_i\}$,
\item $\supp(\phi \tcomp (t_1, \dots, t_n)) =
\supp(t_1) \union \supp(t_2) \union \dots \union \supp(t_n)$,
\end{itemize}
\end{defn}

\begin{defn}
\index{labelling function}
Let $t$ be an \nary\ $\Phi$-term, with $\var(t) = (x_{i_1}, \dots, x_{i_m})$.
The \defterm{labelling function} $\munge(t)$ of $t$ is the function $\m \to \n$
sending $j$ to $i_j$.
\end{defn}

\begin{defn}
\index{equation}
An \defterm{\nary\ $\Phi$-equation} is a pair $(s,t)$ of \nary\ $\Phi$-terms.
A \defterm{$\Phi$-equation} is an \nary\ $\Phi$-equation for some $n \in
\natural$.
\end{defn}

\begin{defn}
\label{def:sr_linterms}
\index{linear!term}\index{strongly regular!term}
An \nary\ term $t$ is \defterm{linear} if $\munge(t)$ is a bijection, and
\defterm{strongly regular} if $\munge(t)$ is an identity.
\index{linear!equation}\index{strongly regular!equation}
An equation $(s,t)$ is \defterm{linear} if both $s$ and $t$ are linear, and
\defterm{strongly regular} if both $s$ and $t$ are strongly regular.
\end{defn}

In other words, a term is linear if every variable is used exactly once, and
strongly regular if every variable is used exactly once in the order $x_1,
\dots, x_n$.
Up to trivial relabellings, an equation is linear if every variable is used
exactly once on both sides, though not necessarily in the same order: an
example is the commutative equation $x_1.x_2 = x_2.x_1$.
An equation is strongly regular if every variable is used exactly once in the
same order on both sides.
An example is the associative equation $x_1 . (x_2 . x_3) = (x_1.x_2).x_3$,
though some care is needed.
Strictly, a $\Phi$-equation is a pair $(n, (s,t))$ where $n \in \natural$ and
$s, t$ are \nary\ $\Phi$-terms.
The equation $(3, ((x_1.x_2).x_3, x_1.(x_2.x_3)))$ is strongly regular, but the
equation $(4, ((x_1.x_2).x_3, x_1.(x_2.x_3)))$ is not.

Classically, an \nary\ equation $(s,t)$ is \defterm{regular} \index{regular} if
$\munge(t)$ and $\munge(s)$ are surjections: however, we will not consider
regular equations further.
The term ``linear'' is borrowed from linear logic, and the term ``strongly
regular'' is due to Carboni and Johnstone (from \cite{c+j1}).

\begin{defn}
\index{presentation!of an algebraic theory}
A \defterm{presentation of a (one-sorted) algebraic theory} is 
\begin{itemize}
\item a signature $\Phi$,
\item a set $E$ of $\Phi$-equations.
\end{itemize}
\index{generating operation}
Elements of $\Phi_n$ are called (\nary) \defterm{generating operations}.
\end{defn}

\begin{defn}
\index{linear!presentation}
\index{strongly regular!presentation}
Let $P = (\Phi, E)$ be a presentation of an algebraic theory.
$P$ is \defterm{linear} if every equation in $E$ is linear, and
\defterm{strongly regular} if every equation in $E$ is strongly regular.
\end{defn}

We will return to the consideration of linear and strongly regular
presentations once we have defined operads.

\begin{defn}
\index{algebra!for a signature}
Let $\Phi$ be a signature.
An \defterm{algebra} for $\Phi$ is
\begin{itemize}
\item a set $A$,
\item for each \nary\ operation $\phi$, a map $\phi_A : A^n \to A$.
\index{primitive operation}
These are called the \defterm{primitive operations} of the algebra $A$.
\end{itemize}
\end{defn}

Let $\Phi$ be a signature, and $A$ a $\Phi$-algebra.
Each \nary\ $\Phi$-term $t$ gives rise to an \defterm{\nary\
derived operation} \index{derived operation} $t_A :A^n \to A$, defined
recursively as follows:

\begin{itemize}
\item if $t = x_i$, then $t_A$ is projection onto the $i$th factor,
\item if $t = \phi\tcomp(t_1, \dots, t_m)$, then $t_A$ is the composite
\[
\xymatrix {
A^n \ar[rr]^{((t_1)_A, \dots, (t_m)_A)} && A^m \ar[r]^{\phi_A} & A
}.
\]
\end{itemize}

Let $\term_n{\Phi}$ denote the set of \nary\ derived operations over $\Phi$.
Then $\term{\Phi}$ is a signature for every signature $\Phi$.
A morphism of signatures $f : \Phi \to \Psi$ induces a map $\bar f : \term\Phi
\to \term \Psi$.
Indeed, $\term$ is an endofunctor on $\Set^\natural$, and in Section
\ref{sec:synclass} we shall show that it is actually a monad.

\begin{defn}
\index{algebra!for a presentation}
Let $P = (\Phi, E)$ be a presentation of an algebraic theory.
A \defterm{$P$-algebra} is a $\Phi$-algebra $A$ such that, for every equation
$(s,t)$ in $E$, the derived operations $s_A, t_A$ are equal.
\end{defn}

An algebra for $\Phi$ is an algebra for $(\Phi, \{\})$.
Conversely, every algebra for $(\Phi, E)$ is an algebra for $\Phi$.

\begin{defn}
\index{morphism!of algebras for a signature}
Let $\Phi$ be a signature, and $A$ and $B$ be $\Phi$-algebras.
A \defterm{morphism of $\Phi$-algebras} $f : A \to B$ is a map $f:A \to B$
which commutes with every primitive operation:
\[
\commsquare {A^n} {f^n} {B^n} {\phi_A} {\phi_B} A f B
\]
for every $n \in \natural$ and every \nary\ primitive operation $\phi$.
If $P = (\Phi, E)$ is a presentation, then a \defterm{morphism of
$P$-algebras} is a morphism of $\Phi$-algebras.
\end{defn}

By an easy induction, a morphism of $\Phi$-algebras will commute with every
derived operation too.

Given a presentation $P$, there is a category $\Alg{P}$ \index{$\Alg{P}$} whose
objects are $P$-algebras and whose arrows are $P$-algebra morphisms.
We shall call a category \C\ a \defterm{variety of algebras} (or simply a
\defterm{variety}) if \C\ is isomorphic to $\Alg P$ for some presentation $P$.

We will need to consider \defterm{closures} of sets of equations; the
idea is that the closure of $E$ contains the members of $E$ and all of their
consequences.

\begin{defn}
\index{grafting!of terms}
Let $t$ be an \nary\ $\Phi$-term, and $t_1, \dots, t_n$ be $\Phi$-terms.
Then the \defterm{graft} $t (t_1, \dots, t_n)$ is the $\Phi$-term defined
recursively as follows.
\begin{itemize}
\item If $t = x_i$, then $t (t_1, \dots, t_n) = t_i$.
\item If $t = \phi(s_1, \dots, s_m)$, where $\phi \in \Phi_m$ and $s_1, \dots,
s_m$ are \nary\ $\Phi$-terms, then 

$t  (t_1, \dots, t_n)
	= \phi  ((s_1 (t_1, \dots, t_n)), \dots, s_m (t_1, \dots, t_n))$.
\end{itemize}
\end{defn}

\begin{defn}
\index{closure}
Let $\Phi$ be a signature and $E$ be a set of $\Phi$-equations.
The \defterm{closure} $\bar E$ of $E$ is the smallest equivalence relation on
$\term \Phi$ which contains $E$ and is closed under grafting of terms:
\begin{itemize}
\item if $(s,t) \in \bar E$, then $(s  (t_1, \dots, t_n), t  (t_1, \dots, t_n)) \in \bar E$ for all $t_1, \dots, t_n$.
\item if $(s_i, t_i) \in \bar E$ for $i = 1, \dots, n$, then $(t  (s_1,
\dots, s_n), t  (t_1, \dots, t_n)) \in \bar E$ for all $t$.
\end{itemize}
\end{defn}

\section{Clones}

Clones attempt to capture theories directly: a clone is to a presentation of an
algebraic theory as a group is to a presentation of that group.

\begin{defn}
\label{def:clone}
\index{clone}
A \defterm{clone} $K$ is
\begin{itemize}
\item a sequence of sets $K_0, K_1, \dots$,
\item for all $m,n \in \natural$, a function $\ccomp : K_n \times (K_m)^n \to
K_m$,
\item for each $n \in \natural$ and each $i \in \{1,\dots,n\}$, an
element $\delta^i_n \in K_n$
\end{itemize}
such that
\begin{itemize}
\item for each $f \in K_n$, $g_1,\dots,g_n \in K_m$, $h_1, \dots, h_{m} \in
K_p$,
\[
f \ccomp(g_1 \ccomp (h_1, \dots, h_m), \dots, g_n \ccomp (h_1, \dots,
h_m))
= (f \ccomp (g_1, \dots, g_n)) \ccomp (h_1, \dots, h_{m})
\]
\item for all $n$, all $i \in {1, \dots, n}$ and all $f_1, \dots, f_n \in K_m$,
\[
\delta^i_n \ccomp (f_1, \dots, f_n) = f_i
\]
\item for all $n$ and $f \in K_n$,
\[
f \ccomp (\delta_n^1, \dots, \delta_n^n) = f
\]
\end{itemize}
\end{defn}

\begin{example}
\index{endomorphism clone}
Let $\C$ be a finite product category, and $A$ be an object of \C.
The \defterm{endomorphism clone} of $A$, $\End(A)$, is defined as follows:
\begin{itemize}
\item $\End(A)_n = \C(A^n, A)$ for each $n \in \natural$, 
\item for all $n \in \natural$ and $i \in \{1, \dots, n\}$, the map
$\delta^i_n$ is the projection of $A^n$ onto its $i$th factor,
\item for all $n, m \in \natural$, all $f \in \End(A)_n$, and all $g_1, \dots,
g_n \in \End(A)_m$, the morphism $f\ccomp(g_1, \dots, g_n)$ is the composite
$fh$, where $h$ is the unique arrow $A^m \to A^n$ induced by the maps $g_1,
\dots, g_n$ and the universal property of $A^n$.
\end{itemize}
\end{example}

\begin{defn}
\index{morphism!of clones}
A \defterm{morphism} of clones $f: K \to K'$ is a map in $\Set^\natural$
which commutes with the composition operations and $\delta$s.
\end{defn}

\begin{defn}
\index{algebra!for a clone}
Let $K$ be a clone, and \C\ a finite product category with specified finite
powers.
An \defterm{algebra} for $K$ in \C\ is an object $A \in \C$ and a morphism of
clones $K \to \End(A)$.
\end{defn}

Equivalently, an algebra for a clone $K$ in a finite product category \C\ with
specified powers is
\begin{itemize}
\item an object $A$ of \C,
\item for each $n \in \natural$ and each $k \in K_n$, a morphism $\hat k : A^n
\to A$
\end{itemize}
such that
\begin{itemize}
\item for all $n \in \natural$ and all $i \in \{1, \dots n\}$, the morphism
$\widehat \delta^i_n$ is the projection of $A^n$ onto its $i$th factor;
\item for all $n, m \in \natural$, all $f \in K_n$, and all $g_1, \dots, g_n
\in K_m$, the diagram
\[
\xymatrix{
	& A^m \ar[ddl]_{\widehat{g_1}} \unar[ddd]^h \ar[ddr]^{\widehat{g_n}}
		\ar@/r3cm/[dddd]^{\widehat{f\ccomp(g_1, \dots, g_n)}} \\ \\
	A & & A \\
	& A^n \ar[ul]^{\widehat \delta^1_n} \ar[ur]_{\widehat \delta^n_n}
		\ar[d]^f \\
	& A
}
\]
commutes, where $h$ is the unique arrow induced by the universal property of
$A^n$.
\end{itemize}

\begin{defn}
\index{morphism!of algebras for a clone}
Let $A$ and $B$ be algebras for a clone $K$ in a finite product category $\C$
with specified finite powers.
A \defterm{morphism} of algebras $A \to B$ is a morphism $F : A \to B$ in \C\
such that the diagram
\[
\commsquare{A^n}{F^n}{B^n}{\hat k}{\hat k} A F B
\]
commutes for all $n \in \natural$ and all $k \in K_n$.
\end{defn}

Algebras for a clone and their morphisms form a category: we call this category
$\Algfctr_\C(K)$, or $\Alg{K}$ in the case where $\C = \Set$.\index{\Algfctr}

Clones can be enriched in any finite product category $\V$ in an obvious
way: the sequence of sets $K_0, K_1, \dots$ becomes a sequence of objects of
$\V$, and so on.

\section{Lawvere theories}

Lawvere theories are a particularly elegant approach to describing algebraic
theories, introduced by Lawvere in his thesis \cite{lawthesis}.
\index{Lawvere, F. William}
Like a clone, a Lawvere theory (sometimes called a \defterm{finite product
theory}) is an object that represents the semantics of the theory directly; in
Lawvere theories, the data are encoded into a category.
Algebras for the theory are then certain functors from the Lawvere theory to
\Set.

\begin{defn}
\index{Lawvere theory}
A \defterm{Lawvere theory} is a category $\T$ whose objects form a
denumerable set $\{\bf{0}, \bf{1}, \bf{2}, \dots \}$, such that $\bf{n}$ is the
$n$-th power of $\bf{1}$.
\index{morphism!of Lawvere theories}
A \defterm{morphism} of Lawvere theories $\T \to \S$ is an identity-on-objects
functor $\T \to \S$ which preserves projection maps.
The category of Lawvere theories and their morphisms is called $\Law$.
\index{\Law}
An \defterm{algebra} for $\T$ is a functor $F:\T \to \Set$ which preserves
finite products.
\index{algebra!for a Lawvere theory}
A \defterm{morphism} of algebras is a natural transformation.
\index{morphism!of algebras for a Lawvere theory}
The \defterm{category of $\T$-algebras} is the full subcategory of $[\T, \Set]$
whose objects are finite-product-preserving functors.
\end{defn}

Lawvere theories encode algebraic theories by storing the \nary\ operations
of the theory as morphisms $\bf{n} \to \bf{1}$.

We can consider algebras for Lawvere theories in categories other than $\Set$:
an algebra for a Lawvere theory $\T$ in a finite product category $\C$ is just
a finite-product-preserving functor $\T \to \C$.
This captures our usual notions of, for instance, topological groups: a
topological group is just an algebra for the Lawvere theory of groups in the
category $\Top$.
Much the same could be said for clones and presentations, of course, but in this
case the definition is especially economical.

We may generalize this definition as follows:

\begin{defn}
\index{finite product theory, multi-sorted}
Let $S$ be a set.
An \defterm{$S$-sorted finite product theory} is a small finite product
category whose underlying monoidal category is strict and whose monoid of
objects is the free monoid on $S$.
Elements of $S$ will be called \defterm{sorts}.
Algebras and morphisms of algebras are defined as above.
\end{defn}

\section{Finitary monads}

\index{monad}
Recall that a \defterm{monad} on a category $\C$ is a monoid object in the
category $[\C, \C]$ of endofunctors on $\C$.
Concretely, a monad is a triple $(T, \mu, \eta)$ where
\begin{itemize}
\item $T: \C \to \C$ is a functor,
\item $\mu : T^2 \to T$ is a natural transformation,
\item $\eta : 1_\C \to T$ is a natural transformation,
\end{itemize}
and $\mu, \eta$ satisfy coherence axioms which are analogues of the usual
associativity and unit laws for monoids, namely
\begin{equation}
\label{eqn:monad_assoc}
\xymatrix{
	& T^3 \ar[dl]_{\mu T} \ar[dr]^{T\mu} \\
	T^2 \ar[dr]_\mu & & T^2 \ar[dl]^\mu \\
	& T
}
\end{equation}
\begin{equation}
\label{eqn:monad_unit}
\xymatrix{
	T \ar[r]^{T\eta} \ar[dr]_{1_T}
	& T^2 \ar[d]^\mu
	& T \ar[l]_{\eta T} \ar[dl]^{1_T} \\
	& T
}
\end{equation}
We shall often abuse notation and refer to the monad $(T, \mu, \eta)$ as simply
$T$.

\begin{defn}
\label{def:monad_morphism}
\index{morphism!of monads}
Let $(T_1, \mu_1, \eta_1), (T_2, \mu_2, \eta_2)$ be monads on a category $\C$.
A \defterm{morphism of monads} $(T_1, \mu_1, \eta_1) \to (T_2, \mu_2, \eta_2)$
is a natural transformation $\alpha : T_1 \to T_2$ such that the diagrams
\begin{equation}
\label{eqn:monad_morph_mu}
\xymatrix{
	T_1^2 \ar[r]^{\mu_1} \ar[d]_{\alpha * \alpha} & T_1 \ar[d]^\alpha \\
	T_2^2 \ar[r]_{\mu_2} & T_2
}
\end{equation}
\begin{equation}
\label{eqn:monad_morph_eta}
\xymatrix{
	1 \ar[r]^{\eta_1} \ar[dr]_{\eta_2} & T_1 \ar[d]^{\alpha} \\
	& T_2
}
\end{equation}
commute.
\end{defn}

Monads on \C\ and monad morphisms form a category $\Mnd(\C)$. \index{$\Mnd(\C)$}
This notion (or rather, a 2-categorical version) was introduced and studied by
Street in \cite{street}.  \index{Street, Ross}

\begin{defn}
\index{filtered!category}
A category $\C$ is \defterm{filtered} if every finite diagram in $\C$ admits a
cocone.
\end{defn}

Equivalently, $\C$ is filtered if:
\begin{itemize}
\item $\C$ is nonempty;
\item for every pair of parallel arrows $\parallelpair A f g B$ in $\C$, there
is an arrow $h: B \to C$ such that $hf = hg$;
\item for every pair of objects $A, B$, there is an object $C$ and arrows
\[
\xymatrix{
A \ar[dr]^f \\
& C \\
B \ar[ur]_g
}
\]
\end{itemize}

Filteredness generalizes the notion of directedness for posets (a
\defterm{directed poset} is a poset in which every finite subset has an upper
bound). \index{directed}
A filtered category which is also a poset is precisely a directed poset.

\begin{defn}
\index{filtered!colimit}
A \defterm{filtered colimit} in a category $\C$ is the colimit of a diagram $D: \I \to \C$, where $\I$ is a filtered category.
\end{defn}

\begin{theorem}
Every object in \Set\ is a filtered colimit of finite sets.
\end{theorem}
\begin{proof}
Let $X \in \Set$, and consider the subcategory $\I$ of $\Set$ whose objects are
finite subsets of $X$ and whose morphisms are inclusions.
This is a directed poset, and thus a filtered category.
$X$ is the colimit of the inclusion of $\I$ into $\Set$.
\end{proof}

\begin{theorem}
Let $\I$ be a small category.
Colimits of shape $\I$ in $\Set$ commute with all finite limits iff $\I$ is
filtered.
\end{theorem}
\begin{proof}
See \cite{m+m}, Corollary VII.6.5.
\end{proof}

\begin{defn}
\index{finitary!functor}
A functor $F: \C \to \D$ is \defterm{finitary} if it preserves filtered
colimits.
\end{defn}

\begin{defn}
\index{finitary!monad}
A monad $(T, \mu, \eta)$ on $\C$ is \defterm{finitary} if $T$ is finitary.
\end{defn}

A finitary monad on \Set\ is determined by its behaviour on finite sets, in the
following sense: since every set $X$ is a filtered colimit of its finite
subsets, then $TX$ has to be the colimit of the images under $T$ of the finite
subsets of $X$.

\section{Equivalences}

Let $(\Phi, E)$ be a presentation of an algebraic theory.
We define $K_{(\Phi, E)}$ to be the clone whose operations are elements of the
quotient signature $(\term\Phi)/\bar E$, with composition given by grafting,
and $\delta^i_n = x_i$ for all $i, n \in \natural$.
By definition of $\bar E$, grafting gives a well-defined family of composition
functions on $K_{(\Phi, E)}$.
Conversely, given a clone $K$, we may define a presentation of an algebraic
theory $(\Phi_K, E_K)$, by taking $(\Phi_K)_n = K_n$ for all $n \in \natural$,
and for all $n, m \in \natural$,
all $k \in K_n$ and all $k_1, \dots, k_n \in K_m$, letting $E_m$ contain the
equation
\[
k(k_1(x_1, \dots, x_m), \dots, k_n(x_1, \dots, x_m))
	= k\ccomp(k_1,\dots, k_n)(x_1, \dots, x_m).
\]

\begin{lemma}
Let $K$ be a clone. 
Then $K_{(\Phi_K, E_K)}$ is isomorphic to $K$.
\end{lemma}
\begin{proof}
See \cite{johnstone}, Lemma 1.7.
\end{proof}

\begin{lemma}
Let $(\Phi, E)$ be a presentation of an algebraic theory.
Let $(\Phi', E')$ be the presentation obtained from the clone $K_{(\Phi, E)}$.
Then the category $\Alg{\Phi, E}$ is isomorphic to the category $\Alg{\Phi',
E'}$
\end{lemma}
\begin{proof}
See \cite{johnstone}, Lemma 1.8.
\end{proof}

\begin{defn}
\index{linear!clone}
\index{strongly regular!clone}
Let $K$ be a clone.
We say that $K$ is \defterm{strongly regular} (resp. \defterm{linear}) if
there exists a strongly regular (resp. linear) presentation $P$ such that
$K = K_{(\Phi, E)}$.
\end{defn}

Given a clone $K$, we construct a Lawvere theory $\T_K$ for which
$\T_K({\bf n},{\bf m}) = (K_n)^m$.
Suppose $f = (f_1, \dots, f_m) \in \T_K(\bf{n}, \bf{m})$ and $g = (g_1, \dots,
g_p) \in \T_K(\bf{m}, \bf{p})$, then the composite $gf$ is $(g_1\ccomp (f_1, \dots, f_m), \dots, g_p \ccomp (f_1, \dots, f_m))$.
By the axioms for a clone, this is a category, with the identity map on
$\bf{n}$ being $(\delta^1_n, \dots, \delta^n_n)$. 
It remains to show that ${\bf n}$ is the $n$th power of $\bf 1$ for every $n
\in \natural$. 
The $i$th projection of $\bf n$ onto $1$ is evidently $\delta^i_n$: we must
show that these have the requisite universal property.
Take $m, n \in \natural$, and $n$ maps $f_1, \dots, f_n : {\bf m} \to {\bf 1}$
in $\T_K$.
The diagram
\[
\xymatrixcolsep{3pc}
\xymatrix{
	& {\bf m} \ar@/u0.8cm/[ddl]_{f_1} \unar[d]^h \ar@/u0.8cm/[ddr]^{f_n} \\
	& {\bf n} \ar[dl]^{\delta^1_n} \ar[dr]_{\delta^n_n} \\
	{\bf 1} & \dots & {\bf 1}
}
\]
commutes if and only if $h = (f_1, \dots, f_n)$, and hence ${\bf n}$ is indeed
the $n$th power of ${\bf 1}$, and so $\T_K$ is a Lawvere theory.

The Lawvere theories so constructed evidently respect isomorphisms of clones.
Furthermore, the diagram
\[
\xymatrixrowsep{3pc}
\xymatrix{
\Clone \ar[rr]^{\T_{(-)}} \ar[dr]_\Algfctr & & \Law \ar[dl]^\Algfctr \\
& \CAT\op
}
\]
commutes up to equivalence:

\begin{theorem}
Let $K$ be a clone.
Then $\Alg K \simeq \Alg{\T_K}$.
\end{theorem}
\begin{proof}
Let $A$ be a $K$-algebra.
We define a $\T_K$-algebra $F_A$ as follows:
\begin{itemize}
\item $F_A {\bf n} = A^n$ for all $n \in \natural$;
\item If $k \in \T_K({\bf n}, 1) = K_n$, then $F_Ak = \hat k$;
\item if $(k_1, \dots, k_n) : {\bf m} \to {\bf n}$ in $\T_K$, then 
$F_A (k_1, \dots, k_n)$ is the unique arrow $A^m \to A^n$ such that the diagram
\[
\xymatrix{
	& A^m \ar[ddl]_{\widehat k_1} \ar[ddr]^{\widehat k_n} \unar[d] \\
	& A^n \ar[dl]^{\widehat{\delta^1_n}} \ar[dr]_{\widehat{\delta_n^n}} \\
	A & \dots & A
}
\]
commutes.
\end{itemize}
Let $f : A \to B$ be a morphism of $K$-algebras.
Then the diagram
\[
\commsquare{A^n}{f^n}{B^n}{\widehat k}{\widehat k}A f B
\]
commutes for all $n \in \natural$ and all $k \in K_n$.
By the universal property of $B^m$, the diagram
\[
\commsquare{A^n}{f^n}{B^n}{F_A(k_1, \dots, k_m)}{F_B(k_1, \dots, k_m)}
	{A^m}{f^m}{B^m}
\]
commutes for all $n, m \in \natural$ and all $(k_1, \dots, k_m) : {\bf n} \to {\bf m}$ in $\T_K$.
Hence $F_f = (f^n)_{{\bf n}}$ is a natural transformation $F_A \to F_B$, and
hence a morphism of $\T_K$-algebras.
This defines a functor $F_{(-)} : \Alg{K} \to \Alg{\T_K}$; we wish to show that
it is an equivalence.

For every $\T_K$-algebra $G$, we may define a $K$-algebra $A$ by setting $A =
G{\bf 1}$ and $\hat k = G({\bf n} \stackrel{k}{\longrightarrow} {\bf 1})$ for
all $k \in K_n$ and all $n \in \natural$.
Then $F_A$ is isomorphic as a $\T_K$-algebra to $G$, and hence the functor
$F_{(-)} : \Alg{K} \to \Alg{\T_K}$ is essentially surjective on objects.
We shall show further that it is full and faithful.
Let $A$ and $B$ be $K$-algebras, and let $\alpha_n : F_A \to F_B$ be
a morphism between their associated $\T_K$-algebras.
Since the diagram
\[
\commsquare{A^n}{\alpha_n}{B^n}{\widehat{\delta^i_n}}{\widehat{\delta^i_n}}
	A {\alpha_1} B
\]
commutes for all $n \in \natural$ and all $i \in \{1, \dots, n\}$, it must be
the case that $\alpha_n = \alpha_1^n$ for all $n$.
Hence, the diagram
\[
\commsquare{A^n}{\alpha_1^n}{B^n}{\widehat k}{\widehat k}A {\alpha_1} B
\]
must commute for all $n \in \natural$ and all $k \in K_n$.
So $\alpha_1$ is a $K$-algebra morphism, and $\alpha_n = F_{\alpha_1}$.
Hence $F_{(-)}$ is full.
Suppose $F_f = F_g$; then $(F_f)_{\bf 1} = (F_g)_{\bf 1}$, so $f = g$.
Hence $F_{(-)}$ is faithful; and hence it is an equivalence of categories.
\end{proof}

Given a Lawvere theory $\T$, we can construct a clone $K_\T$, as follows:
\begin{itemize}
\item Let $(K_\T)_n = \T(\bf{n},\bf{1})$ for all $n \in \natural$.
\item For all $n, m  \in \natural$, all $f \in (K_\T)_n$ and all $g_1, \dots,
g_n \in (K_\T)_m$, let $f \ccomp (g_1, \dots, g_n) = f \fcomp (g_1 + \dots
+ g_n) \fcomp \Delta$,
where $(g_1 + \dots + g_n)$ is the unique map ${\bf mn} \to
{\bf n}$ in $\T$ such that the diagram
\[
\xymatrixcolsep{3pc}
\xymatrix{
& {\bf mn} \unar[dd]|{g_1 + \dots + g_n}
	\ar[dl] \ar[dr] \\
{\bf m} \ar[dd]^{g_1}
& & {\bf m} \ar[dd]_{g_n} \\
& {\bf n} \ar[dl] \ar[dr] \\
{\bf 1} & \dots & {\bf 1}
}
\]
commutes, and $\Delta : {\bf m} \to {\bf mn}$ is the diagonal map (or
equivalently, the image of the codiagonal function ${\bf mn} \to {\bf m}$ under
the contravariant embedding of $\F$ into $\T$).

\item For all $n \in \natural$ and all $i \in {1, \dots, n}$, let $\delta_n^i$
be the $i$th projection $\bf{n} \to \bf{1}$.
\end{itemize}

This extends to a functor $K_{(-)} : \Law \to \Clone$, as follows: given
Lawvere theories $\T_1$ and $\T_2$, and a morphism of Lawvere theories $F: \T_1
\to \T_2$, let $K_F$ be the map of signatures sending $k \in (K_{\T_1})_n =
\T_1(\bf{n}, \bf{1})$ to $Fk \in (K_{\T_2})_n = \T_2(\bf{n}, \bf{1})$.
Since $F$ is a functor, and thus commutes with composition in $\T_1, \T_2$,
then $K_F$ must commute with composition in $K_{\T_1}$ and $K_{\T_2}$.
Since $F$ preserves finite products, it commutes with the projection maps in
$K_{\T_1}$ and $K_{\T_2}$.
Thus, $K_F$ is a morphism of clones.

\begin{theorem}
The functor $K_{(-)}$ is pseudo-inverse to the functor $\T_{(-)}$.
\end{theorem}
\begin{proof}
Since every object in a Lawvere theory is a copower of ${\bf 1}$, a Lawvere
theory $\T$ is entirely determined (up to isomorphism) by the hom-sets $\T({\bf
n}, {\bf 1})$, and thus by $K_\T$.
The theorem follows straightforwardly.
\end{proof}

Given a Lawvere theory $\T$, we construct a monad $(T, \mu, \eta)$ on \Set\ as
follows:
\begin{itemize}
\item If $X$ is a set, let $T X  = \int^{n\in \F} \T({\bf n}, {\bf 1})
\times X^n$.
\item If $x \in X$, then $\eta(x) = (1,x) \in T X$.
\item If $f : {\bf n} \to {\bf 1}$ in $T$ and $(f_i, x^i\seq) \in \T({\bf k_i},
{\bf 1}) \times X^{k_i}$ for $i = 1, \dots, n$, then
\[
\mu(f,((f_1,x^1\seq), \dots, (f_n, x^n\seq)))
	= (f\fcomp(f_1+ \dots+ f_n), x\dseq)
\]
\end{itemize}

\begin{theorem}
The monad so constructed is finitary.
\end{theorem}
\begin{proof}
See \cite{a+r}, Theorems 3.18 and 1.5, and Remarks 3.4(4) and 3.6(6).
\end{proof}

Given a finitary monad $T$ on \Set, we can construct a Lawvere theory $\T$.
Take the full subcategory $\F_T$ of the Kleisli category $\Set_T$ whose objects
are finite sets.
Now let $\T$ be the skeleton of the dual of $\F_T$.
The monad induced by this Lawvere theory is isomorphic to the original monad:
see \cite{a+r}, Remark 3.17 and Theorem 3.18.

The moral of the above theorems is that presentations, clones, Lawvere theories
and finitary monads on \Set\ all capture the same notion, and may be used
interchangeably.
Further, the notion that is captured corresponds to our usual intuitive
understanding of equational algebraic theories.

The equivalence between (finitary monads on \C) and (monads on \C\ that may be
described by a finitary presentation) may actually be generalized to the
case where \C\ is an arbitrary finitely presentable category: see \cite{k+p}.

\chapter{Operads}
\label{ch:operads}

Operads arose in the study of homotopy theory with the work of Boardman and
Vogt \cite{b+v}, and May \cite{may}.
In that field they are an invaluable tool: \cite{mss} describes a diverse range
of applications.
Independently, multicategories (which are to operads as categories are to
monoids) had arisen in categorical logic with the work of Lambek
\cite{lambek}.
Multicategories are sometimes called ``coloured operads''.

We will use multicategories and operads as tools to approach universal algebra:
while operads are not as expressive as Lawvere theories, they can be easily
extended to be so, and the theories that \emph{can} be represented by operads
provide a useful ``toy problem'' to help us get started.

\begin{figure}
\centerline{
	\epsfxsize=3in
	\epsfbox{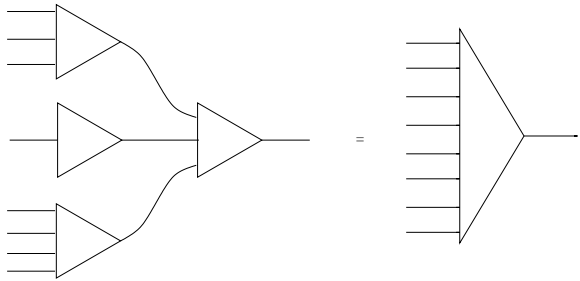}
}
\caption{Composition in a multicategory}
\label{fig:multicat}
\end{figure}

Informally, categories have objects and arrows, where an arrow has one source
and one target; multicategories have objects and arrows with one target but
multiple sources (see Fig. \ref{fig:multicat}); and operads are one-object
multicategories.
Multicategories (and thus operads) have a composition operation that is
associative and unital.

\section{Plain operads}
\label{sec:plain_operads}

\begin{defn}
\index{multicategory!plain}
A \defterm{plain multicategory} (or simply ``multicategory'') $\C$
consists of the following:
\begin{itemize}
\item a collection $\C_0$ of \emph{objects},
\item for all $n \in \natural$ and all $c_1, \dots, c_n,
d \in \C_0$, a set of \emph{arrows} $\C(c_1, \dots, c_n; d)$,
\item for all $n, k_1, \dots k_n \in \natural$ and $c_1^1 \dots, c_{k_n}^n,
d_1, \dots, d_n, e \in \C_0$, a function called \emph{composition}
\[\ocomp
: \C(d_1, \dots, d_n;e)  \times
\C(c^1_1, \dots, c^1_{k_1};d_1) \times \dots \times
\C(c^n_1, \dots, c^n_{k_n};d_n)
\to
\C(c^1_1, \dots, c^n_{k_n};e)\]
\item for all $c \in \C$, an \emph{identity arrow} $1_c \in \C(c;c)$
\end{itemize}
satisfying the following axioms:
\begin{itemize}
\item \emph{Associativity:} $f \ocomp (g\seq \ocomp h\dseq) = (f \ocomp g\seq)
\ocomp h\dseq$ wherever this makes sense (we borrow the ${}\seq$ notation for
sequences from chain complexes)
\item \emph{Units:} $1 \ocomp f = f = f \ocomp (1, \dots, 1)$ for all $f$.
\end{itemize}
A plain multicategory \C\ is \defterm{small} if $\C_0$ forms a set.
In line with the definition above, we shall take all our multicategories to be
locally small: this restriction is not essential.
\end{defn}

We say that an arrow in $\C(c_1, \dots, c_n; d)$ is \defterm{\nary}, or has
\defterm{arity $n$}.
We remark that taking $n = 0$ gives us \emph{nullary} arrows.
This is in contrast to the definition used by some authors, who do not allow
nullary arrows.

\begin{defn}
\index{morphism!of plain multicategories}
A \defterm{morphism} of multicategories $F: \C \to \D$ is a map $F : \C_0
\to \D_0$ together with maps $F : \C(c_1, \dots, c_n ; c) \to \D(F c_1,
\dots, F c_n ; F c)$ which commute with $\ocomp$ and identities.
\index{transformation!between morphisms of plain multicategories}
A \defterm{transformation} of multicategory maps $\alpha : F \to G$ is a family
of arrows $\alpha_c \in \D(F c ; G c)$, one for each $c \in \C$, satisfying the
analogue of the usual naturality squares: for all maps $f: c_1, \dots, c_k \to
c$ in \C, we must have
\[
\alpha_c \ocomp Ff = Gf \ocomp (\alpha_{c_1}, \dots, \alpha_{c_k})
\]
\end{defn}

One is tempted to write this last condition as
\[
\xymatrix{
Fc_1, \dots, Fc_k \ar[rr]^{Ff} \ar[d]_{\alpha_{c_1}, \dots, \alpha_{c_k}}
&& Fc \ar[d]^{\alpha_c} \\
Gc_1, \dots, Gc_k \ar[rr]_{Gf}
&& Gc
}
\]
but care must be taken: in a general multicategory, $\alpha_{c_1}, \dots,
\alpha_{c_k}$ does not correspond to any single map, as it would in a monoidal
category.

Small plain multicategories, their morphisms and their transformations form a
2-category: we shall use the notation \Multicat\ for both this 2-category and
its underlying 1-category. \index{\Multicat}

To simplify the presentation of our first example, we recall the notion of
unbiased monoidal category from \cite{hohc} section 3.1:

\begin{defn}
\label{def:umoncat}
\index{unbiased!monoidal category}
An \defterm{unbiased weak monoidal category} $(\C, \otimes, \gamma, \iota)$
consists of
\begin{itemize}
\item a category $\C$,
\item for each $n \in \natural$, a functor $\otimes_n : \C^n \to \C$ called \defterm{$n$-fold tensor} and written
\[
(a_1, \dots, a_n) \mapsto (a_1 \otimes \dots \otimes a_n)
\]
\item for each $n, k_1, \dots, k_n \in \natural$, a natural isomorphism
\[
\gamma : \otimes_n \ocomp (\otimes_{k_1} \times \dots \times \otimes_{k_n})
\longrightarrow \otimes_{\sum k_i}
\]
\item a natural isomorphism
\[
\iota: 1_A \to \otimes_1
\]
\end{itemize}
satisfying
\begin{itemize}
\item associativity: for any triple sequence $a\tseq$ of objects in $\C$, the
diagram
\[
\xymatrixrowsep{4pc}
\xymatrixcolsep{-10pc}
\xymatrix{
&(
((\otimes a_{1 1}^\bullet) \otimes \dots \otimes (\otimes a_{1 k_1}^\bullet))
\otimes \dots \otimes
((\otimes a_{n 1}^\bullet) \otimes \dots \otimes (\otimes a_{n k_n}^\bullet)))
\ar[dl] \ar[dr] \\
((\otimes a_{1 \bullet}^\bullet)
\otimes \dots \otimes
(\otimes a_{n \bullet}^\bullet))
\ar[dr]
& & ((\otimes a_{1 1}^\bullet) \otimes \dots \otimes
(\otimes a_{n k_n}^\bullet))
\ar[dl] \\
& (a_{1 1}^1 \otimes \dots \otimes a_{n k_n}^{m_{k_n}})
}
\]
commutes.
\item identity: for any $n \in \natural$ and any sequence $a_1, \dots, a_n$ of
objects in $\C$, the diagrams
\[
\xymatrixcolsep{4pc}
\xymatrixrowsep{4pc}
\xymatrix{
(a_1 \otimes \dots \otimes a_n) \ar[r]^{(\iota\otimes\dots\otimes\iota)}
\ar[dr]_1
& ((a_1) \otimes \dots \otimes (a_n)) \ar[d]^\gamma \\
& (a_1 \otimes \dots \otimes a_n)
}
\]
\[
\xymatrixcolsep{4pc}
\xymatrixrowsep{4pc}
\xymatrix{
(a_1 \otimes \dots \otimes a_n) \ar[r]^{\iota}
\ar[dr]_1
& ((a_1 \otimes \dots \otimes a_n)) \ar[d]^\gamma \\
& (a_1 \otimes \dots \otimes a_n)
}
\]
commute.
\end{itemize}

\end{defn}

\begin{example}
\index{multicategory!underlying multicategory of a monoidal category}
\index{unbiased!monoidal category!underlying multicategory}
Let $\C$ be a locally small unbiased weak monoidal category.
The \defterm{underlying multicategory} $\C'$ of $\C$ has
\begin{itemize}
\item objects: objects of $\C$;
\item arrows: $\C'(a_1, \dots, a_n; b) = \C(a_1\otimes\dots\otimes a_n, b)$;
\item composition given as follows: if $f_i \in \C'(a_1^i, \dots, a_{k_i}^i ;
b_i)$ for $i = 1, \dots, n$ and $g \in \C'(b_1, \dots, b_n; c)$, then we define
$g \ocomp (f_1, \dots, f_n)$ as
\[
\xymatrixnocompile{
\bigotimes_{i,j} a^i_j \ar[dd]_{f \ocomp (g_1, \dots, g_n)}
\ar[rr]^{\gamma \otimes \dots \otimes \gamma}
&& \bigotimes_i ( \bigotimes_j a^i_j ) \ar[dd]^{\bigotimes_i {f_i}} \\
\\
c
&& \bigotimes_i b_i \ar[ll]^g
}
\]
\end{itemize}
\end{example}

\begin{defn}
\index{algebra!for a plain multicategory}
Let $M$ and $\C$ be plain multicategories.
An \defterm{algebra} for $M$ in \C\ is a morphism of multicategories $M \to \C$.
\end{defn}

\begin{defn}
Let $M$ be a plain multicategory, and \C\ be an unbiased monoidal category.
An \defterm{algebra} for $M$ in \C\ is a morphism of multicategories from $M$
to the underlying multicategory of $\C$.
\end{defn}

A \defterm{plain operad} (or simply ``operad'') is now a one-object
multicategory. \index{operad!plain}
Morphisms and transformations of operads are defined as for general
multicategories.
\index{morphism!of plain operads}
\index{transformation!between morphisms of plain operads}
As before, we use the notation \Operad\ \index{\Operad} for both the 2-category
of operads, morphisms and transformations, and its underlying 1-category.
Operads are to multicategories as monoids are to categories: just as
with monoids, this allows us to present the theory of operads in a
simplified way.

\begin{lemma}
\index{operad!plain!concrete description}
\label{lem:operad_description}
An operad $P$ can be given by the following data:
\begin{itemize}
\item A sequence $P_0, P_1, \dots$ of sets 
\item For all $n, k_1, \dots, k_n \in \natural$, a function $\ocomp: P_n
\times \prodkn{P} \to P_{\sum k_i}$ 
\item An \defterm{identity element} $1 \in P_1$
\end{itemize}
satisfying the following axioms:
\begin{itemize}
\item \emph{Associativity:} $f \ocomp (g\seq \ocomp h\dseq) = (f \ocomp g\seq)
\ocomp h\dseq$ wherever this makes sense
\item \emph{Units:} $1 \ocomp f = f = f \ocomp (1, \dots, 1)$ for all f.
\end{itemize}
\end{lemma}
\begin{proof}
Using the symbol $*$ for the unique object, let $P_n = P(*, \dots, * ; *)$,
where the input is repeated $n$ times.
The rest of the conditions follow trivially from the definition of
a multicategory.
\end{proof}

\begin{lemma}
Let $P$ and $Q$ be operads.
A morphism $f: P \to Q$ consists of a function $f_n : P_n \to Q_n$ for each $n
\in \natural$ such that, for all $n, k_1, \dots, k_n$, the diagram
\[
\xymatrix{
P_n \times \prodkn P \ar[d]_{f_n \times \prodkn f} \ar[rr]^\ocomp
&& P_{\sum k_i} \ar[d]^{f_{\sum k_i}} \\
Q_n \times \prodkn Q \ar[rr]_\ocomp
&& Q_{\sum k_i}
}
\]
commutes, and that $f_1$ preserves the identity object.

If $f$ and $g$ are morphisms of operads from $P$ to $Q$, then a transformation
from $f$ to $g$ is an element $\alpha \in Q_1$ such that $\alpha \ocomp F p = G p \ocomp (\alpha, \dots, \alpha)$ for all $n \in \natural$ and all $p \in P_n$.
\end{lemma}
\begin{proof}
Trivial.
\end{proof}

\begin{defn}
If a morphism of operads $f: P \to Q$ is such that $f_n$ has some property $X$
for all $n \in \natural$, we say that $f$ is \defterm{levelwise} $X$.
\end{defn}

\begin{example}
\label{ex:plEnd}
\index{endomorphism operad!plain}
Let $A$ be an object of a multicategory $\C$.
The \defterm{endomorphism operad} of $A$ is the full sub-multicategory
$\End(A)$ of $\C$ whose only object is $A$.
In terms of the description in Lemma \ref{lem:operad_description}, $\End(A)_n$
is the set of \nary\ arrows from $A, \dots, A$ to $A$.
Composition is as in $\C$.

In particular, if $\C$ is the underlying multicategory of some monoidal
category $\C'$, then $\End(A)_n = \C'(A \otimes \dots \otimes A, A)$.
This is the case we shall use most frequently.
\end{example}

\begin{example}
\label{ex:symmopd}
\index{\calS!as a plain operad}
\index{operad of symmetries|see{\calS}}
There is an operad $\calS$ for which each $\calS_n$ is the symmetric group
$S_n$.
Operadic composition is given as follows: if $\sigma \in S_n$, and $\tau_i \in
S_{k_i}$ for $i = 1, \dots, n$, then
\[
\sigma\ocomp(\tau_1, \dots, \tau_n) : \sum_{i = 1}^j k_i + m \mapsto 
\sum_{i: \sigma(i) < \sigma(j+1)} k_i + \tau_{j+1}(m)
\]
for all $j \in \{1, \dots, n\}$ and $m \in \{0, \dots, k_{j+1}-1\}$.
Informally, the inputs are divided into ``blocks'' of length $k_1, k_2, \dots,
k_n$, which are then permuted by $\sigma$: the elements of each block are then
permuted by the appropriate $\tau_i$.
For an example, see Figure \ref{fig:symm_comp}.
\begin{figure}[h]
\centerline{
	\epsfxsize=5in
	\epsfbox{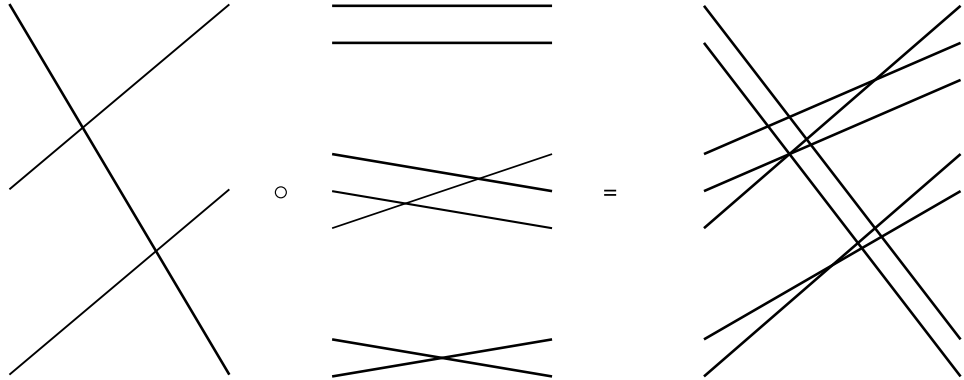}
}
\caption{Composition in the operad $\calS$ of symmetries}
\label{fig:symm_comp}
\end{figure}
\end{example}

\begin{example}
\index{\B}
\index{operad of braids|see{\B}}
There is an operad $\B$ for which each $\B_n$ is the Artin braid group
$B_n$.
Composition is analogous to that for $\calS$: the inputs are divided into
blocks, which are braided, and then the elements of the blocks are braided.
\end{example}

\begin{example}
\index{little $m$-discs operad}
\index{operad!little $m$-discs|see{little $m$-discs operad}}
Fix an $m \in \natural$.
There is an operad $LD$ for which each $LD_n$ is an embedding of $n$ copies of
the closed unit disc $D_m$ into $D_m$.
Composition is by gluing -- see Figure \ref{fig:little_discs}.
\begin{figure}[h]
\centerline{
	\epsfxsize=5in
	\epsfbox{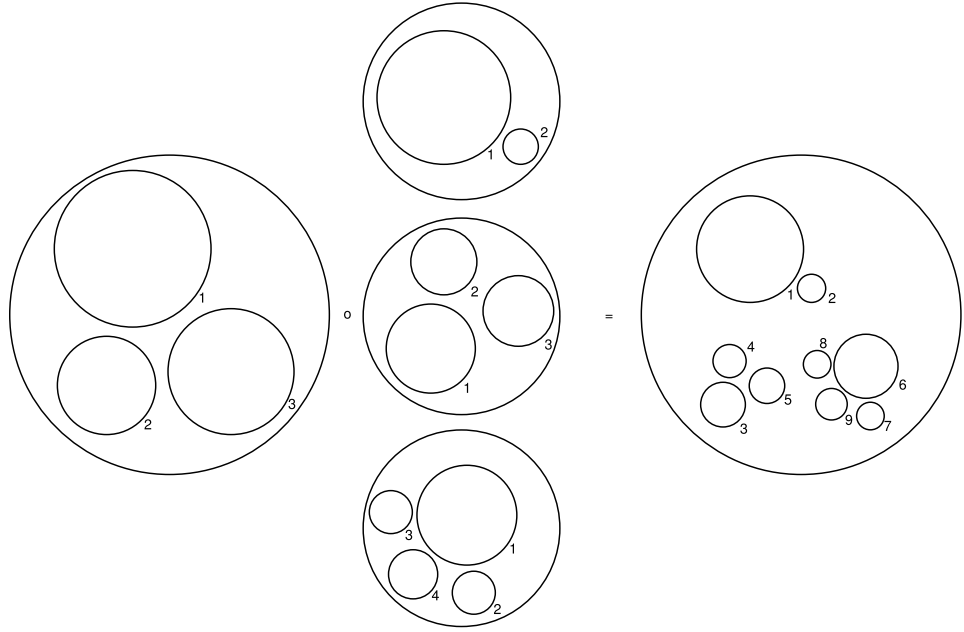}
}
\caption{Composition in the little 2-discs operad}
\label{fig:little_discs}
\end{figure}

$LD$ is known as the \defterm{little $m$-discs operad}.
\end{example}

Since we wish to use operads to represent theories, we need to have some way
of describing the models of those theories.

\begin{defn}
\index{algebra!for a plain operad}
Let $P$ be an operad.
An \defterm{algebra} for $P$ in a multicategory $\C$ is an object $A \in \C$
and a morphism of operads $(\hatmap) : P \to \End(A)$.
\end{defn}
Where $\C$ is a monoidal category, this is equivalent to requiring an object
$A \in \C$, and for each $p \in P_n$ a morphism $\hat p : A^{\otimes n} \to A$
such that $\hat 1 = 1_A$ and $\hat p \ocomp (\hat q_1 \otimes \dots \otimes
\hat q_n) = \widehat{p \ocomp (q_1 \otimes \dots \otimes q_n)}$ for all $p, q_1,
\dots, q_n \in P$.
A third equivalent definition is, for each $n \in \natural$, a map $h_n: P_n
\otimes A^{\otimes n}
\to A$, such that $h_n(p, h_n(q\seq, -)) = h_{\sum k_i}(p \ocomp q\seq,
-)$ for all $p \in P_n$, $q_i \in P_{k_i}$, and $h_1(1, -) = 1_A$.
We leave the proofs of these equivalences as an easy exercise for the reader,
and will make use of whichever formulation is most convenient at the time.

\begin{defn}
\index{morphism!of algebras for a plain operad}
Let $P$ be a plain operad, and $(A, (\hatmap))$ and $(B, (\checkmap))$ be
algebras for $P$ in a multicategory $\C$.
A \defterm{morphism} of algebras is an arrow $F : A \to B$ in $\C$ such that,
for all $n \in \natural$, the diagram
\[
\xymatrix{
	P_n \ar[d]_{(\hatmap)} \ar[r]^{(\checkmap)}
	& \End(B)_n \ar[d]^{-\ocomp(F,\dots,F)} \\
	\End(A)_n \ar[r]^-{F\ocomp -} & \C(A, \dots, A ; B)
}
\]
commutes.
\end{defn}

The definition of morphism may be stated equivalently in terms of any of the
three characterizations of algebras given above.

\section{Symmetric operads}
\label{sec:symm_operads}

\begin{defn}
\label{def:symm_mcat}
\index{multicategory!symmetric}
A \defterm{symmetric multicategory} is a multicategory \C\ and, for every $n
\in \natural$, every $\sigma \in S_n$, and every $A_1, \dots A_n, B \in \C$, a
map
\[
\begin{array}{rccl}
\sigma \act - : & \C(A_1, \dots, A_n; B) & \longrightarrow 
& \C(A_{\sigma 1}, \dots, A_{\sigma n}; B) \\
& f & \longmapsto & \sigma \act f 
\end{array}
\]
such that
\begin{itemize}
\item For each $f \in \C(A_1, \dots, A_n; B)$, $1\act f = f$.
\item For each $\sigma, \rho \in S_n$, and each $f \in \C(A_1, \dots, A_n; B)$,
\[
\rho \act (\sigma \act f) = (\rho\sigma) \act f
\]
\item For each permutation $\sigma \in S_n$, all objects $A^1_1, \dots,
A^n_{k_n}, B_1, \dots, B_n, C \in \C$ and all arrows $f_i \in \C(A^i_1, \dots,
A^i_{k_i}; B_i)$ and $g \in \C(B_1, \dots, B_n; C)$,
\[
(\sigma\act g) \ocomp (f_{\sigma 1}, \dots, f_{\sigma n})
	= (\sigma \ocomp (1,\dots, 1)) \act (g \ocomp (f_1, \dots f_n)).
\]
\item For each $A^1_1, \dots, A^n_{k_n}, B_1, \dots, B_n, C \in \C$,
$\sigma_i \in S_{k_i}$ for $i=1,\dots, n$, and each $f_i \in \C(A^i_1, \dots,
A^i_{k_i}; B_i), g \in \C(B_1, \dots, B_n; C)$,
\[
g \ocomp (\sigma_1 \act f_1, \dots, \sigma_n \act f_n)
= (1 \ocomp (\sigma_1, \dots, \sigma_n)) \act (g\ocomp (f_1, \dots, f_n)).
\]
\end{itemize}
where $\sigma\ocomp(1, \dots, 1)$ and $1\ocomp(\sigma_1, \dots, \sigma_n)$ are
as defined in Example \ref{ex:symmopd}.
\end{defn}
This definition is unusual in that the symmetric groups act on the left rather
than on the right as is more common: however, this change is essential for our
later generalization to \dms\ in Section \ref{sec:dops}.

\begin{defn}
\index{morphism!of symmetric multicategories}
Let $\C_1$ and $\C_2$ be symmetric multicategories.
A \defterm{morphism} (or \defterm{map}) $F$ of symmetric multicategories is a
map $F: \C_1 \to \C_2$ of multicategories such that $F(\sigma \act f) = \sigma
\act F(f)$ for all $n \in \natural$, all \nary\ $f$ in $\C_1$, and all $\sigma
\in S_n$.
\end{defn}

\begin{defn}
\index{algebra!for a symmetric multicategory}
Let $M$ and $\C$ be symmetric multicategories.
An \defterm{algebra} for $M$ in \C\ is a morphism of symmetric multicategories
$M \to \C$.
\end{defn}

\begin{defn}
\index{operad!symmetric}
A \defterm{symmetric operad} is a symmetric multicategory
with only one object.
\end{defn}
In this case, the definition is equivalent to the following:
\begin{defn}
\label{def:symmopd}
\xymatrixrowsep{5pc}
\xymatrixcolsep{6pc}
A \defterm{symmetric operad} is an operad $P$ together with an action of the
symmetric group $S_n$ on each $P_n$, which is compatible with the operadic
composition:
\[
\xymatrixrowsep{2pc}
\xymatrix{
	P_n \times \prod P_{k_i} \ar[r]^{(\sigma \act -) \times 1 \times \dots
		\times 1}
		\ar[d]_{1 \times \sigma_*} 
	& P_n \times \prod P_{k_i} \ar[dd]^{\ocomp} \\
	P_n \times \prod P_{\sigma k_i} \ar[d]_\ocomp
	\\
	P_{\sum k_i} \ar[r]^{(\sigma \ocomp (1, \dots, 1)) \act -}
	& P_{\sum k_i}
}
\xymatrix{
	P_n \times \prod P_{k_i}
		\ar[r]^{1 \times (\rho_1 \act -) \times \dots
			  \times (\rho_n \act -)}
		\ar[dd]_{\ocomp}
	& P_n \times \prod P_{k_i} \ar[dd]^{\ocomp} \\
	\\
	P_{\sum k_i} \ar[r]^{(1 \ocomp (\rho_1, \dots, \rho_n)) \act -}
	& P_{\sum k_i}
}
\]
\[
\xymatrix{
	P_n \ar@/u0.5cm/[r]^{1 \act -} \ar@/d0.5cm/[r]_1 & P_n
}
\]
\end{defn}
Maps of symmetric operads are just maps of symmetric multicategories.

\begin{example}
\index{\calS!as a symmetric operad}
The operad $\calS$ of symmetric groups, as given in Example \ref{ex:symmopd}.
The action of $S_n$ on $\calS_n$ is given by $\sigma \cdot \tau = \tau\sigma^{-1}$.
\end{example}

\begin{example}
\label{ex:symmEnd}
\index{endomorphism operad!symmetric}
Let $\C$ be a symmetric multicategory, and $A \in \C$.
The \defterm{symmetric endomorphism operad} $\End(A)$ of $A$ is the full sub-(symmetric multicategory) of $\C$ whose only object is $A$.
\end{example}

If $\C$ is the underlying symmetric multicategory of a symmetric monoidal
category, then $\End(A)_n = \C(A^{\otimes n}; A)$ for each $n \in \natural$, and
the actions of the symmetric groups are given by composition with the symmetry
maps.

\begin{defn}
\index{algebra!for a symmetric operad}
\index{morphism!of algebras for a symmetric operad}
Let $P$ be a symmetric operad.
An \defterm{algebra} for $P$ in a multicategory $\C$ is an object $A$ and a map $h : P \to \End(A)$ of symmetric operads.
A \defterm{morphism} $(A, h) \to (A', h')$ of $P$-algebras is an arrow $F: A
\to B$ in \C\ such that $h'F = Fh$.
\end{defn}
As with plain operads, the definitions of an algebra for a symmetric operad $P$
and of morphisms between those algebras may be stated in several equivalent
ways.

\section{\Dops}
\label{sec:dops}
The definition of categorification in Chapter \ref{ch:categorification} is
couched in terms of operads.
To generalize it, therefore, we might generalize the definition of operad so
that it is capable of expressing every (one-sorted) algebraic theory.

This generalization is not new: our ``\dops'' were presented by Tronin
under the name ``{\bf FinSet}-operads''. Our Theorem \ref{thm:DOpCat iso Clone}
appears in \cite{tronin}, and Theorem \ref{thm:dop-clone alg} appears as
Theorem 1.2 in \cite{tronin2}.
A fuller treatment was given by T.~Fiore (who called them ``the functional
forms of theories'') in \cite{fiore}.
Tronin's paper constructs an isomorphism between the category of \dops\ and the
category of algebraic clones which commutes with the forgetful functors to
$\Set^\natural$; Fiore's constructs an equivalence between the
category of \dops\ and that of Lawvere theories, and also shows that this
equivalence preserves the categories of algebras.

Let \F\ be a skeleton of the category of finite sets and functions, with
objects the sets \fs 0, \fs 1, \fs 2, \dots, where $\fs n = \{1,2,\dots,n\}$.
\index{\F}
\begin{defn}
\index{multicategory!finite product}
\label{dopdef}
A \defterm{\dm} is:
\begin{itemize}
\item A plain multicategory \C;
\item for every morphism $f : \n \to \m$ in \F, and for all objects $C_1,
\dots, C_n, D \in \C$, a function $f\act - : \C(C_1, \dots, C_n;D) \to
\C(C_{f(1)}, \dots, C_{f(n)}; D)$
\end{itemize}
satisfying the following axioms:
\begin{itemize}
\item the \F-action is functorial: $f\act(g\act p) = (f\ocomp g)\act p$, and
$\id_{\fs n}\act p = p$ wherever these equations make sense;
\item the \F-action and multicategorical composition interact by ``combing
out'':
\[
(f\act p)\cmp{f_1\act p_1, \dots, f_n\act p_n} = (f\cmp{f_1,\dots, f_n)} \act
(p\cmp{p_{f(1)},\dots, p_{f(n)}})
\]
where $(f\cmp{f_1,\dots, f_n})$ is given as follows:

Let $f:\fs n \to \fs m$, and $f_i:{\fs k}_i \to {\fs j_i}$ for $i = 1, \dots,
n$.
Then
\[
\begin{array}{rcccl}
f\cmp{f_1, \dots, f_n} & : & \fs{\sum k_i} & \to & \fs{\sum j_i} \\
f\cmp{f_1,\dots, f_n} & : & \left(\sum_{i=1}^{p-1} k_{f(i)}\right) + h
	& \mapsto & \left(\sum_{i=1}^{f(p)-1} j_i\right) + f_p(h)
\end{array}
\]
for all $p \in \{1, \dots, n\}$ and all $h \in \{1, \dots, k_{f(p)}\}$.
See Figure \ref{combfig}.
The small specks represent inputs to the arrow that are ignored.
\end{itemize}
\begin{figure}[h]
\centerline{
	\epsfysize=2in
	\epsfbox{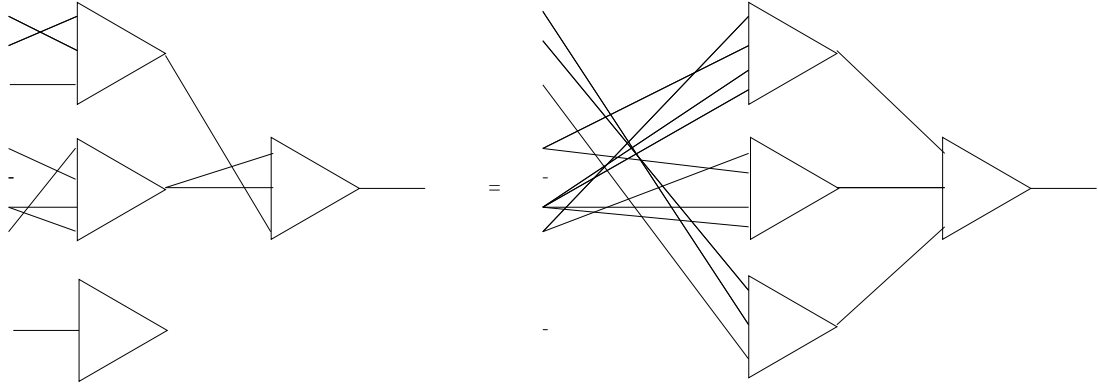}
}
\caption{``Combing out'' the \F-action}
\label{combfig}
\end{figure}
\end{defn}
It is now possible to see why we chose to have our symmetries acting on the left
in Definition \ref{def:symm_mcat}: in this more general case, only a left
action is possible.

\begin{defn}
\index{operad!finite product}
A \defterm{\dop} is a \dm\ with only one object.
\end{defn}

We will see in Section \ref{sec:synclass} that \dops\ are equivalent in
expressive power to Lawvere theories or clones: hence, every finitary algebraic
theory provides an example of a \dop.
As before, the sets $P_n$ contain the \nary\ operations in the theory.
For illustrative purposes, we work out two examples now:
\begin{example}
\index{algebra!for a ring}
\label{ex:Ralg_fp}
Let $R$ be a ring, and $P_n = R[x_1, \dots, x_n]$ (the set of polynomials in
$n$ commuting variables over $R$) for all $n \in \natural$.
If $p \in P_n$ and $q_i \in P_{k_i}$ for $i = 1, \dots, n$, then
\[
(p\ocomp(q_1, \dots, q_n))(x_1, \dots, x_{\sum_{i=1}^n k_i}) =
	p(q_1(x_1, \dots, x_{k_1}), \dots, q_n(x_{(\sum_{i=1}^{n-1} k_i) + 1},
		\dots, x_{\sum_{i=1}^n k_i}))
\]
and if $f: \n \to \m$, then
\[
(f \cdot p)(x_1, \dots, x_m) = p(x_{f(1)}, \dots, x_{f(n)})
\]
\end{example}

\begin{example}
\index{commutative monoids}
\label{ex:comm_monoid_fp}
Let $P_n$ be the set of elements of the free commutative monoid on $n$
variables $x_1, \dots, x_n$.
Elements of $P_n$ are in one-to-one-correspondence with  elements of
$\natural^n$.
We call the $n$th component of $p \in P_n$ the \defterm{multiplicity of the
$n$th argument}.
Composition is defined as follows:
\def\vector#1#2{\left[\begin{array}{c} {#1} \\ \vdots \\ {#2}\end{array}\right]}
\def\veclist#1#2{\vector{#1_1}{#1_{#2}}}
\[
\veclist p n \ocomp \left(\veclist {q^1} {k_1}, \dots, \veclist {q^n}
	{k_n}\right)
= \vector {p_1 q_1^1} {p_n q^n_{k_n}}
\]
and if $f : \n \to \m$,
\[
f \cdot \veclist p n = \vector {\sum_{f(i)=1} p_i} {\sum_{f(i)=m} p_i}
\]
Or, in more familiar notation:
\begin{eqnarray*}
(x_1^{p_1} \dots x_n^{p_n}) \ocomp (x_1^{q^1_1} \dots x_1^{q^1_{k_1}}, \dots
x_{(\sum_{i=1}^{n-1} k_i) + 1}^{q^n_1} \dots x_{\sum_{i=1}^n k_i}^{q^n_{k_n}})
& = &x_1^{p_1 q^1_1} x_2^{p_1 q^1_1} \dots x_{\sum_{i=1}^n k_i}^{p_n q^n_{k_n}}
\\
f\cdot(x_1^{p_1} \dots x_n^{p_n}) & = & x_{f(1)}^{p_1} \dots x_{f(n)}^{p_n} \\
 & = & x_1^{\sum_{f(i)=1} p_i} \dots x_m^{\sum_{f(i)=m} p_i}
\end{eqnarray*}
\end{example}

\begin{example}
\label{ex:fpEnd}
\index{endomorphism operad!finite product}
Let $\C$ be a finite product category, and $A$ be an object of $\C$.
Then there is a \dop\ $\End(A)$, the \defterm{endomorphism operad} of $A$,
where $\End(A)_n = \C(A^n,A)$, and $f\act p$ is $p$ composed with the
appropriate combination of projections to relabel its arguments by $f$.
\end{example}

\begin{defn}
\index{morphism!of \dms}
Let $M$, $N$ be \dms.
A \defterm{morphism} $F: M \to N$ consists of
\begin{itemize}
\item for each object $m \in M$, an object $Fm \in N$;
\item for each $n \in \natural$ and all $m_1, \dots, m_n, m \in M$, a map
\[
F_{m_1, \dots, m_n, m}: M(m_1, \dots m_n; m) \to N(Fm_1, \dots, Fm_n; Fm)
\]
commuting with the \F-action, the unit and composition.
\end{itemize}
\end{defn}
\begin{defn}
\index{algebra!for a \dm}
Let $M$ be a \dm.
An \defterm{algebra} for $M$ in a \dm\ $\C$ is a map of \dms\ $M \to \C$.
An \defterm{algebra} for $M$ in a finite product category $\C$ is a map of
\dms\ from $M$ to the underlying \dm\ of \C.
\Dms\ and their morphisms form a category called \DMCat. \index{\DOpCat}
\end{defn}

In the special case of \dops, these definitions are equivalent to the
following:
\begin{defn}
\index{morphism!of \dops}
Let $P$, $Q$ be \dops.
A \defterm{morphism} $F: P \to Q$ is a sequence of maps $F_i: P_i \to
Q_i$ commuting with the \F\ action, the unit and composition.
\end{defn}
\begin{defn}
\index{algebra!for a \dop}
Let $P$ be a \dop.
An \defterm{algebra} for $P$ in a finite product category $\C$ is an object $A
\in \C$ and a map of \dops\ $P \to \End(A)$.
\end{defn}
\Dops\ and their morphisms form a category called \DOpCat. \index{\DOpCat}

\begin{example}
\index{algebra!for a ring}
The algebras in $\C$ for the operad described in Example \ref{ex:Ralg_fp} are
associative $R$-algebras in $\C$.
\end{example}

\begin{example}
\index{commutative monoids}
The algebras for the operad described in Example \ref{ex:comm_monoid_fp} are
commutative monoid objects in $\C$.
\end{example}

\begin{theorem}
\label{thm:DOpCat iso Clone}
$\DOpCat \cong \Clone$.
\end{theorem}
\begin{proof}
We shall construct a functor $K_{(-)} : \DOpCat \to \Clone$, and show that it is
bijective on objects, full and faithful.

If $P$ is a \dop, let $K_P$ be the following clone:
\begin{itemize}
\item $(K_P)_n = P_n$ for all $n \in \natural$,
\item composition is given by composition in $P$: if $p \in P_n$ and $p_1,
\dots, p_n \in P_m$, then $p\ccomp(p_1, \dots, p_n) \in (K_P)_m$ is $f \act
(p\ocomp(p_1, \dots, p_n)) \in P_m$, where
\begin{equation}
\label{eqn:KPcomp}
\begin{array}{lccl}
f : &\fs{nm} &\to & \fs{m} \\
& x & \mapsto & ((x - 1) \mod m) + 1,
\end{array}
\end{equation}
\item for all $n \in \natural$ and all $i \in \fs{n}$, the projection
$\delta^i_n$ is $f^i_n \act 1$, where
\begin{equation}
\label{eqn:KPproj}
\begin{array}{lccl}
f^i_n : &\fs{1} &\to & \fs{n} \\
& 1 & \mapsto & i.
\end{array}
\end{equation}
\end{itemize}
It is easily checked that $K_P$ satisfies the axioms for a clone given in
Definition \ref{def:clone}.

On morphisms, $K_{(-)}$ acts trivially: morphisms of clones and of \dops\ are
simply maps of signatures commuting with the extra structure, and $K_{(-)}$
preserves the underlying map of signatures.

Let $K$ be a clone.
Let $P_K$ be the \dop\ for which
\begin{itemize}
\item $(P_K)_n = K_n$ for all $n \in \natural$,
\item $1 = \delta^1_1$,
\item $p \ocomp (p_1, \dots, p_n) = p \ccomp (p_1 \ccomp (\delta^1_m, \dots,
\delta^{k_1}_{\sum k_i}), \dots, p_n \ccomp (\delta^{k_1 + \dots + k_{n-1} +
1}_{\sum k_i}, \dots, \delta^{\sum k_i}_{\sum k_i}))$ for all $n$ and $k_1,
\dots, k_n \in \natural$, all $p \in K_n$, and all $p_1 \in K_{k_1}, \dots, p_n
\in K_{k_n}$,
\item $f \act p = p \ccomp (\delta^{f(1)}_m, \dots, \delta^{f(n)}_m)$ for all
$n, m \in \natural$, all $f : \fs{n} \to \fs{m}$ and all $p \in K_n$.
\end{itemize}
We will show that $K_{P_K} = K$ for all $K \in \Clone$, and that $P_{K_P} = P$
for all $P \in \DOpCat$.
Let $K$ be a clone.
Then $(K_{P_K})_n = (P_K)_n = K_n$ for all $n \in \natural$.
If $n, m \in \natural$, $k \in K_n$ and $k_1, \dots, k_n \in K_m$, then the
composite $k\ccomp(k_1, \dots, k_n)$ in $K_{P_K}$ is given by the composite 
$f \act (k\ocomp(k_1, \dots, k_n))$ in $P_K$, where $f$ is given by
(\ref{eqn:KPcomp}) above.
This in turn is given by the composite
\begin{eqnarray*}
(k & \ccomp & (k_1 \ccomp (\delta^1_{nm},
\dots, \delta^{m}_{nm}), \dots, k_n \ccomp (\delta^{(n-1)m + 1}_{nm}, \dots,
\delta^{nm}_{nm}))) \\
& \ccomp & (\delta^1_m, \dots, \delta^m_m, \dots, \delta^1_m, \dots, \delta^m_m)
\end{eqnarray*}
in $K$.
By the associativity law for clones, this is equal to
\[
\begin{array}{rl}
k \ccomp (& k_1 \ccomp (\delta^1_{nm}, \dots, \delta^{m}_{nm})
 \ccomp (\delta^1_m, \dots, \delta^m_m, \dots, \delta^1_m, \dots, \delta^m_m), \\
& \dots, \\
& k_n \ccomp (\delta^{(n-1)m + 1}_{nm}, \dots,\delta^{nm}_{nm}) 
 \ccomp (\delta^1_m, \dots, \delta^m_m, \dots, \delta^1_m, \dots, \delta^m_m))
\end{array}
\]
which in turn may be simplified to $k \ccomp (k_1, \dots, k_n)$ as required.
For every $n \in \natural$ and every $i \in \fs{n}$, the projection
$\delta^i_n$ in $K_{P_K}$ is given by $f^i_n \act 1$, where $f^i_n$ is defined
in (\ref{eqn:KPproj}): this in turn is given by $1 \ocomp
(\delta^i_n) = \delta^1_1 \ccomp (\delta^i_n) = \delta^i_n$.
Hence $K_{P_K} = K$.

Conversely, let $P$ be a \dop\; we shall show that $P_{K_P} = P$.
For every $n \in \natural$, the set $(P_{K_P})_n$ is equal to $P_n$.
The unit element is given by $1 = \delta^1_1 = f^1_1 \act 1 = 1$.
If $p \in P_n$ and $p_i \in P_{k_i}$ for $i = 1, \dots, n$, then the composite
$p\ocomp(p_1, \dots, p_n)$ is given by $p \ccomp (p_1 \ccomp (\delta^1_m, \dots,
\delta^{k_1}_{\sum k_i}), \dots, p_n \ccomp (\delta^{k_1 + \dots + k_{n-1} +
1}_{\sum k_i}, \dots, \delta^{\sum k_i}_{\sum k_i}))$ in $K_P$, which in turn
is given by (after simplification) $p \ocomp (p_1, \dots, p_n)$ in $P$.
Hence $P = P_{K_P}$, and $K_{(-)}$ is bijective on objects.
The reasoning above also suffices to show that $K_{(-)}$ is well-defined on
morphisms and full (since preserving a \dop\ structure amounts exactly to
preserving the associated clone structure).
Since the morphisms of both categories are simply maps of signatures with extra
properties and $K_{(-)}$ commutes with the forgetful functors to
$\Set^\natural$,
then $K_{(-)}$ is faithful.
Hence $K_{(-)}$ is an isomorphism of categories, and $\DOpCat \cong \Clone$.
\end{proof}

\begin{theorem}
\label{thm:dop-clone alg}
Let $P$ be a \dop.
Then $\Alg P \cong \Alg{K_P}$.
\end{theorem}
\begin{proof}
Let $(A, (\hatmap))$ be a $P$-algebra.
Since the elements of the finite product endomorphism operad $\End(A)$ are
endomorphisms of $A$, and composition is given by composition of morphisms,
then $K_{\End(A)} = \End(A)$, the endomorphism clone of $A$.
Since the functor $K_{(-)} : \DOpCat \to \Clone$ is an isomorphism, a morphism
of \dops\ $P \to \End(A)$ is exactly a map of clones $K_P \to \End(A)$.
Hence an algebra for $P$ is exactly an algebra for $K_P$.
A morphism between $P$-algebras is a morphism between their underlying objects
that commutes with $\hat p$ for every $p \in P_n$ and every $n \in \natural$;
this is true iff it commutes with $\hat k$ for every $k \in (K_P)_n$ and every
$n \in \natural$.
\end{proof}

\section{Adjunctions}
\label{sec:adj}

In the next few sections, we shall show that there is a chain of monadic
adjunctions
\index{adjunctions!for operads}
\begin{equation}
\label{eqn:adjs}
\xymatrix{
  \DOpCat \ar@<1.2ex>[d]^\Ufpsig
	 \ar@/r1.4cm/[dd]^\Ufpp
	 \ar@/r3cm/[ddd]^\Ufp \\
  \SymmOperad \ar@<1.2ex>[d]^\Usigp \ar@<1.2ex>[u]^\Ffpsig_{\dashv}
	 \ar@/r1.4cm/[dd]^\Usig \\
  \Operad \ar@<1.2ex>[u]^\Fsigp_{\dashv} \ar@<1.2ex>[d]^\Up
	 \ar@/l1.4cm/[uu]^\Ffpp \\
  \Set^\natural \ar@<1.2ex>[u]^\Fp_{\dashv} \ar@/l1.4cm/[uu]^\Fsig 
  	 \ar@/l3cm/[uuu]^\Ffp
}
\end{equation}
The notation is chosen such that $F^x_y \dashv U_x^y$, and $U_x^y U_y^z =
U_x^z$.
The notation is inspired by the exponential notation used for hom-objects: the
source category of one of these functors is determined by its superscript, and
the target category is determined by its subscript.
The ``pl'' stands for ``plain''.
A similar chain of adjunctions (for PROPs rather than operads) was discussed in
\cite{baez_ualg}, pages 51--59.

We refer to the monad $U_x^yF^x_y$ as $T^x_y$.
The right adjoints $\Up, \Usigp$ and $\Ufpsig$ are found by forgetting
respectively the compositional structure, the symmetric structure, and the
actions of all non-bijective functions, and will not be described further.
By standard properties of adjunctions, the composite functors are adjoint:
$\Fsig \dashv \Usig$ etc.

\section{Existence and monadicity}

All the left adjoints in (\ref{eqn:adjs}) are examples of a more general
construction.
We shall now investigate this general case, and show that the adjunction which
arises is always monadic.

But first, we have so far only asserted that $\Ufp, \Ufpsig$ and $\Ufpp$ have
left adjoints.
We shall show that these left adjoints must exist for general reasons.

Let $\FP$ be the category of small categories with finite products and
product-preserving functors.

\begin{lemma}
\label{lem:frees_exist}
Let $\Csmall$ and $\Dsmall$ be small finite-product categories, let $\C$ be
cartesian closed and have all small colimits, and let $Q: \Csmall \to \Dsmall$
preserve finite products.
Then the adjunction
\[
\adjunction{[\Dsmall,\C]}{[\Csmall,\C]}{Q_!}{Q^*},
\]
where $Q^*$ is composition with $Q$ and $Q_! = \Lan_Q$, restricts to an
adjunction
\[
\adjunction{\FP(\Dsmall,\C)}{\FP(\Csmall,\C)}{Q_!}{Q^*},
\]
\label{lem:leftadj}
\end{lemma}
\begin{proof}
\def\AC{[\Csmall, \C]}
\def\BC{[\Dsmall, \C]}
\def\FPAC{\FP(\Csmall, \C)}
\def\FPBC{\FP(\Dsmall, \C)}
Certainly $Q^*$ restricts in this way, since $Q$ preserves finite products.
$\FPAC$ and $\FPBC$ are full subcategories of $\AC$
and $\BC$, so if we can show that $Q_!$ restricts to a functor $\FPAC \to
\FPBC$, then it is automatically left adjoint to the restriction of $Q^*$.

Let $X : \Csmall \to \C$ preserve finite products.
We must show that $Q_! X : \Dsmall \to \C$ preserves finite products.
We shall proceed by showing that $Q_! X$ preserves terminal objects and binary
products.

Recall that
\[
(Q_! X)(b) \cong \int^a \Dsmall(Q a, b) \times X a
\]
for all $b \in \Dsmall$.
Hence, using $1$ for the terminal objects in $\Dsmall$ and $\C$,
\begin{eqnarray*}
(Q_! X)(1) & \cong & \int^a \Dsmall(Q a, 1) \times X a \\
& \cong & \int^a 1 \times X a \\
& \cong & \int^a X a \\
& \cong & X 1 \\
& \cong & 1
\end{eqnarray*}
since $X$ preserves finite products and the colimit of a diagram $D$ over a
category with a terminal object $1$ is simply $D 1$.

Now, let $b_1, b_2 \in \Dsmall$.
Then
\begin{eqnarray}
\label{eq:defQ!}
\nonumber\lefteqn{(Q_! X)(b_1 \times b_2)} \\
& \cong & \int^a \Dsmall(Q a, b_1 \times b_2) \times X a
\\
\label{eq:defprod1}
& \cong & \int^a \Dsmall(Q a, b_1) \times \Dsmall(Q a, b_2) \times X a \\
\label{eq:density1}
& \cong & \int^a
	\left(\int^{c_1} \Dsmall(Q c_1, b_1) \times \Csmall(a, c_1) \right)
	\times
	\left(\int^{c_2} \Dsmall(Q c_2, b_2) \times \Csmall(a, c_2) \right)
	\times Xc \\
\label{eq:dist1}
& \cong & \int^{a, c_1, c_2} \Dsmall(Q c_1, b_1) \times \Dsmall(Q c_2, b_2) \times
	\Csmall(a, c_1) \times \Csmall(a, c_2) \times X a \\
\label{eq:defprod2}
& \cong & \int^{a, c_1, c_2} \Dsmall(Q c_1, b_1) \times \Dsmall(Q c_2, b_2) \times
	\Csmall(a, c_1 \times c_2) \times X a \\
\label{eq:dist2}
& \cong & \int^{c_1, c_2} \Dsmall(Q c_1, b_1) \times \Dsmall(Q c_2, b_2) \times
	\left( \int^a \Csmall(a, c_1 \times c_2) \times X a \right) \\
\label{eq:density2}
& \cong & \int^{c_1, c_2} \Dsmall(Q c_1, b_1) \times \Dsmall(Q c_2, b_2) \times
	X(c_1 \times c_2) \\
\label{eq:defprod3}
& \cong & \int^{c_1, c_2} \Dsmall(Q c_1, b_1) \times \Dsmall(Q c_2, b_2) \times
	X c_1 \times X c_2 \\
\label{eq:Xprod}
& \cong & \int^{c_1, c_2} \Dsmall(Q c_1, b_1) \times X c_1 \times
	\Dsmall(Q c_2, b_2) \times X c_2 \\
\label{eq:rearrange}
& \cong & \left (\int^{c_1} \Dsmall(Q c_1, b_1) \times X c_1 \right) \times
	\left( \int^{c_2} \Dsmall(Q c_2, b_2) \times X c_2 \right) \\
\label{eq:dist3}
& \cong & (Q_!X)(b_1) \times (Q_!X)(b_2)
\end{eqnarray}
(\ref{eq:defQ!}) is the definition of $Q_!$;
(\ref{eq:defprod1}), (\ref{eq:defprod2}) and (\ref{eq:defprod3}) are from the
definition of products;
(\ref{eq:density1}) and (\ref{eq:density2}) are applications of the Density
Formula;
(\ref{eq:dist1}), (\ref{eq:dist2}) and (\ref{eq:dist3}) use the
distributivity of products over colimits in \C\ (since \C\ is cartesian
closed),
and (\ref{eq:Xprod}) uses the fact that $X$ preserves finite products.

So $Q_!$ preserves terminal objects and binary products, and hence all finite
products.
\end{proof}

\begin{corollary}
The functors $\Ufpsig, \Ufpp$ and $\Ufp$ all have left adjoints.
\end{corollary}

\begin{lemma}
\label{lem:monadj}
Let $S$ be a set, whose elements we will call \defterm{sorts}.
Let $T$ and $T'$ be $S$-sorted finite product theories, such that $T'$ is a
subcategory of $T$ and the inclusion of $T'$ into $T$ preserves finite
products.
Let $\Alg T$ be the category of $T$-algebras and morphisms in some finite
product category $\C$, and $\Alg{T'}$ be the category of $T'$-algebras and
morphisms in $\C$.
Then the free/forgetful adjunction
\[ 
\adjunction {\Alg{T'}} {\Alg{T}} F U
\]
is monadic, provided the left adjoint $F$ exists.
\end{lemma}
\begin{proof}
\def\prods#1{\prod #1 s_i}
\def\prodr#1{\prod #1 r_j}
\def\prodars#1{\prod #1_{s_i}}
\def\prodarr#1{\prod #1_{r_j}}
We will make use of Beck's theorem to prove monadicity: precisely, we shall
make use of the version in \cite{catwork} VI.7.1, which states that $U$ is
monadic if it has a left adjoint and it strictly creates coequalizers for
$U$-absolute coequalizer pairs.
Recall that a functor $G: \C \to \D$ \defterm{strictly creates coequalizers}
for a diagram $\parallelpair A f g B$ in $\C$ if, for every coequalizer $e: GB
\to E$ of $Gf$ and $Gg$ in $\D$, there are a unique object $E'$ in $\C$ and a
unique arrow $e' : B \to E'$ such that $GE' = E$ and $Ge' = e$, and moreover
that $e'$ is a coequalizer of $\parallelpair A f g B$.

Let $\parallelpair A f g B$ be a $U$-absolute coequalizer pair in $\Alg{T}$,
and $e: U B \to E$ be the coequalizer of $\parallelpair{UA}{Uf}{Ug}{UB}$.
We wish to extend $E$ to a functor $E': T \to \C$.
Define $E'$ to be equal to $E$ on objects.
On arrows, we shall define $E'$ using the universal property of $E$ and the
$U$-absolute property of $\parallelpair A f g B$.

For each arrow $\phi : s_1 \times \dots \times s_n \to r_1 \times \dots
\times r_m$ in $T$ (where $s_i, r_j \in S$),
consider the diagram
\begin{equation}
\label{eqn:coeq}
\xymatrixrowsep{4pc}
\xymatrixcolsep{4pc}
\xymatrixnocompile{
	\prods{U A} \parallelars{\prodars{U f}}{\prodars{Ug}} \ar[d]_{A\phi}
	& \prods{U B} \ar[r]^{\prodars e} \ar[d]_{B\phi}
	& \prods{E} \unar[d]^{E'\phi}\\
	\prodr{U A} \parallelars{\prodarr{f}}{\prodarr{g}}
	& \prodr{U B} \ar[r]^{\prodarr e}
	& \prodr E
}
\end{equation}
in $\C$, where $e: UB \to E$ is a coequalizer for $\parallelpair {UA}
{Uf} {Ug} {UB}$.

Since $\parallelpair A f g B$ is a $U$-absolute coequalizer pair, $\prodars{e}
: \prods{UB} \to \prods{E}$ is a coequalizer.
Since $f$ and $g$ are $T$-homomorphisms, (\ref{eqn:coeq}) serially commutes, so
$(\prodarr e)\phi$ factors uniquely through $\prodars e$.
Define $E'\phi$ to be this map, as shown (and note that $E'\phi = E\phi$ if
$\phi$ is in $T'$).
This definition straightforwardly makes $E'$ into a functor $T \to \C$.
Since $E$ is a $T'$-algebra, and products in $T$ are the same as products in
$T'$, we may deduce that $E': T\to C$ preserves finite products, and thus is
a $T$-algebra.
Clearly, $E'$ is the unique extension of $E$ to a $T$-algebra such that $e$ is
a $T$-algebra morphism.
It remains to show that $e$ is a coequalizer map for $\parallelpair A f g B$ in
$\Alg{T}$.

Suppose $\fork A f g B d D$ is a fork in $\Alg T$.
Then $\fork {UA} {Uf} {Ug} {UB} {Ud} {UD}$ is a fork in $\Alg{T'}$, so $Ud$
factors through $e$; say $Ud = he$.
We must show that $h$ is a $T$-homomorphism.
As before, take $\phi : s_1 \times \dots \times s_n \to r_1 \times \dots \times
r_m$ in $T$, and consider the diagram
\begin{equation}
\xymatrixcolsep{4pc}
\xymatrixnocompile{
	\prods{U A} \parallelars{\prodars{U f}}{\prodars{U g}} \ar[dd]_{A\phi}
	& \prods{U B} \ar[r]^{\prodars{e}}
		\ar[dd]_{B\phi} \ar[dr]_{\prodars{Ud}}
	& \prods{E'} \ar@/r1cm/[dd]^{E'\phi} \unar[d]_{\prodars{h}} \\
	& & \prods{U D} \ar@/r1cm/[dd]^{D\phi} \\
	\prodr{U A} \parallelars{\prodarr f}{\prodarr g}
	& \prodr{U B} \ar[r]^{\prodarr e} \ar[dr]_{\prodarr{Ud}}
	& \prodr{E'} \unar[d]_{\prodarr{h}} \\
	& & \prodr{D}
}
\end{equation}

We must show that the curved square on the far right commutes.
Now $(D\phi) \fcomp (\prodars d) = (D\phi) \fcomp (\prodars h) \fcomp (\prodars
e)$, and $(\prodarr{Ud})\fcomp(UB\phi) = h\fcomp e\fcomp(UB\phi) = h \fcomp
(E\phi) \fcomp (\prodars e)$, since $e$ is a $T$-algebra homomorphism.
But $(\prodarr {Ud})\fcomp \phi = \phi \fcomp (\prod d_i)$, so $D\phi \fcomp
(\prodars{h}) \fcomp (\prodars{e}) = h \fcomp (E\phi) \fcomp (\prodars{e})$.
And $\prodars{e}$ is (regular) epic, so $h \fcomp (E\phi) = (D\phi) \fcomp
(\prodars{h})$.

So $h$ is a $T$-algebra homomorphism.
Hence $U$ strictly creates coequalizers for $U$-absolute coequalizer pairs,
and hence is monadic.
\end{proof}

This result could also have been deduced from the Sandwich Theorem of Manes:
see \cite{manes} Theorem 3.1.29 (page 182).

\begin{theorem}
\label{thm:opd monadic over setN}
\index{\Operad!monadicity over $\Set^\natural$}
All the adjunctions in diagram \ref{eqn:adjs}, namely $\Ffpsig \dashv \Ufpsig, \Ffpp \dashv \Ufpp, \Ffp \dashv
\Ufp, \Fsigp \dashv \Usigp, \Fsig \dashv \Usig$ and $\Fp \dashv \Up$, are
monadic.
\end{theorem}
\begin{proof}
Each category mentioned is a category of algebras for some $\natural$-sorted
theory, and the monadicity of each adjunction mentioned is obtained by a simple
application of Lemma \ref{lem:monadj}. For instance, symmetric operads are
algebras for the theory presented by
\begin{itemize}
\item \emph{operations}: one of the appropriate arity for each composition
operation in Definition \ref{def:symmopd}, and an operation $\sigma \cdot -$
for each $n \in \natural$ and each $\sigma$ in $S_n$.
\item \emph{equations}: one for each instance of the axioms in Definition
\ref{def:symmopd}, and an equation $(\sigma \cdot -)\fcomp (\rho \cdot -) =
\sigma\rho \cdot -$ for each $\sigma, \rho \in S_n$ and every $n \in \natural$.
\end{itemize}
\end{proof}

\section{Explicit construction of $\Fp$ and $\Fsigp$}
\label{sec:explicit_Fs}

The previous section showed that $\Fp$ and $\Fsigp$ exist for general reasons,
but it will be useful later to have an explicit construction of these functors.
For this reason, we shall now explicitly construct functors $\Set^\natural \to
\Operad$ and $\Operad \to \SymmOperad$, and prove that they are left adjoint to
$\Up$ and $\Usigp$.

\begin{defn}
\label{def:srterm}
\index{tree!strongly regular}
\index{strongly regular!tree|see{tree, strongly regular}}
Let $\Phi$ be a signature.
An \defterm{\nary\ strongly regular tree labelled by $\Phi$} is an element of
the set $\tr_n \Phi$, which is recursively defined as follows:
\begin{itemize}
\item $|$ is an element of $\tr_1 \Phi$.
\item If $\phi \in \Phi_n$, and $\tau_1 \in \tr_{k_1} \Phi, \dots, \tau_n \in
\tr_{k_n} \Phi$, then $\phi \cmp{\tau_1, \dots, \tau_n} \in
\tr_{\sum k_i} \Phi$.
\end{itemize}
\end{defn}

In graph-theoretic terms, all our trees are planar and rooted.
They need not be level.

We shall abuse notation and write $\phi$ instead of $\phi \cmp{|, \dots, |}$,
for $\phi \in \Phi_n$.

\begin{figure}
\centerline{
	\epsfxsize=5in
	\epsfbox{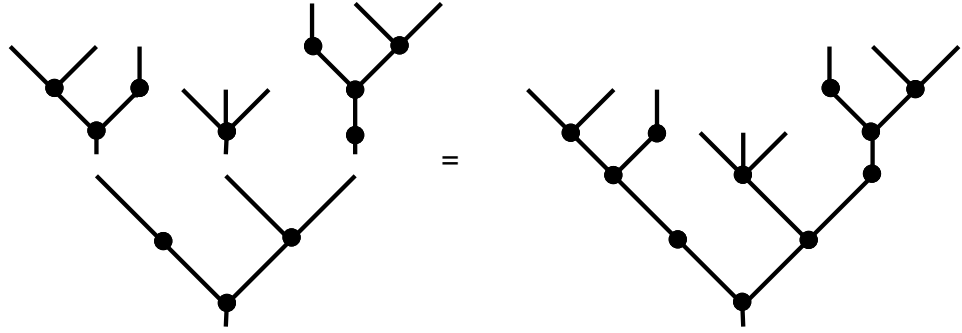}
}
\caption{Grafting of trees}
\label{fig:grafting}
\end{figure}

Given a signature $\Phi$, the objects of the plain operad $(\Fp \Phi)_n$ are
the elements of $\tr_n \Phi$, and composition is given by grafting of trees:
\begin{itemize}
\index{grafting!of trees}
\item $| \cmp \tau$ = $\tau$
\item If $\tau_1 \in \tr_{k_1} \Phi, \dots, \tau_n \in \tr_{k_n} \Phi$, then
\[
(\phi \cmp{\tau_1,\dots, \tau_n}) \cmp{\sigma_1, \dots, \sigma_{\sum k_i}}
 = \phi \cmp{\tau_1 \cmp {\sigma_1, \dots, \sigma_{k_1}}, \dots, 
	\tau_n \cmp {\sigma_{(\sum k_i) - k_n + 1}, \dots, \sigma_{\sum k_i}}}
\]
\end{itemize}
See Figure \ref{fig:grafting}.

The unary tree $|$ is thus the identity in $(\Fp \Phi)$.

$\Fp$ acts on arrows as follows.
Let $f: \Phi \to \Psi$ be a map of signatures.
Then: 
\begin{itemize}
\item $(\Fp f) | = |$
\item $(\Fp f) (\phi \ocomp (\tau_1, \dots, \tau_n))
= (f \phi) \ocomp ( (\Fp f)\tau_1, \dots, (\Fp f)\tau_n )$
\end{itemize}
It is readily verified that with this definition $\Fp$ is a functor
$\Set^\natural \to \Operad$.

We define natural transformations $\eta: 1_{\Set^\natural} \to \Up\Fp$ and
$\epsilon : \Fp\Up \to 1_{\Operad}$ as follows:
\begin{eqnarray}
\eta_\Phi(\phi) & = &\phi \cmp{|, \dots, |} \\
\epsilon_P(|) & = &1_P \\
\epsilon_P(\phi\cmp{\tau_1, \dots, \tau_n)}
	& = & \phi \cmp {\epsilon_P(\tau_1), \dots, \epsilon_P(\tau_n)}
\end{eqnarray}
where $P \in \Operad, \Phi \in \Set^\natural, \phi \in \Phi$, and  $\tau_1, \dots, \tau_n$ are arrows of $P$.

In other words, $\epsilon_P$ is given by applying composition in $P$ to the
formal tree $\Fp \Up P$.

\begin{lemma}
$(\Fp, \Up, \eta, \epsilon)$ is an adjunction.
\end{lemma}
\begin{proof}
We proceed by checking the triangle identities.
We require to show that
\begin{eqnarray}
\label{eqn:triarrFp}
\xymatrix{
	\Fp \ar[r]^{\Fp\eta} \ar[dr]_{1_{\Fp}}
	& \Fp\Up\Fp \ar[d]^{\epsilon\Fp} \\
	& \Fp
} \\
\label{eqn:triarrUp}
\xymatrix{
	\Up \ar[r]^{\eta\Up} \ar[dr]_{1_{\Up}}
	& \Up\Fp\Up \ar[d]^{\Up\epsilon} \\
	& \Up
}
\end{eqnarray}
commute.
For (\ref{eqn:triarrFp}), we proceed by induction on trees.
We shall suppress all subscripts on natural transformations in the interest of
legibility.
For the base case:
\begin{eqnarray*}
\epsilon\Fp ( \Fp\eta (|)) & = & \epsilon\Fp (| \cmp |) \\
	& = & | \cmp{\epsilon(|)} \\
	& = & | \\
	& = & 1_\Fp(|).
\end{eqnarray*}
For the inductive step, let $\Phi$ be a signature, $\phi$ be an \nary\ element of $\Phi$, and $\tau_1, \dots, \tau_n$ be trees labelled by $\Phi$.
Then:
\begin{eqnarray*}
(\epsilon\Fp) ( (\Fp\eta) (\phi \cmp{\tau_1, \dots, \tau_n} ))
	& = & \epsilon\Fp (\phi \cmp{\tau_1, \dots, \tau_n}\cmp{|, \dots, |}) \\
	& = & \phi \cmp{\tau_1, \dots, \tau_n}
		\cmp{(\epsilon\Fp)(|), \dots, (\epsilon\Fp)(|)} \\
	& = & \phi \cmp{\tau_1, \dots, \tau_n} \\
	& = & 1_\Fp(\phi \cmp{\tau_1, \dots, \tau_n})
\end{eqnarray*}

Hence $\epsilon\Fp \fcomp \Fp\eta = 1\Fp$, as required.

For (\ref{eqn:triarrUp}), let $P$ be a plain operad, and let $p$ be an \nary\
arrow in $P$.
\begin{eqnarray*}
(\Up\epsilon) ((\eta\Up) (p)) & = & \Up\epsilon (p\cmp{|, \dots, |}) \\
	& = & p\cmp{(\Up\epsilon)(|), \dots, (\Up\epsilon)(|)} \\
	& = & p\cmp{1, \dots, 1} \\
	& = & p \\
	& = & 1_\Up(p)
\end{eqnarray*}
So $\Up\epsilon \fcomp \eta\Up = 1_\Up$, as required.
\end{proof}

\def\Stimes{\calS \times}
\index{\calS!relation to $\Fsigp$}
We now consider the ``free symmetric operad'' functor $\Fsigp$.
We shall explicitly define a functor $\Stimes - : \Operad \to \SymmOperad$ and
show that it is left adjoint to $\Usigp$, and hence isomorphic to $\Fsigp$.

If $P$ is a plain operad, an element of $(\Stimes P)_n$ is a pair $(\sigma, p)$,
where $p \in P_n$ and $\sigma \in S_n$; i.e., $(\Stimes P)_n = S_n \times P_n$.
Composition is given as follows:
\[
(\sigma, p) \cmp {(\tau_1, q_1), \dots, (\tau_n, q_n)}
 = (\sigma \cmp{\tau_1, \dots, \tau_n}, p \cmp{q_{\sigma(1)}, \dots,
q_{\sigma(n)}})
\]
The symmetric group action is given by $\rho \cdot (\sigma, p) = (\rho
\sigma, p)$. 

\begin{figure}[h]
\centerline{
	\epsfxsize=5in
	\epsfbox{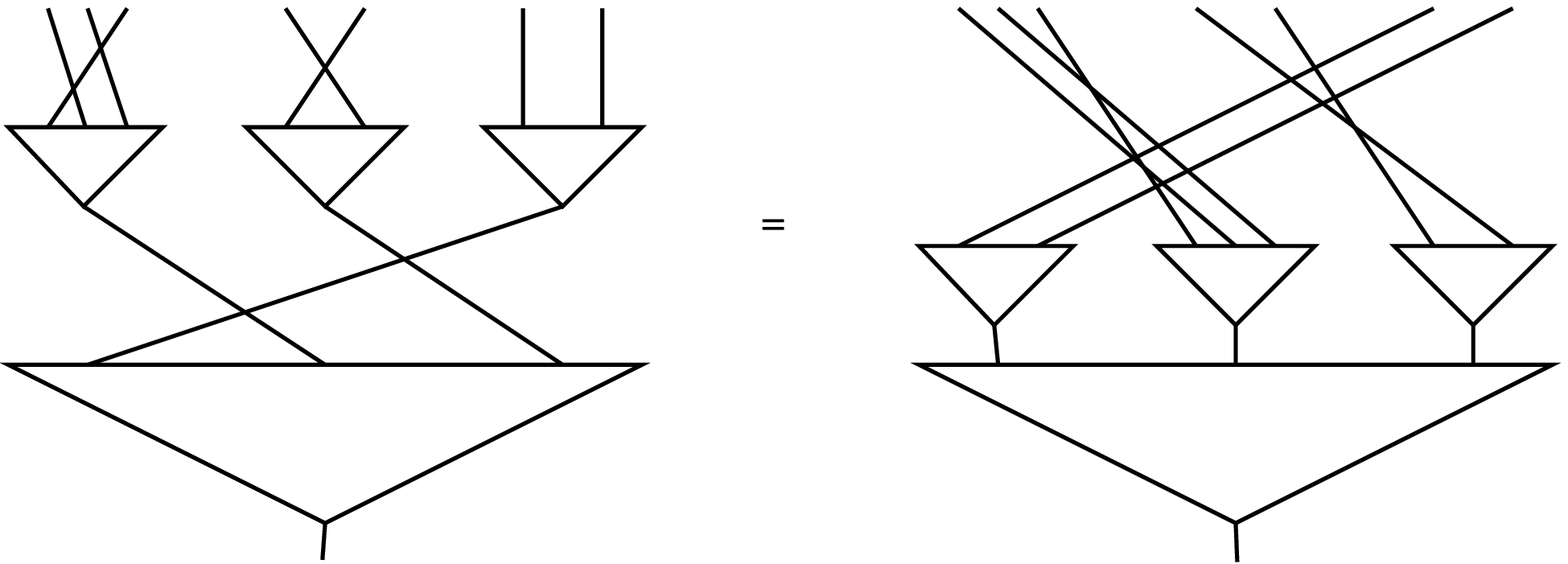}
}
\caption{Composition in $\Stimes P$}
\label{fig:freesym}
\end{figure}

\begin{lemma}
\label{lem:Stimes=Fsigp}
$(\Stimes -) $ is left adjoint to $\Usigp$.
The unit of the adjunction is given by 
\begin{eqnarray*}
\eta: 1 & \to & \Usigp (\Stimes -) \\
\eta_P: p & \mapsto & (1, p),
\end{eqnarray*}
and the counit is given by 
\begin{eqnarray*}
\epsilon : (\Stimes -) \Usigp & \to & 1 \\
\epsilon_P : (\sigma, p) & \mapsto & \sigma \cdot p.
\end{eqnarray*}
\end{lemma}
\begin{proof}
As before, we proceed by checking the triangle identities.
First, let $P$ be a plain operad, $p$ be an \nary\ arrow in $P$, and $\sigma
\in S_n$.
Then $(\sigma, p)$ is an element of $\Stimes P$.
\begin{eqnarray*}
(\epsilon(\Stimes -)) (\Stimes \eta) (\sigma, p))
	& = & (\epsilon(\Stimes -)) (\sigma , (1, p)) \\
	& = & \sigma \cdot (1, p) \\
	& = & (\sigma \fcomp 1, p) \\
	& = & (\sigma, p)
\end{eqnarray*}

Now let $P'$ be a symmetric operad, and $p'$ be an \nary\ arrow in $P'$.
\begin{eqnarray*}
(\Usigp\epsilon) (\eta\Usigp(p')) & = & (\Usigp\epsilon)((1, p')) \\
	& = & 1 \cdot p' \\
	& = & p'
\end{eqnarray*}

So the triangle identities are indeed satisfied, and $(\Stimes -)\dashv \Usigp$.
\end{proof}
Hence, $\Fsigp = \Stimes -$.

\begin{defn}
\label{def:permtree}
\index{tree!permuted}
Let $\Phi$ be a signature.
An \defterm{\nary\ permuted tree labelled by $\Phi$} is an element of $(\Fsigp
\Fp \Phi)_n = (\Fsig \Phi)_n$.
\index{tree!finite product}
An \defterm{\nary\ finite product tree labelled by $\Phi$} is an element of
$(\Ffp \Phi)_n$.
\end{defn}
By Lemma \ref{lem:Stimes=Fsigp}, a permuted tree is a pair $(\sigma, t)$, where
$t \in \tr_n \Phi$ and $\sigma \in S_n$, and (by analogous reasoning) an \nary\
finite product tree is a pair $(f, t)$, where $t \in \tr_m \Phi$ and $f: \m \to
\n$.

\section{Syntactic characterization of the forgetful functors}

There is also a syntactic characterization of the forgetful functor $\Usigp$.
Given a symmetric operad $P$, we take the signature given by all operations in
$P$ (in other words, the signature $\Usig P$).
We then impose all the plain-operadic equations that are true in $P$, and
take the plain operad corresponding to this strongly regular theory.
This operad is $\Usigp P$.

We start by making this precise.

\begin{defn}
\label{def:sreqn}
\index{equation!plain-operadic}
Let $\Phi$ be a signature.
A \defterm{plain-operadic equation over $\Phi$ in $n$ variables} is an element
of $(\nterms \Phi)^2$ (that is, a pair of \nary\ strongly regular trees over
$\Phi$), and a \defterm{plain-operadic equation over $\Phi$} is an element of
$\alleqns \Phi$.
\end{defn}

We shall show that a plain-operadic equation over $\Phi$ is the same thing as a strongly regular equation over $\Phi$.

\begin{defn}
\label{def:opd_presentation}
\index{presentation!for a plain operad}
Let $P$ be a plain operad.
A \defterm{presentation} for $P$ is a signature $\Phi$, a signature $E$, and
maps $e_1, e_2 : \Fp E \to \Fp \Phi $ such that, for some $\phi$,
\[
\NoCompileMatrices
\fork {\Fp E} {e_1} {e_2} {\Fp \Phi} \phi P
\]
is a coequalizer.
We say that a regular epi $\phi : \Fp \Phi \to P$ \defterm{generates} $P$, or
that $\phi$ (or, where the choice of $\phi$ is clear, $\Phi$) is a
\defterm{generator} of $P$.
\end{defn}
Presentations and generators for symmetric and \dops\ are defined analogously.

We will see how these ``presentations'' are related to presentations of
algebraic theories in Section \ref{sec:synclass}.
We now wish to describe the family of all strongly regular equations that are
true in a given symmetric operad $P$: we will then show that this, together
with the signature given by $\Usig P$, is a presentation for $\Usigp P$ as
claimed.

\begin{defn}
Let $P$ be a plain operad, and $\phi : \Fp \Phi \to P$ be a generator for $P$.
Let $E$ be a subsignature of $(\Up \Fp \Phi)^2$, so that each $E_n$ is a set of
\nary\ $\Phi$-equations. 
Let $i$ be the inclusion map $\xymatrix{E \incar[r] & (\Up \Fp \Phi)^2 }$, and
$\pi_1, \pi_2$ be the projection maps $(\Up \Fp \Phi)^2 \to \Up \Fp \Phi$.
Then $P$ \defterm{satisfies all equations in $E$} if the diagram
\[
\NoCompileMatrices
\fork{\Fp E}{\trans{\pi_1 i}}{\trans{\pi_2 i}}{\Fp \Phi} \phi P
\]
is a fork.
\end{defn}
We say that a symmetric or \dop\ satisfies a signature of equations if the
analogous condition holds in \SymmOperad\ or \DOpCat.

Recall the notion of the ``kernel pair'' of a morphism:

\begin{defn}
\index{kernel pair}
Let $f: A \to B$ in some category $\C$.
The \defterm{kernel pair} of $f$ is the pair $\parallelpair W p q A$ of maps in
the pullback square
\[
\NoCompileMatrices
\commsquare {W \pullback} p A q f A f B
\]
if this pullback exists.
\end{defn}

\def\srtrees#1{\Fp\Up#1}
\begin{lemma}
\label{lem:sreqns}
Let $\epsilon$ be the counit of the adjunction $\Fp \dashv \Up$.
Let $\parallelpair Q {\pi_1} {\pi_2} {\Fp{(\Usig P)}}$ be the kernel
pair of the component
\[
\epsilon_{\Usigp P} : \Fp \Usig P = \Fp \Up \Usigp P \to \Usigp P,
\]
of $\epsilon$.
Let $h$ be the unique map $Q \to (\Fp\Usig P)^2$ induced by $\pi_1, \pi_2$.
Then the image of $\Up h$ is the largest signature of plain-operadic $\Usig
P$-equations satisfied by $P$.
\end{lemma}
\begin{proof}
$Q, \pi_1, \pi_2$ are given by the diagram
\[
\NoCompileMatrices
\commsquare {Q \pullback} {\pi_1} {\Fp\Usig P}
	 {\pi_2} \epsilon
	 {\Fp\Usig P} \epsilon {\Usigp P}
\]
As a right adjoint, $\Up$ preserves pullbacks; we take the standard
construction of pullbacks in $\Set^\natural$ as subobjects of products, in
which case $h$ is an inclusion map.
An element of $(\Up Q)_n$ is then a pair $(e_1, e_2)$ of \nary\ strongly regular
$\Usig P$ trees such that $\epsilon(e_1) = \epsilon(e_2)$.
Hence, $Q$ is a signature of plain-operadic $\Usig P$-equations satisfied by
$P$.
Conversely, let $E$ be a signature of plain-operadic $\Usig P$-equations
satisfied by $P$, and let $(e_1, e_2)$ be an element of $E_n$:
then $\epsilon(e_1) = \epsilon(e_2)$ and so $(e_1, e_2)$ is an element of $(\Up
Q)_n$.
\end{proof}

\def\pibar#1{\overline{\Up \pi_{#1}}}
\begin{corollary}
\label{cor:srpres}
Let $R$ be the plain operad generated by $\Usig P$, satisfying exactly those
plain-operadic equations satisfied by $P$.
Then
\[
\xymatrixcolsep{4pc}
\parallelpair{\Fp\Up Q}{\pibar 1}{\pibar 2}{\Fp \Usig P}
\]
is a presentation for $R$, where the overbars refer to transposition with
respect to the adjunction $\Fp \dashv \Up$.
\end{corollary}
We recall some standard results.

\begin{lemma}
\label{lem:mon=>regep}
The counit of a monadic adjunction is componentwise regular epi.
\end{lemma}
\begin{proof}
See \cite{acc} 20.15.
\end{proof}

\begin{lemma}
\label{lem:fork.epi}
If $\fork X f g Y h Z$ is a coequalizer in some category, and if $e: W \to X$
is epi, then $\fork W {fe} {ge} Y h Z$ is a coequalizer.
\end{lemma}
\begin{proof}
Suppose $\fork W {fe} {ge} Y i A$ is a fork.
Then $ife = ige$, so $if = ig$ since $e$ is epi.
So $\fork X f g Y i A$ is a fork, and hence $i$ factors uniquely through $h$.
\end{proof}

\begin{lemma}
\label{lem:coeq-kerpair}
In categories with all kernel pairs, every regular epi is the coequalizer of
its kernel pair.
\end{lemma}
\begin{proof}
Let $\C$ have all kernel pairs, and $\fork A f g B e C$ be a coequalizer
diagram in $\C$.
Let $\parallelpair W p q B$ be the kernel pair of $e$.
We will show that $\fork W p q B e C$ is a coequalizer diagram.
Since $ef = eg$, we may uniquely factor $(f,g)$ through W:
\[
\xymatrix{
A \ar[ddr]_g \ar[rrd]^f \unar[dr]|i \\
& W \pullback \ar[r]_p \ar[d]^q
& B \ar[d]_e \\
& B \ar[r]^e
& C
}
\]
Suppose $\fork W p q B h D$ is a fork.
$hp = hq$, so $hpi = hqi$, so $hf = hg$.
By the universal property of $e$, we may factor $h$ uniquely through $e$.
So $\fork W p q B e C$ is a coequalizer diagram, as required.
\end{proof}

\def\epUP{\epsilon_{\Usigp P}}
\begin{lemma}
Let $P$ be a symmetric operad, and let $Q, \pi_1, \pi_2$ be as in Lemma
\ref{lem:sreqns}.
Then the coequalizer of the diagram
\[
\parallelpair Q {\pi_1} {\pi_1} {\Fp \Usig P}
\]
is $\Usigp P$.
\end{lemma}
\begin{proof}
Let $\epsilon'$ be the unit of the adjunction $\Fsigp \dashv \Usigp$.
This adjunction is monadic, so $\epUP$ is regular epi by Lemma
\ref{lem:mon=>regep}.
By Lemma \ref{lem:coeq-kerpair}, $\epUP$ is the coequalizer of its kernel pair, i.e.
\[
\NoCompileMatrices
\fork Q {\pi_1} {\pi_2} {\allterms{\Usig P}} {\epUP} {\Usigp P}
\]
is a coequalizer diagram.
\end{proof}

\begin{theorem}
\index{\Usigp!explicit construction}
Let $P$ be a symmetric operad.
Then $\Usigp P$ is the plain operad whose operations are those in $P$,
satisfying exactly those plain-operadic equations which are true in $P$.
\end{theorem}
\begin{proof}
The adjunction $\Fp \dashv \Up$ is monadic, so if $\epsilon'$ is its
counit, then $\epsilon'_Q : \Fp\Up Q \to Q$ is (regular) epi by Lemma
\ref{lem:mon=>regep}.
Hence, by Lemma \ref{lem:fork.epi},
\def\pieps#1{\pi_{#1} \epsilon'_Q}
\[
\fork {\Fp\Up Q} {\pieps 1} {\pieps 2} {\Fp\Usig P} {\epUP} {\Usigp P}
\]
is a coequalizer.
But $\pieps 1 = \overline{\Up \pi_1}$, and similarly
$\pieps 2 = \overline{\Up \pi_2}$.
Hence, by Corollary \ref{cor:srpres}, $\Usigp P$ is the plain operad generated by $\Usig P$, satisfying all plain-operadic equations true in $P$.
Since $\Usig P = \Up \Usigp P$, the \nary\ operations of $\Usigp P$ are exactly
the \nary\ operations of $P$.
\end{proof}

We may generalize this as follows:

\begin{theorem}
Let $\C$ be a category with pullbacks, and $T$ be a monad on $\C$.
Let $(\halg T A a) \in \C^T$.
Let $\parallelpair E {\phi_1} {\phi_2} {TA}$ be the kernel pair of $a$ in $\C$.
Then
\[
\fork {F_T E} {\bar \phi_1} {\bar \phi_2} {F_T A} {\epsilon_{(A,a)}} {(A,a)}
\]
is a coequalizer in $\C^T$, where $F_T : \C \to \C^T$ is the free functor, and $\epsilon$ is the counit of the adjunction $F_T \dashv U_T$.
\end{theorem}
\begin{proof}
As above.
\end{proof}

\begin{corollary}
\index{\Ufpsig!explicit construction}
Let $P$ be a \dop.
Then $\Ufpsig P$ is the symmetric operad whose operations are given by those of
$P$, satisfying all linear equations that are true in $P$, and $\Ufpp P$ is the
plain operad whose operations are given by those in $P$, satisfying all
strongly regular equations that are true in $P$.
\end{corollary}

\section{Operads and syntactic classes of theories}
\label{sec:synclass}

We have defined notions of algebras for plain, symmetric and \dops.
We might ask how these are related to the algebraic theories of Chapter
\ref{ch:theories}: are the algebras for an operad $P$ algebras for some
algebraic theory $\T_P$?
If so, what can we say about the theories that so arise?

We will show the following:
\begin{itemize}
\item Plain operads are equivalent in expressive power to \emph{strongly
regular} theories.
\item Symmetric operads are equivalent in expressive power to \emph{linear}
theories.
\item \Dops\ are equivalent in expressive power to general algebraic theories.
\end{itemize}

The first equivalence is proved in \cite{hohc}.
The second has long been folklore (see, for instance, \cite{baez_ualg} page
50), but as far as I know no proof has appeared before.
An (independently found) proof does appear in an unpublished paper of Ad\'amek
and Velebil, who also consider the enriched case.
\index{Ad\'amek, Ji\v r\'i} \index{Velebil, Ji\v r\'i}
The third equivalence was proved in two stages by Tronin, in \cite{tronin} and
\cite{tronin2}.

Recall the definitions of strongly regular and linear terms from
Definition \ref{def:sr_linterms}, and the definitions of strongly regular,
permuted and finite product trees (Definitions \ref{def:srterm} and
\ref{def:permtree}).

Let $\Phi$ be a signature.
We will show that there is an isomorphism between the set $(\Tfp \Phi)_n$ and
the set of \nary\ words in $\Phi$, and that this isomorphism restricts to
further isomorphisms as follows:

\begin{equation}
\label{eq:termtypes}
\xymatrix{
(\Tfp \Phi)_n            & \cong
& \{\mbox{\nary\ words in $\Phi$}\} \\
(\Tsig \Phi)_n \incar[u] & \cong
& \{\mbox{\nary\ linear words in $\Phi$}\} \incar[u] \\
(\Tp \Phi)_n \incar[u] & \cong
& \{\mbox{\nary\ strongly regular words in $\Phi$}\} \incar[u] \\
}
\end{equation}
The maps in the left-hand column can be viewed as inclusions between different
sets of finite product trees, or equivalently as maps arising from the units of
the adjunctions $\Ffpsig \dashv \Ufpsig$ and $\Fsigp \dashv \Usigp$.

Let $\Phi$ be a signature.
Observe that trees in $\Phi$ give rise to terms according to the following
recursive algorithm:
\begin{itemize}
\item Let $\tau$ be an \nary\ strongly regular tree, and  $Y = (y_1, y_2,
\dots, y_n)$ a finite sequence of variables.
The \defterm{term $\term(\tau,Y)$ arising from $\tau$ with working alphabet
$Y$} is given as follows:
\begin{itemize}
\item If $\tau = |$, then $\term(\tau,Y) = y_1$.
\item If $\tau = \phi\ocomp(\tau_1, \dots, \tau_n)$, then
\[
\term(\tau,Y) =
\phi\tcomp(\term(\tau_1, (y_1, \dots, y_{i_1})), \dots, \term(\tau_n, (y_{1 +
i_{n-1}}, \dots, y_{i_n}))),
\]
where $i_1$ is the arity of $\tau_1$, and $i_j - i_{j-1}$ is the arity of
$\tau_j$ for $j > 1$.

\end{itemize}
\item The \defterm{term $\term(\tau)$ arising from $\tau$} is $\term(\tau,
(x_1, x_2, \dots, x_n))$.
\item Let $\sigma\cdot \tau$ be a permuted tree.
Then $\term(\sigma\cdot \tau) = \term(\tau, (x_{\sigma 1}, x_{\sigma 2},
\dots, x_{\sigma n}))$.
\item Let $f \cdot \tau$ be a finite product tree.
Then $\term(f\cdot \tau) = \term(\tau, (x_{f(1)}, x_{f(2)}, \dots, x_{f(n)})$
\end{itemize}

\begin{defn}
Let $t$ be a $\Phi$-term.
We define a plain-operadic tree $\tree(t)$ recursively:
\begin{itemize}
\item if $t$ is a variable, let $\tree(t) = |$.
\item if $t = \phi \tcomp (t_1, \dots, t_n)$, let
$\tree(t) = \phi \ocomp (\tree(t_1), \dots, \tree(t_n))$.
\end{itemize}
\end{defn}

\begin{lemma}
\label{lem:term<->tree}
\index{tree!finite product!equivalence to terms}
Every $\Phi$-term $t$ is equal to $\term(f \cdot \tau)$ for a unique finite
product tree $(f\cdot \tau)$.
\end{lemma}
\begin{proof}
We will show
\begin{enumerate}
\item
\label{term->tree->term}
if $t$ is a $\Phi$-term, then $\term(\munge(t)\cdot \tree(t)) = t$;
\item
\label{tree->term->tree}
if $(f \cdot \tau)$ is a finite product tree, then
$f = \munge(\term(f \cdot \tau))$ and $\tau = \tree(\term(f \cdot \tau))$.
\end{enumerate}

(\ref{term->tree->term}) Let $t$ be a $\Phi$-term.
Let $f = \munge(t)$, and $\tau = \tree(t)$.
Then $\term(f \cdot \tau)$ is $\term(\tau, (x_{f(1)}, \dots, x_{f(n)}))$.
We proceed by induction.
\begin{itemize}
\item if $t = x_i$, then $\term(f\cdot \tau)$ is $\term(|, (x_i))$, which is
$x_i$.
\item if $t = \phi \tcomp (t_1, \dots, t_n)$, where each $t_i$ has arity $k_i$,
then
\begin{eqnarray*}
\term(f \cdot \tau) &= &\term(\phi \ocomp (\tree(t_1), \dots, \tree(t_n)),
	(x_{f(1)}, \dots, x_{f(n)})) \\
& = & \phi \tcomp (\term(\tree(t_1), (x_{f(1)}, \dots, x_{f(k_1)})), \dots, \\
& & \phantom{\phi \tcomp (} \term(\tree(t_n), (x_{f((\sum_{i=1}^{n-1} k_i) + 1)},
		\dots, x_{f(\sum_{i=1}^n k_i)}))) \\
& = & \phi \tcomp (t_1, \dots, t_n) \\
& = & t.
\end{eqnarray*}
\end{itemize}

(\ref{tree->term->tree})
Let $\tau$ be an \nary\ plain-operadic tree in $\Phi$, and $f$ a function of
finite sets with codomain $\m$.
Let $t = \term(f \cdot \tau)$.
We proceed as usual by induction on $\tau$.
\begin{itemize}
\item If $\tau = |$, then $t = x_{f(1)}$; then $\tree(t) = | = \tau$ and
$\munge(t)$ is the function $\fs 1 \to \m$ sending 1 to $f(1)$, i.e. $\munge(t)
= f$.
\item If $\tau = \phi \ocomp (\tau_1, \dots, \tau_n)$, then
\begin{eqnarray*}
t &= &\term(\phi \ocomp (\tau_1, \dots, \tau_n), (x_{f(1)},
	\dots, x_{f(\sum k_i)})) \\
&= &\phi \tcomp(\term(\tau_1,(x_{f(1)}, \dots, x_{f(k_1)})), \dots,
	\term(\tau_n, (x_{f((\sum_{i=1}^{n-1} k_i) + 1)}) \dots, 
		x_{f(\sum_{i=1}^n k_i)})))
\end{eqnarray*}
By induction, $\var(t) = (x_{f(1)}, \dots, x_{f(\sum k_i)})$, so $\munge(t) =
f$, and $\tree(t) = \phi \ocomp(\tau_1, \dots, \tau_n) = \tau$ as required.
\end{itemize}
\end{proof}

We have now established the isomorphism in the top line of \ref{eq:termtypes}.
If we use this isomorphism to identify \dops\ with finitary monads on \Set, we
may view the functor $\Ffpp$ as the well-known functor sending a plain operad
to its associated monad on \Set.

\begin{lemma}
\label{lem:syntactic_semantic_equiv}
\index{tree!permuted!equivalence to linear terms}
\index{tree!plain-operadic!equivalence to strongly regular terms}
Let $t$ be a $\Phi$-term.
Then $t$ is linear iff $t = \term(\tau)$ for some permuted tree $\tau$, and
strongly regular iff $t = \term(\tau)$ for some strongly regular tree $\tau$.
\end{lemma}
\begin{proof}
In Lemma \ref{lem:term<->tree}, we factored every $\Phi$-term $t$ into a
strongly regular tree $\tree(t)$ and a labelling function $\munge(t)$.
By definition, $t$ is linear iff $\munge(t)$ is a bijection, which occurs iff
$t = \term(\sigma \cdot \tau)$ for some plain-operadic tree $\tau$ and some
bijection $\sigma$.
Hence, the linear terms and permuted trees are in one-to-one correspondence.
Similarly, strongly regular terms and plain-operadic trees are in one-to-one
correspondence.
\end{proof}

The commutativity of \ref{eq:termtypes} now follows from our explicit
construction of $\Fsig$ and $\Fp$ in Section \ref{sec:adj}.

\begin{lemma}
\label{lem:linpresfactor}
Let $(\Phi, E)$ be a presentation of an algebraic theory.
Then $(\Phi, E)$ is linear if and only if the projection maps $\parallelpair E
{\pi_1} {\pi_2} \term(\Phi)$ may be factored through the map
$\xymatrix{\Tsig \Phi \ar[r]^\eta & \Tfp \Phi \ar[r]^{\sim} & \term{\Phi}}$:
\[
\xymatrix{
E \parallelarsdir{\pi_1}{\pi_2}{rr} \uparallelarsdir{}{}{dr}
& & \Tfp \Phi \ar[r]^-{\sim} & \term \Phi \\
& \Tsig \Phi \ar[r]^-{\sim} \ar[ur]_{\eta}
& \{\mbox{\rm linear $\Phi$-terms}\} \incar[ur]
}
\]
\end{lemma}
\begin{proof}
By definition, the presentation is linear iff $\pi_1, \pi_2$ factor
through the signature of linear $\Phi$-terms.
By Lemma \ref{lem:syntactic_semantic_equiv}, this signature is isomorphic to
$\Tsig \Phi$, so we are done.
\end{proof}

\begin{theorem}
Let $Q \in \DOpCat$.
Then
\begin{enumerate}
\item $Q$ is strongly regular iff there exists a $P \in \Operad$ such that $Q \cong \Ffpp P$;
\item $Q$ is linear iff there exists a $P \in \SymmOperad$ such that $Q \cong
\Ffpsig P$;
\end{enumerate}
\end{theorem}
\begin{proof}
We will consider the linear case; the strongly regular case is proved
analogously.

If $Q$ is linear, then there exists a linear presentation $\parallelpair E {}{}
{\Ffp \Phi}$ for $Q$.
We may regard $E$ as a subobject of the signature of $\Phi$-equations.
By assumption, $E$ consists only of linear equations; by Lemma
\ref{lem:syntactic_semantic_equiv}, every $(s,t) \in E$ is $(\term(\sigma_1
\cdot \tau_1), \term(\sigma_2 \cdot \tau_2))$ for some pair $(\sigma_1 \cdot
\tau_1, \sigma_2 \cdot \tau_2)$ of permuted trees.
So the diagram $\parallelpair E {} {} {\Ffp \Phi}$ in \DOpCat\ is the image
under $\Ffpsig$ of a diagram
\[
\parallelpair {E'} {} {} {\Fsig \Phi}
\]
in $\SymmOperad$.
This diagram has a coequalizer: call it $P$.
The functor $\Ffpsig$ is a left adjoint, and thus preserves coequalizers:
hence, $Q$ is the image under $\Ffpsig$ of $P$.

Now suppose $Q = \Ffpsig P$ for some symmetric operad $P$.
We may take the canonical presentation of $P$:
\[
\xymatrixcolsep{4pc}
\fork{\Fsig\Usig\Fsig\Usig P}{\epsilon\Fsig\Usig}{\Fsig\Usig\epsilon}
	{\Fsig \Usig P}\epsilon P
\]
and apply $\Ffpsig$ to it:
\[
\xymatrixcolsep{4pc}
\fork{\Ffp(\Usig\Fsig\Usig P)}{\Ffpsig\epsilon\Fsig\Usig}
	{\Ffp\Usig\epsilon}{\Ffp \Usig P}{\Ffpsig\epsilon}{\Ffpsig P = Q}
\]
Since $\Ffpsig$ is a left adjoint, it preserves coequalizers, so the transpose
of this parallel pair is a presentation for $Q$.
Take this transpose:
\[
\xymatrixcolsep{4pc}
\xymatrix{
	\Usig\Fsig\Usig P
	\parallelars{\overline{\epsilon\Fsig\Usig}}
		    {\overline{\Fsig\Usig\epsilon}}
	& \Usig\Fsig\Usig P \ar[r]^{\eta'}
	& {\Ufp \Ffp \Usig P} \ar[r]^{\Ufp\Ffpsig\epsilon}
	& \Ufp Q
}
\]
where $\eta'$ is the unit of the adjunction $\Ffp \dashv \Ufp$, and the bars
refer to transposition with respect to the adjunction $\Fsig \dashv \Usig$.

This is in precisely the form required for Lemma \ref{lem:linpresfactor}.
\end{proof}

\begin{example}
\index{monoids}
\index{pointed sets}
\index{commutative monoids}
The theories of monoids and pointed sets are strongly regular, because the \dops\ corresponding to these theories are in the image of $\Ffpp$; the theory of
commutative monoids is linear but not strongly regular, because the \dop\ whose
algebras are commutative monoids is in the image of $\Ffpsig$ but not in the
image of $\Ffpp$.
\end{example}

There is a little more to be said about these classes of theories.

\begin{defn}
\index{wide pullback}
A \defterm{wide pullback} is a limit of a (possibly infinite) diagram of the
form
\[
\xymatrixrowsep{1pc}
\xymatrixcolsep{3pc}
\xymatrix{
	\bullet \ar[ddr] \\
	\bullet \ar[dr] \\
	\vdots & \bullet \\
	\bullet \ar[ur] \\
	\vdots
}
\]
\end{defn}

\begin{defn}
\index{cartesian!natural transformation}
A natural transformation $\alpha : F \to G$ is \defterm{cartesian} if every
naturality square
\[
\commsquare {FA} {Ff} {FB} {\alpha_A} {\alpha_B} {GA} {Gf} {GB}
\]
is a pullback square.
\end{defn}

\begin{defn}
\index{cartesian!monad}
A monad $(T, \mu, \eta)$ is \defterm{cartesian} if $T$ preserves pullbacks and
$\mu, \eta$ are cartesian natural transformations.
\end{defn}

\begin{theorem}
\label{thm:opds<=>cart. monad}
A plain operad is equivalent to a cartesian monad on \Set\ equipped with a
cartesian map of monads to the free monoid monad.
\end{theorem}
\begin{proof}
See \cite{hohc} 6.2.4.
Let 1 be the terminal plain operad; algebras for 1 are monoids.
Since 1 is terminal, every plain operad $P$ comes equipped with a map $!: P \to
1$.
This induces a cartesian map of monads $T_! : T_P \to T_1$, and $T_1$ is the
free monoid monad.
\end{proof}

\begin{lemma}
Let $T, S$ be endofunctors on a category $\Acat$, let $\alpha : T \to S$ be a
cartesian natural transformation, and let $S$ preserve wide pullbacks.
Then $T$ preserves wide pullbacks.
\end{lemma}
\begin{proof}
This follows from the facts that wide pullbacks are products in slice
categories and that the functor $f^* : \Acat/B \to \Acat/A$ induced by a map $f
: A \to B$ is product-preserving.
\end{proof}

\begin{corollary}
\label{cor:wpb}
Let $P$ be a plain operad.
Then the functor part of the  monad $T_P$ arising from $P$ preserves wide
pullbacks.
\end{corollary}

\begin{defn}
\index{familial representability}
A functor $F : \C \to \Set$ is \defterm{familially representable} if $F$ is a
coproduct of representable functors.
A monad $(T, \mu, \eta)$ on $\Set$ is \defterm{familially representable} if $T$
is familially representable.
\end{defn}

\begin{theorem} (Carboni-Johnstone)
Let $\C$ be a complete, locally small, well-powered category with a small
cogenerating set, and let $F : \C \to \Set$ be a functor.
The following are equivalent:
\begin{enumerate}
\item $F$ is familially representable;
\item $F$ preserves wide pullbacks.
\end{enumerate}
\end{theorem}
\begin{proof}
See \cite{c+j1}, Theorem 2.6.
\end{proof}

\begin{corollary}
The monad associated to a strongly regular theory is familially representable.
\end{corollary}

However, the inclusion is only one-way: there exist cartesian monads $(T, \mu,
\eta)$ such that $T$ is familially representable but the induced theory
is not strongly regular.
For instance, take the theory of involutive monoids:

\begin{defn}
\index{involutive monoid}
\index{monoid with involution|see{involutive monoid}}
An \defterm{involutive monoid} (or \defterm{monoid with involution}) is a
monoid $(M,.,1)$ equipped with an involution $i : M \to M$,
satisfying $i(a.b) = i(b).i(a)$.
\end{defn}

The theory of involutive monoids is familially representable, but not strongly
regular --- see \cite{c+j2}.

\section{Enriched operads and multicategories}
\label{sec:enriched}

In the previous sections we considered operads $P$ where $P_0, P_1,
\dots \in \Set$, and composition was given by functions.
It is possible to consider operads where $P_0, P_1, \dots$ lie in some other
category; the resulting objects are called \emph{enriched operads}.
Enriched operads have many applications and a rich theory: for instance,
topologists often consider operads enriched in \Top\ or in some category of
vector spaces.
Our treatment here will be brief, sufficient only to set up the definitions of
Chapter \ref{ch:categorification}: for more on enriched operads, see
\cite{mss}.

Throughout this section, let $(\V, \otimes, I, \alpha, \lambda, \rho, \tau)$ be
a symmetric monoidal category.

\begin{defn}
\index{multicategory!plain!enriched}
A \defterm{$\V$-multicategory} $\C$ consists of the following:
\begin{itemize}
\item a collection $\C_0$ of \emph{objects},
\item for all $n \in \natural$ and all $c_1, \dots, c_n,
d \in \C_0$, an object $\C(c_1, \dots, c_n; d) \in \V$ called the
\emph{arrows} from $c_1, \dots, c_n$ to $d$,
\item for all $n, k_1, \dots k_n \in \natural$ and $c_1^1 \dots, c_{k_n}^n,
d_1, \dots, d_n, e \in \C_0$, an arrow in $\V$ called \emph{composition}
\[\ocomp
: \C(d_1, \dots, d_n;e)  \otimes
\C(c^1_1, \dots, c^1_{k_1};d_1) \otimes \dots \otimes
\C(c^n_1, \dots, c^n_{k_n};d_n)
\to
\C(c^1_1, \dots, c^n_{k_n};e)\]
\item for all $c \in \C$, a \emph{unit} $u_c: I \to \C(c;c)$
\end{itemize}
satisfying the following axioms:
\begin{itemize}
\item \emph{Associativity:} For all $b\tseq, c\dseq, d\seq, e \in \C$, the
following diagram commutes:
\[
\xymatrixcolsep{1.5pc}
\xymatrix{
	*\txt{$\C(d\seq; e) \otimes \C(c^1\seq; d_1) \otimes \dots \otimes
	\C(c^n\seq; d_n)$ \\
	$\otimes \C(b^1_{1\bullet}; c^1_1) \otimes
	\dots \otimes \C(b^{n}_{k_n\bullet}; c^n_{k_n})$}
		\ar[r]^-\ocomp
		\ar[dd]^-\ocomp
		&
	\C(d\seq; e) \otimes \C(b^1\ldseq; d_1) \otimes \dots \otimes
	\C(b^n\ldseq; d_n)
		\ar[dd]^-\ocomp
		\\ \\
	\C(c\dseq; e) \otimes \C(b^1_{1\bullet}; c^1_1) \otimes
	\dots \otimes \C(b^{n}_{k_n\bullet}; c^n_{k_n})
		\ar[r]^-\ocomp
		&
	\C(b\tseq ; e)
}
\]
\item \emph{Units:} For all $c\seq, d \in \C$, the following diagram commutes:
\[
\xymatrixcolsep{1.3pc}
\xymatrix{
	\C(c\seq ; d)
		 \ar[r]^\lambda
		 \ar[ddrr]^1
		 \ar[d]_{\rho^n}
		 &
	I \otimes \C(c\seq; d)
		\ar[dr]^{u \otimes 1}
		\\
	\C(c\seq ; d) \otimes I \otimes \dots \otimes I
		\ar[dr]_{1\otimes u^{\otimes n}}
		&
		&
	\C(d;d) \otimes \C(c\seq;d)
		\ar[d]^{\ocomp}
		\\
	&
	\C(c\seq ; d) \otimes \C(c_1;c_1) \otimes \dots \otimes \C(c_n;c_n)
		\ar[r]^-{\ocomp}
		&
	\C(c\seq ; d)
}
\]
\end{itemize}
(We suppress the symmetry maps in \V\ for clarity).
\end{defn}

\begin{defn}
\label{def:symmVmcat}
\index{multicategory!symmetric!enriched}
A \defterm{symmetric \V-multicategory} is a \V-multicategory \C\ and, for every
$n \in \natural$, every $\sigma \in S_n$, and every $a_1, \dots a_n, b \in \C$,
an arrow
\[
\sigma \cdot - : \C(a_1, \dots, a_n; b) \longrightarrow 
\C(a_{\sigma 1}, \dots, a_{\sigma n}; b)
\]
in \V\ such that
\begin{itemize}
\item for each $n \in \natural$ and each $a_1, \dots, a_n, b \in \C$, the arrow
$1_n \act - : \C(a_1, \dots, a_n;b) \to \C(a_1, \dots, a_n;b)$ is the identity
arrow on $\C(a_1, \dots, a_n;b)$,
\item for each $\sigma, \rho \in S_n$,
\[
(\rho \cdot - )(\sigma \cdot -) = (\rho\sigma) \cdot -
\]
\item for each $n, k_1, \dots, k_n \in n$, each $\sigma \in S_n$ and $\rho_i
\in S_{k_i}$ for $i = 1, \dots, n$, and for all $a^1_1, \dots, a^n_{k_n}, b_1,
\dots, b_n, c \in \C$, the diagram
\[
\xymatrixcolsep{-5pc}
\xymatrix{
& \C(b_1, \dots, b_n; c) \otimes \bigotimes_{i=1}^n{\C(a^i_1,\dots, a^i_n; b_i)}
\ar[dl]_{(\sigma\act-) \otimes \bigotimes_{i=1}^n{(\rho_i\act -)}}
\ar[ddr]^\ocomp \\
\C(b_{\sigma 1}, \dots, b_{\sigma n};c)\otimes
	\bigotimes_{i=1}^n{\C(a^i_{\rho_i 1}, \dots, a^i_{\rho_i n};b_i)}
\ar[dd]_{1\otimes \sigma_*} \\
& & \C(a^1_1, \dots, a^n_{k_n}; c) \ar[ddl]^{\sigma \ocomp(\rho_1,\dots,\rho_n)}
\\
\C(b_{\sigma 1}, \dots, b_{\sigma n};c) \otimes \bigotimes_{i=1}^n \C(a^{\sigma
i}_{\rho_{\sigma i} 1}, \dots, a^{\sigma i}_{\rho_{\sigma i} n} ; b_{\sigma i})
\ar[dr]_\ocomp \\
& \C(a^{\sigma 1}_{\rho_{\sigma 1} 1}, \dots, a^{\sigma n}_{\rho_{\sigma n}
k_n}; c)
}
\]
commutes, where $\sigma\ocomp(\listn{\rho_#1})$ is as defined in Example
\ref{ex:symmopd}.
\end{itemize}
\end{defn}
In the case $\V = \Set$ (with the cartesian monoidal structure), this is
equivalent to Definition \ref{def:symm_mcat}.

Let \F\ be a skeleton of the category of finite sets and functions, with
objects the sets \fs 0, \fs 1, \fs 2, \dots, where $\fs n = \{1,2,\dots,n\}$.
\begin{defn}
\index{multicategory!finite product!enriched}
A \defterm{\dvm} is
\begin{itemize}
\item A plain \V-multicategory \C;
\item for every function $f : \n \to \m$ in \F, and for all objects $C_1,
\dots, C_n, D \in \C$, a morphism $f\act - : \C(C_1, \dots, C_n;D) \to
\C(C_{f(1)}, \dots, C_{f(n)}; D)$ in \V\,
\end{itemize}
satisfying the conditions given in Definition \ref{def:symmVmcat}, where
$(f\cmp{f_1,\dots, f_n})$ is as given in Definition \ref{dopdef}.
\end{defn}
In the case $\V = \Set$, this is equivalent to Definition \ref{dopdef}.

\begin{defn}
\index{operad!enriched}
A (\defterm{plain}, \defterm{symmetric}, \defterm{finite product})
\defterm{\V-operad} is a (plain, symmetric, finite product) \V-multicategory
with only one object.
\end{defn}

\begin{defn}
\index{morphism!of enriched multicategories}
Let $\C$, $\D$ be plain \V-multicategories.
A \defterm{morphism} $F: \C \to \D$ is
\begin{itemize}
\item for each object $C \in \C$, a choice of object $FC \in \D$,
\item for each $n \in \natural$ and all collections of objects $A_1, \dots,
A_n, B \in \C$, an arrow

$ \C(A_1, \dots, A_n;B) \to \D(FA_1, \dots, FA_n; FB) $ in $\V$
\end{itemize}
such that
\begin{itemize}
\item for all $A \in \C$, the diagram
\[
\xymatrix{
 & I \ar[ddl]_u \ar[ddr]^u \\ \\
\C(A;A) \ar[rr]^F & & \D(FA;FA)
}
\]
commutes;
\item for all $n, k_1, \dots, k_n \in \natural$ and all $C, B_1, \dots, B_n, A^1_1, \dots, A^n_{k_n} \in \C$, the diagram
\[
\xymatrixrowsep{4pc}
\xymatrix{
\C(B\seq;C) \otimes \bigotimes_{i=1}^n \C(A^i\seq;B_i)
\ar[r]^-{\ocomp} \ar[d]_{F\otimes \dots \otimes F}
& \C(A\dseq;C) \ar[d]^F \\
\D(FB\seq;FC)\otimes \bigotimes_{i=1}^n \D(F A^i\seq;FB_i)
\ar[r]^-\ocomp
& \D(FA\dseq;FC)
}
\]
commutes.
\end{itemize}
\end{defn}

Suppose that $\V$ is cocomplete.
Let $Q$ be a (plain, symmetric, finite product) $\V$-operad, and $A$ an object of $\V$.
Let $Q \kel A$ denote the coend
\[
\int^{n\in \C} Q_n \times A^n
\]
where \C\ is
\begin{itemize}
\item the discrete category on $\natural$ if $Q$ is a plain operad;
\item a skeleton $\Bcat$ of the category of finite sets and bijections if $Q$
is a symmetric operad;
\index{\ensuremath{\Bcat}}
\item a skeleton $\F$ of the category of finite sets and all functions if $Q$
is a \dop.
\end{itemize}
This notation is taken from Kelly's papers \cite{kelly} and \cite{kellycalc} on
clubs.
\index{Kelly, Max}

The various endomorphism operads defined in Examples \ref{ex:plEnd},
\ref{ex:symmEnd} and \ref{ex:fpEnd} transfer straightforwardly to the
\V-enriched setting.
An \defterm{algebra} for a (plain, symmetric, finite product) \V-operad $P$ in
a (plain, symmetric, finite product) \V-multicategory $\C$ is an object $A \in
\C$ and a morphism $(\hatmap) : P \to \End(A)$ of the appropriate type.
\index{algebra!for an enriched operad}
Equivalently, an algebra for $P$ in $\C$ is an object $A \in \C$ and a morphism
$h: P \kel A \to A$ such that the diagram
\[
\xymatrix{
P \kel P \kel A \ar[r]^{1 \kel h} \ar[d]_{\ocomp} & P \kel A \ar[d]^h \\
P \kel A \ar[r]^h & A
}
\]
commutes, and $h(1_P, -)$ is the identity on $A$.

\begin{remark}
\index{multicategory!internal}
\index{operad!internal}
There is another possibility, that of considering \emph{internal}
multicategories in the category \V, which gives a different notion: now
$\C_0$ is an object in \V\ rather than a collection.
An internal operad in \V\ is an internal multicategory \C\ such that $\C_0$ is
terminal in \V.
We shall not consider internal multicategories or operads further.
\end{remark}

We shall in particular consider the case $\V = \Cat$, and \Cat-operads again
have a simple concrete description:
\begin{lemma}
\index{operad!enriched!in \Cat}
A {(plain) \Cat-operad} $Q$ is a sequence of categories $Q_0, Q_1,
\dots$, a family of composition functors $\ocomp:Q_n \times Q_{k_1}
\times \ldots \times Q_{k_n} \to Q_{\sum k_i}$ and an identity 
$1_Q \in Q_1$, satisfying (strict) functorial versions of the axioms given in
\ref{lem:operad_description}.
\end{lemma}

\begin{lemma}
A {symmetric \Cat-operad} is a plain \Cat-operad $Q$ with a left group
action of each symmetric group $S_n$ on the corresponding category $Q_n$,
strictly satisfying equations as in Definition \ref{def:symmopd}.
\end{lemma}

\begin{lemma}
A {\dco} is a plain \Cat-operad $Q$ equipped with functors $f \act - :
Q_n \to Q_m$ for each function $f: \n \to \m$ of finite sets, strictly
satisfying equations as in Definition \ref{dopdef}.
\end{lemma}

All of these lemmas can be established by a straightforward check of
the definitions.

Just as 2-category theory has a special flavour distinct from the theory of
$\V$-categories in the case $\V=\Cat$, so the theories of $\Cat$-operads and
$\Cat$-multicategories have unique features:

\begin{defn}
\label{laxmorphdef}
\index{morphism!of algebras for a \dco}
Let $Q$ be a \dco, and let $Q \kel A \stackrel\alpha\to A , Q \kel B
\stackrel\beta\to B$ be algebras for $Q$ in $\Cat$. A
\defterm{lax morphism of $Q$-algebras} $A \to B$ consists of a 1-cell $F: A \to
B$ and a 2-cell $\phi: \beta F \to F \alpha$ satisfying the
following conditions:
\begin{equation}
\label{eq:laxmorph1}
\xymatrixrowsep{3pc}
\xymatrix{
	A \ar[r]^F \ar[d]_\eta
	  \ddlowertwocell<-20>_1{=} 
	  \drtwocell\omit{=}
	& B \ar[d]^\eta
	  \dduppertwocell<20>^1{=} \\ 
	Q\kel A \ar[r]^{1\kel F} \ar[d]_\alpha
	  \drtwocell\omit{\phi}
	& Q\kel B \ar[d]^\beta
		& & = & A \rtwocell^F_F{=} & B \\
	A \ar[r]_F & B
}
\end{equation}
\begin{equation}
\label{eq:laxmorph2}
\xymatrixcolsep{1.4pc}
\xymatrixrowsep{3pc}
\xymatrix{
	Q \kel Q \kel A \ar[d]_\mu \ar[r]^{1\kel 1\kel F}
	  \drtwocell\omit{=}
	& Q \kel Q \kel B \ar[d]^\mu
	& & & Q \kel Q \kel A \ar[dl]_\mu \ar[r]^{1\kel 1\kel F}
		\ar[d]_{1\kel \alpha}
		\drtwocell\omit{\phi}
	& Q \kel Q \kel B  \ar[dr]^\mu \ar[d]^{1\kel \beta} \\
	Q \kel A \ar[r]^{1\kel F} \ar[d]_{\alpha}
		\drtwocell\omit{\phi}
	& Q \kel B \ar[d]^{\beta}
	& =
	& Q \kel A \ar[dr]_{\alpha}
	  \rtwocell\omit{=}
	& Q \kel A \ar[r]^{1\kel F} \ar[d]_{\alpha}
		\drtwocell\omit{\phi}
	& Q \kel B \ar[d]^{\beta}
	  \rtwocell\omit{=}
	& Q \kel B \ar[dl]^{\beta} \\
	A \ar[r]_F
	& B
	& & & A \ar[r]_F
	& B
}
\end{equation}
\begin{equation}
\label{eq:laxmorph3}
\xymatrixcolsep{1.4pc}
\xymatrixnocompile{
	& Q_m \times A^n \ar[dl]_{1 \times f^*} \ar[dr]^{1 \times F^n}
		\ddtwocell\omit{=}
	& & & & Q_m \times A^n \ar[dl]_{1 \times f^*} \ar[dr]^{1 \times F^n}
		\ar[dd]^{f_* \times 1}
	\\
	Q_m \times A^m \ar[dd]^\alpha \ar[dr]^{1 \times F^m}
	& & Q_m \times B^n \ar[dl]_{1 \times f^*} \ar[dd]^{f_* \times 1}
	& & Q_m \times A^m \ar[dd]^\alpha
		\drtwocell\omit{=}
	& & Q_m \times B^n \ar[dd]^{f_* \times 1} \\
		\drtwocell\omit{\phi}
	& Q_n \times B^m \ar[dd]^\beta
		\drtwocell\omit{=}
	& & =
	& & Q_n \times A^n \ar[dl]_\alpha \ar[dr]_{1 \times F^n}
		\urtwocell\omit{=}
		\ddtwocell\omit{\phi} \\
	A \ar[dr]_F
	& & Q_n \times B^n \ar[dl]^\beta
	& & A \ar[dr]_F
	& & Q_n \times B^n \ar[dl]^\beta \\
	& B
	& & & & B
}
\end{equation}
for all functions $f : \m \to \n$.

A morphism $(F,\phi)$ is \defterm{weak} if $\phi$ is invertible, and
\defterm{strict} if $\phi$ is an identity.
\index{weak!morphism of algebras for a \dco}
\index{strict!morphism of algebras for a \dco}

\index{morphism!of algebras for a symmetric \Cat-operad}
\index{lax!morphism of algebras for a symmetric \Cat-operad}
\index{weak!morphism of algebras for a symmetric \Cat-operad}
\index{strict!morphism of algebras for a symmetric \Cat-operad}
Lax morphisms for algebras of plain \Cat-operads are required to satisfy
\ref{eq:laxmorph1} and \ref{eq:laxmorph2}, and lax morphisms for algebras of
symmetric \Cat-operads are required to satisfy \ref{eq:laxmorph1},
\ref{eq:laxmorph2} and the restriction of \ref{eq:laxmorph3} to the case where
$f$ is a bijection.
\end{defn}

We shall make use of a more explicit formulation in the plain case.
\begin{lemma}
\label{lem:wkfunc_explicit}
Let $Q$ be a plain $\Cat$-operad, and let $(A,h)$ and $(B,h')$ be $Q$-algebras.
A lax map of $Q$-algebras $(A,h) \to (B,h')$ is a pair $(G, \psi)$,
where $G:A \to B$ is a functor and $\psi$ is a sequence of natural
transformations $\psi_i : h'_i (1 \times G^i) \to G h_i$, called the
\defterm{coherence maps}, such that the following equation holds, for all $n,
k_1, \dots, k_n \in \natural$:
\[
\xymatrixrowsep{4pc}
\xymatrixcolsep{4pc}
\xymatrix{
	Q_n \times \prodkn Q \times A^{\sum k_i}
		\ar[d]_{1\times 1^n \times G^{\sum k_i}}
		\ar[r]^-{\prodkn{h}}
		\drtwocell\omit{^*{!(-1,-1.5)\object{\prodkn{\psi}}}}
	& Q_n \times A^n
		\ar[d]|{1\times G^n}
		\ar[r]^{h_n}
		\drtwocell\omit{^*{!(-1,-1.5)\object{\psi_n}}}
	& A
		\ar[d]^G
		\\
	Q_n \times \prodkn Q \times A^{\sum k_i}
		\ar[r]_-{\prodkn{h'}}
	& Q_n \times B^n
		\ar[r]_{h'}
	& B \\
}
\hskip 1in
\]
\begin{equation}
\label{weakfunctordef}
\hskip 2in
\midequals
\xymatrix{
	Q_n \times \prodkn Q \times A^{\sum k_i}
		\ar[d]_{1\times 1^n \times G^{\sum k_i}}
		\ar[r]^-{h_{\sum k_i}}
		\drtwocell\omit{^*{!(-1.5,-2.5)\object{\psi_{\sum k_i}}}}
	& A
		\ar[d]^G
		\\
	Q_n \times \prodkn Q \times B^{\sum k_i}
		\ar[r]_-{h'_{\sum k_i}}
	& B \\
}
\end{equation}
and the diagram
\begin{equation}
\label{weakfunctordef2}
\xymatrixrowsep{3pc}
\xymatrixcolsep{3pc}
\xymatrix{
	G a \ar[r]^{\delta'_1} \ar[d]_1
	& h'(1_P, G a) \ar[d]^{\psi_1} \\
	G a \ar[r]_{G \delta_1}
	& G h(1_P, a)
}
\end{equation}
commutes.
The morphism is weak if every $\psi$ is invertible, and
strict if every $\psi$ is an identity.
\end{lemma}
\begin{proof}
This can be established by a straightforward check of the definition.
\end{proof}

\begin{defn}
\label{transfdef}
\index{transformation!between morphisms of algebras for a \Cat-operad}
Let $Q, A, B$ etc.\ be as above, and let $(F,\phi), (G,\gamma)$ be lax
morphisms of $Q$-algebras $A \to B$.
A \defterm{$Q$-transformation} $F \to G$ is a natural transformation $\sigma:
F \to G$ such that
\begin{equation}
\xymatrixrowsep{1pc}
\xymatrix{
	Q \kel A \rtwocell^{1\kel F}_{1\kel G}{\sigma} \ar[dd]_\alpha
		\ddrtwocell\omit{\gamma}
	& Q \kel B \ar[dd]^\beta
	& & Q \kel A \ar[r]^{1 \kel F} \ar[dd]_\alpha
		\ddrtwocell\omit{\phi}
	& Q \kel B \ar[dd]^\beta \\
	& & = \\
	A \ar[r]_G
	& B
	& & A \rtwocell^F_G{\sigma}
	& B
}
\end{equation}
\end{defn}

\begin{lemma}
\label{inv <=> Q-inv}
A $Q$-transformation $\sigma: (F, \phi) \to (G, \psi)$ is invertible as a
natural transformation if and only if it is invertible as a $Q$-transformation.
\end{lemma}
\begin{proof}
``If'' is obvious: we concentrate on ``only if''. It is enough to show that
$\sigma^{-1}$ is a $Q$-transformation, which is to say that
\begin{equation}
\xymatrix{
	h(q, G\adot)
		\ar[r]^\psi
		\ar[d]_{h(q,\sigma^{-1}_{\adot})} &
	Gh(q,\adot)
		\ar[d]^{\sigma^{-1}_{h(q,\adot)}} \\
	h(q, F\adot)
		\ar[r]^\phi &
	Fh(q,\adot)
}
\end{equation}
commutes for all $(q,\adot) \in Q \kel A$, and this follows
from the fact that ${\sigma_{h(q,\adot)}} \fcomp \phi = \psi \fcomp
{h(q,\sigma_{\adot})}$.

\end{proof}

\Dcos, their morphisms and transformations form a 2-category
called \DCOCat.
Similarly, there is a 2-category $\CatOperad$  of plain \Cat-operads, their
morphisms and transformations, and a 2-category $\CatSymmOperad$, of symmetric
operads, their morphisms and transformations.
\begin{theorem}
There is a chain of monadic adjunctions
\index{adjunctions!for \Cat-operads}
\begin{equation}
\label{eqn:cat_adjs}
\xymatrixrowsep{3pc}
\xymatrix{
  \DCOCat \ar@<1.2ex>[d]^\Ufpsig
	 \ar@/r1.7cm/[dd]^\Ufpp
	 \ar@/r3cm/[ddd]^\Ufp \\
  \CatSymmOperad \ar@<1.2ex>[d]^\Usigp \ar@<1.2ex>[u]^\Ffpsig_{\dashv}
	 \ar@/r1.7cm/[dd]^\Usig \\
  \CatOperad \ar@<1.2ex>[u]^\Fsigp_{\dashv} \ar@<1.2ex>[d]^\Up
	 \ar@/l1.7cm/[uu]^\Ffpp \\
  \Cat^\natural \ar@<1.2ex>[u]^\Fp_{\dashv} \ar@/l1.7cm/[uu]^\Fsig 
  	 \ar@/l3cm/[uuu]^\Ffp
}
\end{equation}
\end{theorem}
\begin{proof}
This follows from Lemmas \ref{lem:frees_exist} and \ref{lem:monadj}, via an
application of the argument of Theorem \ref{thm:opd monadic over setN}.
\end{proof}

Since operads can be considered as one-object multicategories, a \Cat-operad
$P$ (of whatever type) is really a 2-dimensional structure.
We will therefore refer to the objects and morphisms of the categories $P_i$
as \defterm{1-cells} and \defterm{2-cells} of $P$, respectively.
\index{1- and 2-cells}

\section{Maps of algebras as algebras for a multicategory}
\label{sec:Pbar}

Let $P$ be a plain operad.
We form a multicategory $\bar P = \cat{2} \times P$, where $\cat{2}$ is
the category $(\xymatrix{\cdot\ar[r] &\cdot})$.
We may describe $\bar P$ as follows: there are two objects, labelled
0 and 1; the hom-sets $\bar P(0,\dots, 0; 0)$ and $\bar P(x_1, \dots, x_n; 1)$
are copies of $P_n$, for $x_i \in \{0,1\}$, and $\bar P(x_1, \dots, x_n; 0) =
\emptyset$ if any of the $x_i$s are 1.
Composition is given by composition in $P$.
An algebra for $\bar P$ is a pair $A_0, A_1$ of $P$-algebras, and a morphism of $P$-algebras $A_0 \to A_2$.
See \cite{markl}, Example 2.4 for more details.

We can extend this construction by defining a multicategory $\bar{\bar P} =
\cat{3} \times P$, whose algebras are composable pairs of maps of
$P$-algebras, a multicategory $\bar {\bar {\bar P}} = \cat{4} \times P$
whose algebras are composable triples of maps of $P$-algebras, and so on.
With the obvious face and degeneracy maps, these multicategories form a
cosimplicial object in the category of plain multicategories.

The same construction can be performed for symmetric and enriched operads, and
the result continues to hold.

\chapter{Factorization Systems}
\label{ch:factsys}

The theory of factorization systems was introduced by Freyd and Kelly in
\cite{freyd+kelly} (though it was implicit in work of Isbell in the 1950s).
We shall use it in subsequent chapters to define the weakening of an algebraic
theory.
Here, we recall the basic definitions and some relevant theorems.

The material in this chapter is standard, and may be found in (for instance)
\cite{borceux} or \cite{acc}; for an alternative perspective and some more
historical background (as well as the interesting generalization to
\emph{weak} factorization systems), see \cite{k+t}.

\begin{defn}
\index{orthogonal}
Let $e: a \to b$ and $m: c \to d$ be arrows in a category \C.
We say that $e$ is \defterm{left orthogonal} to $m$, written $e \orth m$, if,
for all arrows $f: a \to c$ and $g: b \to d$ such that $mf = ge$, there exists
a unique map $t: b \to c$ such that the following diagram commutes:
\[
\xymatrix{
	a \ar[r]^{\forall f} \ar[d]_e
	& c \ar[d]^m \\
	b \unar[ur]^{\exists! t} \ar[r]_{\forall g}
	& d
}
\]
\end{defn}

\begin{defn} \label{def:FS}
\index{factorization system}
Let \C\ be a category.
A \defterm{factorization system} on \C\ is a pair $(\E, \M)$ of classes of maps
in \C\ such that
\begin{enumerate}
\item \label{axm:factor} for all maps $f$ in \C, there exist $e \in \E$ and $m
\in \M$ such that $f = m  e$;
\item \label{axm:iso_closure} \E\ and \M\ contain all identities, and are
closed under composition with isomorphisms on both sides;
\item \label{axm:orthogonal} $\E \orth \M$, i.e. $e \orth m$ for all $e \in
\E$ and $m \in \M$.
\end{enumerate}
\end{defn}
\begin{example}
\label{ex:epi-mono}
\index{factorization system!epi-mono}
Let $\C = \Set$, $\E$ be the epimorphisms, and $\M$ be the monomorphisms.
Then $(\E,\M)$ is a factorization system.
\end{example}
\begin{example} More generally, let $\C$ be some variety of algebras,
$\E$ be the regular epimorphisms (i.e., the surjections), and $\M$ be the
monomorphisms.
Then $(\E,\M)$ is a factorization system.
\end{example}
\begin{example}
\index{factorization system!bijective on objects/full and faithful!on \Digraph}
\label{ex:digraphFS}
Let $\C = \Digraph$, the category of directed graphs and graph morphisms.
Let $\E$ be the maps bijective on objects, and \M\ be the full and faithful
maps.
Then $(\E,\M)$ is a factorization system.
\end{example}
In deference to Example \ref{ex:epi-mono}, we shall use arrows like
$\xymatrix{{} \booar[r] & {}}$ to denote members of \E\ in commutative
diagrams, and arrows like $\xymatrix{{} \lffar[r] & {}}$ to denote members of
$\M$, for whatever values of \E\ and \M\ happen to be in force at the time.

We will use without proof the following standard properties of factorization
systems:
\begin{lemma}
\label{lem:fact_standard}
Let \C\ be a category, and $(\E, \M)$ be a factorization system on \C.
\begin{enumerate}
\item $\E \cap \M$ is the class of isomorphisms in \C.
\item \label{fact_unique} The factorization in \ref{def:FS} (\ref{axm:factor})
is unique up to unique isomorphism.
\item The factorization in \ref{def:FS} (\ref{axm:factor}) is functorial, in
the following sense: if the square
\[
\commsquare A f B g h C {f'} D
\]
commutes, and $f = me, f' = m'e'$, then there is a unique morphism $i$ making
\[
\xymatrix{
	A \ear[r]^{e} \ar[d]_g
	& {} \mar[r]^{m} \unar[d]_i
	& B \ar[d]^h \\
	C \ear[r]^{e'}
	& {} \mar[r]^{m'}
	& D
}
\]
commute.
Thus, given a choice of $e \in \E$ and $m \in \M$ for each $f$ in \C (such that
$f = me$), we may construct functors $\E_*, \M_* : [\cat{2}, \C] \to [\cat{2},
\C]$:
\begin{eqnarray*}
& \E_* : f \mapsto e \\
& \E_* : (g,h) \mapsto (g,i) \\
& \M_* : f \mapsto m \\
& \M_* : (g,h) \mapsto (i,h).
\end{eqnarray*}
These functors are determined by \E\ and \M\ uniquely up to unique isomorphism.
\item \E\ and \M\ are closed under composition.
\item $\E\orthset = \M$ and ${}\orthset \M = \E$, where $\E\orthset = \{f
\mbox{ in } \C : e \orth f \mbox{ for all } e \in \E\}$ and $\orthset\M = \{f
\mbox{ in } \C : f \orth m \mbox{ for all } m \in \M\}$.
\end{enumerate}
\end{lemma}
Proofs of these statements may be found in \cite{acc} section 14.

We will also use the following fact:
\begin{lemma}
\label{lem:monad}
Let \C\ be a category with a factorization system $(\E, \M)$.
Let $T$ be a monad on \C\ and let $\Ebar = \{ f \arin\ \C : U f \in \E\}$
and $\Mbar = \{f \arin\ \C : U f \in \M\}$, where $U$ is the forgetful functor
$\C^T \to \C$.
Then $(\Ebar, \Mbar)$ is a factorization system on $\C^T$ if $T$ preserves
\E-arrows.
\end{lemma}
\begin{proof}
This is established in \cite{acc}, Proposition 20.24: however, we shall provide
a proof for the reader's convenience.
We shall establish the axioms listed in Definition \ref{def:FS}.

\ref{axm:factor}.
Take an algebra map
\[
\algmap f T A a B b
\]
Applying axiom \ref{axm:factor} to the factorization system $(\E, \M)$, we
obtain a decomposition $f = m  e$, where $e: A \to I$ and $m: I \to B$.
We wish to lift this to a decomposition of $f$ as an algebra map.
In other words, we need a map $i: TI \to I$ making the diagram
\[
\xymatrix{
	TA \ear[r]^{Te} \ar[d]_a
	& TI \ar[r]^{Tm} \unar[d]_i
	& TB \ar[d]^b \\
	A \ear[r]^e
	& I \mar[r]^m
	& B
}
\]
commute, such that $(I,i)$ is a $T$-algebra.
Since $T$ preserves \E-arrows, $T e \orth m$, and we may obtain $i$ by
applying this orthogonality to the diagram \[
\xymatrixrowsep{3pc}
\xymatrix{
	TA \ar[r]^{a} \ar[d]_{T e}
	& A \ar[r]^e
	& I \ar[d]^m \\
	TI \ar[r]_{Tm} \unar[urr]^{\exists!i}
	& TB \ar[r]_b
	& B.
}
\]
It remains to show that $(I,i)$ is a $T$-algebra.
For the unit axiom, consider the diagram
\[
\xymatrix{
	A \ear[r]^e \ar[d]^{\eta_A} \ar@/l1cm/[dd]_1
	& I \mar[r]^m \ar[d]^{\eta_I} & B \ar[d]^{\eta_B} \ar@/r1cm/[dd]^1\\
	TA \ear[r]^{Te} \ar[d]_a & TI \ar[r]^{Tm} \ar[d]^i & TB \ar[d]^b \\
	A \ear[r]_e & E \mar[r]_m & B
}
\]
The top squares commute by naturality, and the outside triangles commute since $(A,a)$ and $(B,b)$ are $T$-algebras.
Hence the diagram
\[
\xymatrix{
	& I \mar[dr]^m \unar[dd] \\
	A \ear[ur]^e \ear[dr]_e & & B \\
	& I \mar[ur]_m
}
\]
commutes if the dotted arrow is either $1_I$ or $i\eta_I$.
By orthogonality, $i\eta_I = 1_I$.

For the multiplication axiom, observe that the diagrams
\[
\xymatrixcolsep{3pc}
\xymatrixrowsep{3pc}
\xymatrix{
	T^2 A \ear[r]^{T^2 e} \ar[d]_{\mu_A}
	& T^2 I \ar[r]^{T^2 m} \ar[d]^{\mu_I}
	& T^2 B \ar[d]^{\mu_B} \\
	TA \ear[r]^{Te} \ar[d]_a & TI \ar[r]^{Tm} \ar[d]^i & TB \ar[d]^b \\
	A \ear[r]_e & E \mar[r]_m & B
}
\hskip 1cm
\xymatrix{
	T^2 A \ear[r]^{T^2 e} \ar[d]_{T a} \ar@/l1cm/[dd]_{a\mu_A}
	& T^2 I \ar[r]^{T^2 m} \ar[d]^{T i}
	& T^2 B \ar[d]^{T b} \ar@/r1cm/[dd]^{b\mu_B} \\
	TA \ear[r]^{Te} \ar[d]_a & TI \ar[r]^{Tm} \ar[d]^i & TB \ar[d]^b \\
	A \ear[r]_e & E \mar[r]_m & B
}
\]
both commute.
So the diagram
\[
\xymatrixcolsep{3pc}
\xymatrixrowsep{3pc}
\xymatrix{
	T^2 A \ar[r]^{\mu_A} \ear[d]_{T^2 e}
		& TA \ar[r]^a & A \ar[r]^e & I \mar[d]^m \\
	T^2 I \ar[r]_{T^2 m} \unar[urrr]
		& T^2 B \ar[r]_{\mu_B} & TB \ar[r]_b & B
}
\]
commutes if we take the dotted arrow to be either $i\mu_I$ or $i(Ti)$.
By orthogonality, $i\mu_I = i(Ti)$.

\ref{axm:iso_closure}.  
The image under $U$ of an isomorphism in $\C^T$ is an isomorphism in \C.
The class \E\ contains all isomorphisms in \C, so $\Ebar = U^{-1}(\E)$ contains
all isomorphisms in $\C^T$.
By similar reasoning, $\Ebar$ is closed under composition with isomorphisms,
and $\Mbar$ also satisfies these conditions.

\ref{axm:orthogonal}.
We wish to show that $\Ebar \orth \Mbar$.
Take $T$-algebras
\[
\valg T A a,
\valg T B b,
\valg T I i,
\valg T J j
\]
and algebra maps
\[
\algmap e T A a I i,
\algmap m T J j B b,
\algmap f T A a J j,
\algmap g T I i B b
\]
where the first two maps are in $\Ebar$ and $\Mbar$ respectively.
Suppose that $ge = mf$.
Now, $e \in \E$ and $m \in \M$, so $e \orth m$, and there is a unique map $t$
in \C\ such that
\[
\commsquare A f J e m {I \unar[ur]^{\exists! t}} g B
\]
commutes.
We wish to show that $t$ is a map of $T$-algebras.

Consider the diagram
\[
\xymatrix{
	TA \ar[r]^{Te} \ar[d]_a \ar@/u0.7cm/[rr]^{Tf}
	&TI \ar[r]^{Tt} \ar[d]_i \ar@/u0.7cm/[rr]^{Tg}
	& TJ \ar[r]^{Tm} \ar[d]^j
	& TB \ar[d]^b \\
	A \ar[r]^e \ar@/d0.7cm/[rr]_f
	& I \ar[r]^t \ar@/d0.7cm/[rr]_g
	& J \ar[r]^m
	& B
}
\]
We wish to show that the middle square commutes: the assumptions tell us that
all other squares commute.
Recall that $T e \orth m$, and apply orthogonality to the square
\[
\commsquare{TA}{j Tf} J {Te} m {TI \unar[ur]^{\exists!u}} {gi} B
\]
Now,
\[
\begin{array}{rcll}
j(Tf)	&= & fa		& \mbox{($f$ is a map of $T$-algebras)} \\
	&= & tea	& \mbox{(Definition of $t$)} \\
	&= & ti(Te)	& \mbox{($e$ is a map of $T$-algebras)} \\
\end{array}
\]
and $mti = gi$ by definition of $t$, so $ti = u$ by uniqueness.
Similarly,
\[
\begin{array}{rcll}
gi	&= & b(Tg)	& \mbox{($g$ is a map of $T$-algebras)} \\
	&= & b(Tm)(Tt)	& \mbox{(Definition of $t$)} \\
	&= & mj(Tt)	& \mbox{($m$ is a map of $T$-algebras)} \\
\end{array}
\]
and $j(Tf) = j(Tt)(Te)$ by definition of $t$, so $j(Tt) = u$ by uniqueness.
Hence $j(Tt) = ti$, and $t$ is a map of $T$-algebras.

By construction, $t$ is unique.
So $e \orth m$ in $\C^T$, so $\Ebar \orth \Mbar$.
All the axioms are satisfied, and so $(\Ebar, \Mbar)$ is a factorization
system on $\C^T$.
\end{proof}
\begin{example} 
\index{factorization system!bijective on objects/full and faithful!on \Cat}
Let $(\E, \M)$ be the factorization system on \Digraph\ described in Example
\ref{ex:digraphFS} above, and let $T$ be the free category monad.
\Cat\ is monadic over \Digraph, and $T$ preserves the property of being
bijective on objects.
Hence, this gives a factorization system $(\Ebar, \Mbar)$ on
\Cat\ where $\Ebar$ is the collection of bijective-on-objects functors, and
$\Mbar$ is the collection of full and faithful functors.
\end{example}
\begin{example}
Similarly, there is a factorization system on $\Digraph^\natural$, where \E\ is
the class of maps that are pointwise bijective on objects, and \M\ is the
class of maps that are pointwise full and faithful.
This lifts to a factorization system $(\Ebar, \Mbar)$ on $\Cat^\natural$, in
which $\Ebar$ is the class of pointwise bijective-on-objects arrows, and
$\Mbar$ is the class of pointwise full-and-faithful arrows.
\end{example}
\begin{example} 
\label{ex:catopdFS}
\index{factorization system!bijective on objects/full and faithful!on \CatOperad}
Let $\C = \Cat^\natural$, \E\ be the pointwise bijective-on-objects maps, and
\M\ be those that are pointwise full and faithful.
Since \CatOperad\ is monadic over $\Cat^\natural$ and the monad preserves
bijective-on-objects maps, this gives a factorization system $(\Ebar, \Mbar)$
on \CatOperad\ where $\Ebar$ is the class of levelwise bijective-on-objects
maps, and $\Mbar$ is the class of levelwise full and faithful ones.
Similarly, there is a factorization system $(\Ebar', \Mbar')$ on
\CatSymmOperad\, where  $\Ebar'$ is the class of bijective-on-objects maps, and
$\Mbar'$ is the class of levelwise full and faithful ones.
\end{example}

We shall need one final piece of background:
{
\def\SetX{\ensuremath{\Set^X}}
\def\SetXT{\ensuremath{(\SetX)^T}}
\begin{theorem}
\label{thm:reg epis split in Set^X^T}
If $X$ is a set and $T$ is a monad on \SetX\ then the regular epis in
\SetXT\ are the pointwise surjections.
In other words, the forgetful functor $U:\SetXT \to \SetX$ preserves and
reflects regular epis.
\end{theorem}
\begin{proof}
See again \cite{acc} section 20, in particular Definition 20.21 and Proposition
20.30.
\end{proof}
}

\chapter{Categorification}
\label{ch:categorification}

\section{Desiderata}

Many categorifications of individual theories have been proposed in the
literature.
We aim to replace these with a general definition, which should satisfy the
following criteria insofar as possible:
\begin{itemize}
\item \textbf{Broad:} it should cover as large a class of theories as
possible.
\item \textbf{Consistent with earlier work:} where a categorification of a
given theory is known, ours should agree with this categorification or be
demonstrably better in some way.
\item \textbf{Canonical:} it should be free of arbitrary tunable parameters
(and if possible should be given by some universal property).
\end{itemize}

We shall return to these criteria in Section \ref{sec:evaluation} and evaluate
how close we have come to achieving them.

Our strategy is as follows: we start with a na\"ive version of categorification
for strongly regular theories, which closely parallels Mac Lane and Benabou's
categorification of the theory of monoids.
This will be an \emph{unbiased} categorification, which treats all operations 
equally, without regarding any as ``primitive'': for instance, if
$P$ is the terminal operad (whose strict algebras are monoids), then the weak
$P$-categories will have tensor products of all arities, not just 0 and 2. 
We then re-express our definition of categorification in terms of factorization
systems, which allows us to generalize our definition in two directions
simultaneously: to symmetric operads, and to operads with presentations.
We then use this new definition to recover the classical theory of symmetric
monoidal categories (at which several other proposed general definitions of
categorification fail), and investigate what it yields in the case of some
other linear theories.

\section{Categorification of strongly regular theories}
\label{sec:naivecatn}

The idea is to consider the strict models of our theory as algebras for an
operad, then to obtain the weak models as (strict) algebras for a weakened
version of that operad (which will be a \Cat-operad).
We weaken the operad using a similar approach to that used in Penon's
definition of $n$-category, as described in \cite{penon}.
A non-rigorous summary of Penon's construction can be found in
\cite{cheng+lauda}.

Throughout this section, let $P$ be a plain (\Set-)operad.

Let $D_*: \Operad \to \CatOperad$ be the functor which
takes discrete categories levelwise; i.e., $(D_*P)_n$ is the discrete category
on the set $P_n$.
In terms of the ``$n$-cell'' terminology introduced in Chapter \ref{ch:factsys},
the 1-cells of $(D_* P)_n$ are \nary\ arrows in $P$, and the only 2-cells are
identities.

\begin{defn}
\label{def:weakening}
\index{unbiased!weakening of a plain operad}
\index{weakening!of a plain operad}
The \defterm{unbiased weakening of $P$}, \WkP, is the following \Cat-operad:
\begin{itemize}
\item \emph{1-cells:} 1-cells of $D_* \Fp\Up P$;
\item \emph{2-cells:} if $A, B \in (\Fp\Up P)_n$, there is a single 2-cell $A
\to B$ if $\epsilon(A) = \epsilon(B)$ (where $\epsilon$ is the counit
of the adjunction $\Fp \dashv \Up$), and no 2-cells $A \to B$ otherwise;
\item \emph{Composition of 2-cells:} the composite of two arrows $A \rightarrow
B \rightarrow C$ is the unique arrow $A \rightarrow C$, and in particular, the
arrows $A \rightarrow B$ and $B \rightarrow A$ are inverses;
\item \emph{Operadic composition:} on 1-cells, as in $\Fp\Up P$, and on 2-cells,
determined by the uniqueness property.
\end{itemize}
\end{defn}

\begin{figure}[h]
\centerline{
	\epsfxsize=3in
	\epsfbox{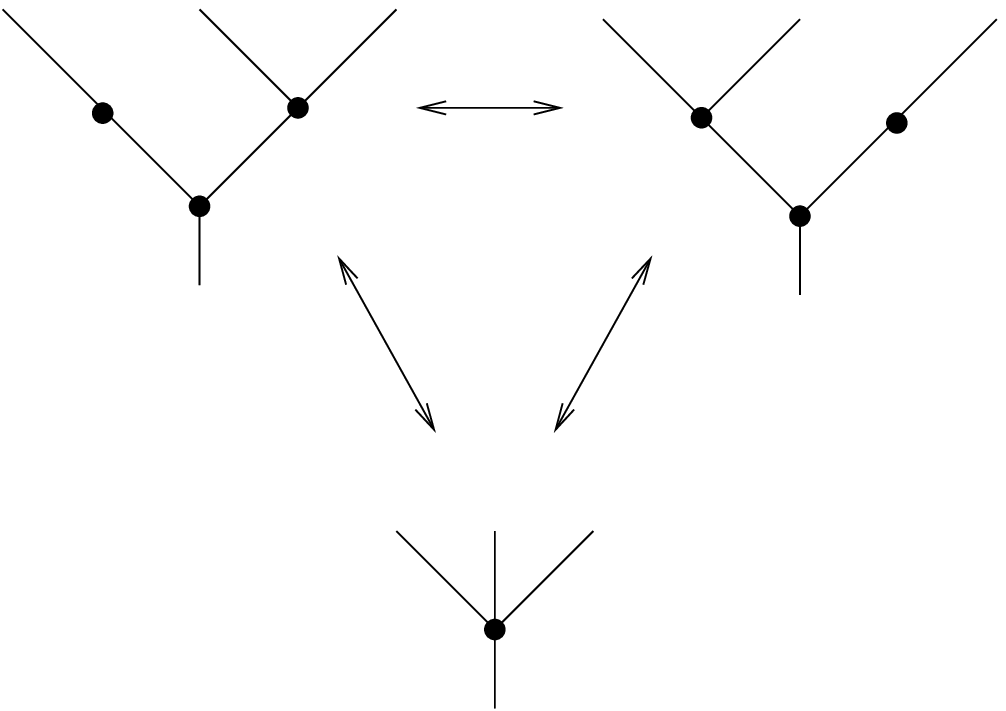}
}
\caption{Part of ${\Wk 1}_3$}
\label{Wk(1)pic}
\end{figure}

See Fig. \ref{Wk(1)pic}, which illustrates a fragment of the unbiased weakening
of the terminal operad 1.
Since $1_n$ is a singleton set for every $n \in \natural$, then $\Wk 1_n$ is
the indiscrete category whose objects are unlabelled \nary\ strongly regular
trees for all $n \in \natural$.
We may embed the discrete category on each $P_n$ in $\WkP_n$, via the map $p
\mapsto p \kel (|, \dots, |)$.
We shall occasionally abuse notation and consider some $p \in P_n$ as a 1-cell
of $\WkP_n$.

\begin{theorem}
\label{thm:WkP_characterization}
$\WkP$ is the unique \Cat-operad with the following properties:
\begin{itemize}
\item $\WkP$ has the same 1-cells as $D_* \Fp\Up P$;
\item we may extend the counit $\epsilon_P : \Fp\Up P \to  P$ to a map of
\Cat-operads $\WkP \to D_* P$, which is full and faithful levelwise.
\end{itemize}
\end{theorem}
\begin{proof}
Immediate.
\end{proof}

We may now make the following definition:

\begin{defn}
\label{def:wkpcat}
\index{weak!$P$-category}
A \defterm{weak $P$-category} is an algebra for \WkP.
\end{defn}

In the case $P = 1$, this reduces exactly to Leinster's definition of unbiased
monoidal category in \cite{hohc} section 3.1.
There, two 1-cells $\phi$ and
$\psi$ have the same image under $\epsilon$ iff they have the same arity, so
the categories $\Wk{1}_i$ are indiscrete.
If $h: \WkP \kel A \to A$ is a weak $P$-category, we refer to the image under
$h$
of a 2-cell $q \to q'$ in $\WkP$ as $\delta_{q,q '}$.
This is clearly a natural transformation $h(q,-) \to h(q',-)$.
As a special case, we write $\delta_q$ for $\delta_{q, \epsilon(q)}$ (where we consider $\epsilon(q)$ as a 1-cell of $\WkP$ as described above).

\begin{defn}
\index{strict!$P$-category}
A \defterm{strict $P$-category} is an algebra for $D_* P$.
\end{defn}
Equivalently, a strict $P$-category is a weak $P$-category in which every
component of $\delta$ is an identity arrow.

\begin{defn}
\index{weak!$P$-functor}
\index{strict!$P$-functor}
Let $(A, h)$ and $(B, h')$ be weak $P$-categories. A \defterm{weak
$P$-functor} from $(A, h)$ to $(B, h')$ is a weak map of $\WkP$-algebras.
A \defterm{strict $P$-functor} from $(A, h)$ to $(B, h')$ is a strict map of
$\WkP$-algebras.
\end{defn}
Equivalently, a strict $P$-functor is a weak $P$-functor for which all the
coherence maps are identities.
These definition are natural generalizations of the definition of weak and
strict unbiased monoidal functors given in \cite{hohc} section 3.1.

\begin{defn}
\index{transformation!between $P$-functors}
Let $(F, \phi)$ and $(G, \psi)$ be weak $P$-functors $(A, h) \to (B, h')$. A
\defterm{$P$-transformation} $\sigma: (F, \phi) \to (G, \psi)$ is a
$\WkP$-transformation $(F, \phi) \to (G, \psi)$, in the sense of Definition
\ref{transfdef}.
\end{defn}

Note that there is only one possible level of strictness here.

There is a 2-category, \WkPCat, whose objects are weak $P$-categories, whose
1-cells are weak $P$-functors, and whose 2-cells are $P$-transformations.
Similarly, there is a 2-category \StrPCat\ of strict $P$-categories, strict
$P$-functors, and $P$-transformations, which can be considered a sub-2-category
of \WkPCat.

\begin{defn}
A \defterm{$P$-equivalence} is an equivalence in the 2-category \WkPCat.
\end{defn}

\begin{lemma}
\label{inv <=> P-inv}
Let $P$ be a plain operad, $(A, h)$ and $(B, h')$ be weak $P$-categories, and
$(F, \phi), (G, \psi) : (A, h) \to (B, h')$ be weak $P$-functors.
A $P$-transformation $\sigma: (F, \phi) \to (G, \psi)$ is invertible as a
$P$-transformation if and only if it is invertible as a natural
transformation.
\end{lemma}
\begin{proof}
This is a straightforward application of Lemma \ref{inv <=> Q-inv}.
\end{proof}

\section{Examples}

Unfortunately, few well-studied theories are strongly regular.
We will consider the following examples:
\begin{enumerate}
\item the trivial theory (in other words, the theory of sets);
\item the theory of pointed sets;
\item the theory of monoids;
\item the theory of $M$-sets, for a monoid $M$.
\end{enumerate}
While we could easily invent a new strongly regular theory to categorify, this
would not help us to see how well our definition of weakening accords with our
intuitions.
Further examples will be considered later, when the machinery to categorify
theories-with-generators and linear theories has been developed.

We will first need to introduce an auxiliary definition:

\begin{defn}
\index{trivial monad}
Let \C\ be a category, and $(T, \mu, \eta)$ be a monad on $\C$.
We say that $(T, \mu, \eta)$ is \defterm{trivial} if $\eta$ is a natural
isomorphism.
\end{defn}

\begin{lemma}
\label{lem:triv<=>initial}
The identity monad on $\C$ is initial in the category $\Mnd(\C)$ of monads on
\C, with the unique morphism of monads $(1_\C, 1, 1) \to (T, \mu, \eta)$ being
$\eta$.
\end{lemma}
\begin{proof}
First we show that $\eta$ is a morphism of monads in the sense of Street 
(Definition \ref{def:monad_morphism}).
One axiom corresponds to the outside of the diagram
\[
\xymatrix{
1 \ar[r]^\eta \ar[d]_1 & T \ar[r]^{\eta T} \ar[dr]^1 \ar[d]_1 & T^2 \ar[d]^\mu\\
1 \ar[r]_\eta & T \ar[r]_1 & T
}
\]
commuting; all the inner segments commute (the top right triangle by the unit
axiom for monads), so the outside must commute.
The other axiom corresponds to the diagram
\[
\xymatrix{
	1 \ar[r]^1 \ar[dr]_\eta & 1 \ar[d]^\eta \\
	& T
}
\]
and this commutes trivially.
Hence, $\eta$ is a morphism of monads $1 \to T$.

Now suppose that $\alpha : 1 \to T$ is a morphism of monads.
From the unit axiom for monad morphisms,  the diagram
\[
\xymatrix{
	1 \ar[r]^1 \ar[dr]_\eta & 1 \ar[d]^\alpha \\
	& T
}
\]
must commute, so $\eta = \alpha$.
\end{proof}

\begin{corollary}
\label{cor:triv<=>initial}
A monad $(T, \mu, \eta)$ on $\C$ is trivial if and only if it is isomorphic to
the identity monad on $\C$.
\end{corollary}
\begin{proof}
If $T$ is isomorphic to the identity monad, then by Lemma
\ref{lem:triv<=>initial} the isomorphism concerned must be $\eta$, so $\eta$
must be invertible.
It is readily checked that if $\eta$ is invertible, then $\eta^{-1}$ must be a
morphism of monads, so if $T$ is trivial then it is isomorphic to the identity
monad.
\end{proof}

\begin{example} \emph{The trivial theory:}
\label{ex:trivial_theory}
\index{sets}
Let $0$ be the initial operad, whose algebras are sets.
An unbiased weak $0$-category is a category equipped with a specified trivial
monad, for the following reason.
$0$ has only one operator (call it $I$), of arity one. Hence $(\Fp\Up 0)_1
\cong \natural$, and all other $(\Fp\Up 0)_n$'s are empty.
All derived operations in the theory of sets are composites of identities, and
thus equivalent to the identity.
So all objects of $\Wk{0}_1$ are isomorphic.
Hence, $I \cong \mbox{id}$.
All diagrams commute: in particular, those giving the monad and monad morphism
axioms commute, so in any weak $0$-category $(\C, (\hatmap))$, the functor
$\hat I$ is a monad, and the isomorphism $\hat I \to 1_\C$ is an isomorphism of
monads.
By Corollary \ref{cor:triv<=>initial}, $\hat I$ must be trivial.

Conversely, suppose $T$ is a trivial monad on a category $\C$.
We wish to show that $\mu$ is also invertible, and thus that $\C$ is an
unbiased weak $0$-category.
From the monad axioms, we have that
\[
\xymatrix{
	T \ar[r]^{\eta T} \ar[dr]_{\mbox{id}}
	& T^2 \ar[d]^{\mu} \\
	& T
}
\]
commutes.
But $\eta T$ is invertible, so $\mu$ must be its inverse.
So $\mbox{id} \cong T \cong T^2 \cong T^3 \cong \ldots$, and all diagrams
commute.
Hence $\C$ is an unbiased weak $0$-category.

If $(\C, S)$ and $(\D, T)$ are weak $0$-categories, then a weak $0$-functor
$(\C, S) \to (\D, T)$ is a functor $F: \C \to \D$ and a natural isomorphism 
$\phi : TF \to FS$, such that the equations 
\[
\xymatrix{
	& \C \ar[dl]_S \ar[dd]^S  \ar[rr]^F \ddrrtwocell\omit{^\phi^{-1}}
	& & \D \ar[dd]^T \\
	\C \ar[dr]_S  \\
	& \C \ar[rr]_F \uuuppertwocell\omit{<-3>\mu} & & D
}
\xymatrix{{} \\ = \\ {}}
\xymatrix{
	& \C \ar[dl]_S \ar[rr]^F \drtwocell\omit{^\phi^{-1}}
	& & \D \ar[dd]^T \ar[dl]_T \\
	\C \ar[dr]_S  \ar[rr]^F 
	& {} \drtwocell\omit{^\phi^{-1}} & D \ar[dr]_T \\
	& \C \ar[rr]_F  & & D \uuuppertwocell\omit{<-3>\mu}
}
\]
and
\[
\xymatrix{
	\C \ddlowertwocell<-10>_{1_\C}{^\eta} \ar[dd]^S \ar[rr]^F
		\ddrrtwocell\omit{^\phi^{-1}}
	& & \D \ar[dd]^T \\ \\
	\C \ar[rr]_F & & \D
}
\xymatrix{{} \\ = \\ {}}
\xymatrix{
	\C  \ddrrtwocell\omit{^\phi^{-1}} \ar[rr]^F \ar[dd]_{1_\C}
	& & \D \dduppertwocell<10>^T{^\eta} \ar[dd]_{1_\D} \\ \\
	\C \ar[rr]_F & & \D
}
\]
are satisfied.
\end{example}

\begin{example} \emph{Pointed sets:}
\label{ex:pointed_sets}
\index{pointed sets}
Let $P$ be the operad with a single element of arity 0 (call it $*$) and a
single element of arity 1 (the identity).
Strict algebras for $P$ in $\Set$ are pointed sets.
The set $(\Fp\Up P)_0$ is countable (it has elements $*, I*, I^2*, I^3*,
\ldots$, and so is $(\Fp\Up P)_1$ (it has elements $\mbox{id}, I, I^2,
\ldots$).
So an unbiased weak $P$-category is a category $\C$ equipped with a
distinguished object $\hat *$ and a trivial monad $\hat I$.

If $(\C, (\hatmap))$ and $(\D, (\barmap)$ are unbiased weak $P$-categories,
then a weak $P$-functor $(\C, (\hatmap)) \to (\D, (\barmap))$ is a triple $(F,
\phi, \psi)$, where $F$ and $\phi$ are as in Example \ref{ex:trivial_theory},
$\psi : \bar * \to F \hat *$ is an isomorphism, and there is exactly one
natural isomorphism $\bar I^n \bar * \to F \hat I^m \hat *$ composed from
$\phi$s and $\psi$s for each $m$ and each $n \in \natural$.
\end{example}

\begin{example} \emph{Monoids:}
\label{ex:monoids}
\index{monoids}
An unbiased weak $1$-category is precisely an unbiased weak monoidal
category in the sense of Definition \ref{def:umoncat}.
An unbiased weak $1$-functor is an unbiased weak monoidal functor.
For a proof, see \cite{hohc} Theorem 3.2.2.
\end{example}

\begin{example} \emph{$M$-sets:}
\label{ex:M-sets}
\index{$M$-sets}
Let $M$ be a monoid, and $N$ be the operad such that
\begin{eqnarray*}
N_1 & = & M \\
N_i & = & \emptyset \mbox{ whenever $i \neq 1$}
\end{eqnarray*}
with composition of arrows of arity 1 given by the multiplication in $M$.
An algebra for $N$ in \Set\ is an $M$-set.
An unbiased weak $N$-category is a category $\C$ with a functor $\hat m :
\C \to \C$ for each $m \in \M$.
For every equation $m_1 m_2 \dots m_i = n_1 n_2 \dots n_j$ that is true in $M$,
there is a natural isomorphism $\delta_{m_1 \dots m_i}^{n_1 \dots n_j} : \hat
m_1 \hat m_2 \dots \hat m_i \to \hat n_1 \hat n_2 \dots \hat n_j$.
If $e$ is the identity element in $M$, then $\hat e$ is a trivial monad.
All diagrams involving these natural isomorphisms commute.
Hence, an unbiased weak $N$-category is a category $\C$ together with a weak
monoidal functor $M \to \End(\C$).
If $(\C, (\hatmap))$ and $(\D, (\barmap))$ are unbiased weak $N$-categories, an
unbiased weak $N$-functor is a functor $F : \C \to \D$ together with natural
transformations $\phi_m : \bar m F \to F \hat m$ for all $m \in M$, such that
if $m_1 m_2 \dots m_i = n_1 n_2 \dots n_j$ in $M$, there is precisely one
natural isomorphism $\bar m_1 \dots \bar m_i F \to F \hat n_1 \dots \hat n_j$
that can be formed by composing $\delta$s and $\phi$s.
\end{example}

\section{A more general approach: factorization systems}

Recall from Definition \ref{def:opd_presentation} the definition of a
presentation and a generator for an operad.
We will define a categorification of any symmetric operad equipped with a
generator, generalizing the unbiased categorification defined in Section
\ref{sec:naivecatn}.
In particular, we shall consider categorification with respect to the component
of the counit $\epsilon_P: \Fsig\Usig P \to P$ at a symmetric operad $P$; this
is a generator for $P$ since both
\[
\fork{\Fsig\Usig \Fsig\Usig P}{\epsilon \Fsig\Usig }{\Fsig\Usig  \epsilon}{\Fsig\Usig P} {\epsilon}{P}
\]
and
\[
\label{item:regular_pres}
\fork{\Fsig\Usig P\times_P \Fsig\Usig P}{\pi_1}{\pi_2} {\Fsig\Usig P}{\epsilon}{P}\]
are coequalizer diagrams (the latter by Lemma \ref{lem:coeq-kerpair}).
We will then show that the categorification is independent of our choice of
generator, in the sense that the symmetric \Cat-operads which arise are
equivalent (and thus have equivalent categories of algebras).
\begin{defn}
\index{weakening!of a symmetric operad}
\label{def:weakening_pres}
Let $\Phi$ be a signature, $P$ be a symmetric operad, and $\phi: \Fsig\Phi \to
P$ be a regular epi in \SymmOperad.
Then the \defterm{weakening} (or \defterm{categorification}) $\Wkwrt P \phi$ of
$P$ with respect to $\phi$ is the (unique-up-to-isomorphism) symmetric
\Cat-operad such that the following diagram commutes:
\[
\xymatrix{
	D_*\Fsig\Phi \ar[rr]^{D_* \phi} \booar[dr]_b
	&& D_* P \\
	& \Wkwrt P \phi \lffar[ur]_f
}
\]
where $f$ is full and faithful levelwise, $b$ is levelwise bijective on
objects, and $D_*$ is the levelwise discrete category functor $\SymmOperad \to
\CatSymmOperad$.
The existence and uniqueness of $\Wkwrt P \phi$ follow from Lemma
\ref{lem:fact_standard} applied to the factorization system on \CatSymmOperad\
described in \ref{ex:catopdFS} above.
\end{defn}

\begin{defn}
\label{def:wkPphicat}
\index{weak!$P$-category}
Let $\phi, \Phi$ and $P$ be as above.
A \defterm{$\phi$-weak $P$-category} is an algebra for $\Wkwrt P \phi$.
\end{defn}
Note that any strict algebra for $P$ can be considered as a $\phi$-weak
$P$-category (for any $\phi$), via the map $\xymatrix{\Wkwrt P \phi \lffar[r] &
D_* P}$.

\begin{defn}
\index{weak!$P$-functor}
Let $\phi, \Phi$ and $P$ be as above.
A \defterm{$\phi$-weak $P$-functor} is a weak map of $\Wkwrt P \phi$-algebras.
\end{defn}

\begin{defn}
\index{weakening!of a symmetric operad}
Let $P$ be a symmetric operad, and $\parallelpair {\Fsig E}{e_1}{e_2}{\Fsig
\Phi}$ be a presentation for $P$, with $\phi: \Fsig \Phi \to P$ being the
regular epi in Definition \ref{def:opd_presentation}.
The \defterm{weakening of $P$ with respect to $(\Phi, E)$} is the weakening of
$P$ with respect to $\phi$.
\end{defn}

\begin{defn}
\label{def:wkp}
\index{unbiased!weakening of a symmetric operad}
The \defterm{unbiased} weakening of $P$ is the weakening arising from
the counit $\epsilon: \Fsig\Usig P \to P$ of the adjunction $\Fsig \dashv
\Usig$.
Call this symmetric \Cat-operad \WkP.
\end{defn}

\begin{lemma}
\label{lem:explicit_wkp}
Let $\phi, \Phi$ and $P$ be as above.
Then, for every $n \in \natural$, the category $\Wkwrt P \phi_n$ is the
equivalence relation $\sim$ on the elements of $(\Fsig \Phi)_n$, where $t_1 \sim t_2$ if $\phi(t_1) = \phi(t_2)$.
\end{lemma}
\begin{proof}
Let $n \in \natural$, and $t_1, t_2 \in \Wkwrt P \phi_n$.
The objects of $\Wkwrt P \phi_n$ are the elements of $(\Fsig \Phi)_n$, by
construction.
Since $\phi_n$ factors through a full functor $\xymatrix{ \Wkwrt P
\phi_n \lffar[r] & (D_*P)_n }$ and $(D_*P)_n$ is the discrete category on
$P_n$, there is an arrow $t_1 \to t_2$ in $\Wkwrt P \phi_n$ iff $\phi(t_1) =
\phi(t_2)$.
Since this functor is also faithful, such an arrow must be unique.
Hence $\Wkwrt P \phi_n$ is a poset; it is readily checked that it is also an
equivalence relation.
\end{proof}

An obvious question is how this notion of weakening is related to the version
defined for plain operads in Section \ref{sec:naivecatn}.
In light of Theorem \ref{thm:WkP_characterization}, it is clear that the
plain-operadic version can be re-phrased as in Definition \ref{def:wkp} above,
but with the factorization occurring in \CatOperad\ rather than \CatSymmOperad.
We may generalize it to give a definition of the weakening of a plain operad
$P$ with respect to a generator $\phi$:

\begin{defn}
\label{def:phi-wkpcat_plain}
Let $P$ be a plain operad, $\Phi$ be a signature, and $\phi: \Fp \Phi \to P$ be
a regular epi.
The \defterm{weakening $\Wkwrt P \phi$ of $P$ with respect to $\phi$} is the
plain $\Cat$-operad given by the bijective on objects/levelwise full and
faithful factorization
\[
\xymatrix{
	D_* \Fp \Phi \ar[rr]^{D_* \phi} \booar[dr] & & D_* P \\
	& \Wkwrt P \phi \lffar[ur]
}
\]
in \CatOperad.
A \defterm{$\phi$-weak $P$-category} is an algebra for $\Wkwrt P \phi$.
\end{defn}

But do the weak algebras for a strongly regular theory $T$ change if we
consider $T$ as a linear theory instead?
We now answer that question in the negative.

\begin{theorem}
\label{thm:plain_as_symmetric}
\index{weakening!of a plain operad!considered as a symmetric operad}
Let $P$ be a plain operad, let $\Phi$ be a signature, and let $\phi : \Fp \Phi
\to P$ be a regular epi.
Then $\Wkwrt{\Fsigp P}{\Fsigp \phi} \cong \Fsigp (\Wkwrt{P}{\phi})$ in the
category $\CatSymmOperad$.
\end{theorem}
\begin{proof}
First note that $\Wkwrt{\Fsigp P}{\Fsigp \phi}$ is well-defined: $\Fsigp$ is a
left adjoint, and hence preserves colimits, so $\Fsigp \phi$ is a regular epi
in $\SymmOperad$.

$\Wkwrt{\Fsigp P}{\Fsigp \phi}$ is defined by its universal property, so it is
enough to show that the \Cat-operad $\Fsigp (\Wkwrt P \phi)$ also has this
property.
Specifically, it is enough to show that if
\[
\xymatrix{
D_* \Fp \Phi \ar[rr]^{D_* \phi} \booar[dr]_{b} && D_* P \\
& \Wkwrt P \phi \lffar[ur]_f
}
\]
is the bijective-on-objects/full-and-faithful factorization of $\phi$, then in
the diagram
\[
\xymatrix{
D_* \Fsig \Phi \ar[rr]^{D_* \Fsigp \phi} \booar[dr]_{\Fsigp b}
&& D_* \Fsigp P \\
& \Fsigp (\Wkwrt P \phi) \lffar[ur]_{\Fsigp f}
}
\]
the arrow $\Fsigp b$ is bijective on objects and the arrow $\Fsigp f$ is
levelwise full and faithful (note that $D_* \Fsigp = \Fsigp D_*$).
But this follows straightforwardly from the explicit construction of $\Fsigp$
in Section \ref{sec:explicit_Fs}.
\end{proof}

\begin{corollary}
Let $P$ be a plain operad, $\phi : \Fp \Phi \to P$ generate $P$, and $A$ be
a $\phi$-weak $P$-category in the sense of Definition
\ref{def:phi-wkpcat_plain}.
Then $A$ is an $\Fsigp \phi$-weak $\Fsigp P$-category in the sense of
Definition \ref{def:wkPphicat}.
Conversely, every $\Fsigp \phi$-weak $\Fsigp P$-category is a weak
$P$-category.
\end{corollary}
\begin{proof}
A $\Fsigp \phi$-weak $\Fsigp P$-category is a category $A$ and a
morphism $\Wkwrt{\Fsigp P}{\Fsigp \phi} \to \End(A)$ of symmetric \Cat-operads.
By Theorem \ref{thm:plain_as_symmetric}, this is equivalent to a morphism
$\Fsigp (\Wkwrt P \phi) \to \End(A)$ in $\CatSymmOperad$, which is equivalent
by the adjunction $\Fsigp \dashv \Usigp$ to a morphism of plain
\Cat-operads $\Wkwrt P \phi \to \Usigp \End(A)$.
This is exactly a $\phi$-weak $P$-category.
\end{proof}

Note that we had to apply $\Fsigp$ to $\phi$ to obtain a generator for $\Fsigp
P$.
This means that the theorem does not tell us that the unbiased categorification
is unaffected by whether we consider our theory as a linear or a strongly
regular one.
In fact, it is not the case that $\Wk{\Fsigp P} \cong \Fsigp (\WkP)$ in
general.

\begin{example}
\index{monoids}
\index{\calS!as $\Fsigp 1$}
Consider the terminal plain operad $1$ whose algebras are monoids.
$\Fsigp 1$ is the operad $\calS$ of Example \ref{ex:symmopd}, for which each
$\calS_n$ is the symmetric group $S_n$.
Then the objects of $\Wk{\calS}_n$ are $n$-leafed permuted trees with each node
labelled by a permutation, whereas the 1-cells of $(\Fsigp \Wk{1})_n$ are
unlabelled permuted trees.
These two sets are not canonically isomorphic.
Hence, there is no canonical isomorphism between $\Wk{\calS}$ and $\Fsigp
\Wk{1}$.
\end{example}

However, we can make a weaker statement: the two candidate unbiased weakenings
are \emph{equivalent} in the 2-category $\CatSymmOperad$.
We shall return to this point in Corollary \ref{cor:unbiased_symm}.

\section{Examples}

\begin{example}
\label{ex:trivial_theory_biased}
\index{sets}
Consider the trivial theory (given by the initial operad 0), with the empty
generating set.
A weak algebra for this theory (with respect to this generating set) is simply
a category.
$\Fsig$ is a left adjoint, and hence preserves colimits, so $\Fsig\emptyset$ is
the initial operad, and the coequalizer $\phi : \Fsig\emptyset \to 0$ is
therefore the identity.
Hence $\Wkwrt 0 \phi$ is also the initial operad, and so a $\phi$-weak
0-category is just a category.
A $\phi$-weak $0$-functor is just a functor.
\end{example}
\begin{example}
\label{ex:ptd_sets}
\index{pointed sets}
Consider the operad $P$ of Example \ref{ex:pointed_sets}, generated by one
nullary operation $*$.
Let $\phi$ be the associated regular epi.
Then $\Wkwrt P \phi$ has one nullary object and no objects of any other arity; the only arrow is the identity on the unique nullary object.
In fact, $\Wkwrt P \phi = D_* P$.
So a weak algebra for this theory and this generating set is a category \C\
with a distinguished object $\hat * \in \C$.
A $\phi$-weak $P$-functor from $(\C, (\hatmap))$ to $(\D, (\barmap))$ is a
functor $F : \C \to \D$ and an isomorphism $\bar * \stackrel{\sim}{\to} \hat *$.
\end{example}
\begin{example}
\label{ex:ptd_sets_ABCD}
\index{pointed sets}
Consider again the operad $P$ of Example \ref{ex:pointed_sets}, this time
generated by four nullary operations $A,B,C,D$ (which are all set equal to each
other).
Let $\phi$ be the associated regular epi.
Then $\Wkwrt P \phi_0$ is the indiscrete category on the four objects $A, B, C,
D$, and $\Wkwrt P \phi_i$ is empty for all other $i \in \natural$.
Hence a $\phi$-weak $P$-category is a category \C\ containing four specified
objects $\hat A, \hat B, \hat C$ and $\hat D$.
These four objects are isomorphic via specified isomorphisms $\delta_{AB},
\delta_{AC}, \delta_{AD}$ etc, and all diagrams involving these isomorphisms
commute:
\[
\xymatrixcolsep{6pc}
\xymatrixrowsep{4pc}
\xymatrix{
\hat A \ar[r]^{\delta_{AB}} \ar[dr]_(0.3){\delta_{AD}} \ar[d]_{\delta_{AC}}
& \hat B \ar[d]^{\delta_{BD}} \ar[dl]^(0.3){\delta_{BC}} \\
\hat C \ar[r]_{\delta_{CD}} & \hat D
}
\]
and $\delta_{XY}\delta_{YX} = 1_{\hat X}$ for all $X,Y \in \{A,B,C,D\}$.

Let $(\C, (\hatmap))$ and $(\D, (\barmap))$ be $\phi$-weak $P$-categories.
A $\phi$-weak $P$-functor $(\C, (\hatmap)) \to (\D, (\barmap))$ consists of
\begin{itemize}
\item a functor $F : \C \to \D$,
\item an isomorphism $\phi_{XY} : \bar X \stackrel{\sim}{\to} F
\hat X$ for all $X \in \{A, B, C, D\}$,
\end{itemize}
such that, for all $X, Y \in \{A, B, C, D\}$, there is precisely one
isomorphism $\bar X \to F \hat Y$ formed by compositions of $\delta$s and
$\phi$s.
\end{example}

\section{Symmetric monoidal categories}
\index{commutative monoids}
\label{sec:symmmoncats}
\def\Csmc{\ensuremath{(\C, \otimes, I, \alpha, \lambda, \rho, \tau)}}
\def\Csmcprime{\ensuremath{(\C, \otimes', I', \alpha', \lambda', \rho', \tau')}}
\def\Csmcpprime{\ensuremath{(\C', \otimes', I', \alpha', \lambda', \rho',
	\tau')}}
\def\can{{\mbox{\upshape{can}}}}
\def\Q{{\Wkwrt P \phi}}

Consider the terminal symmetric operad $P$, whose algebras in \Set\ are
commutative monoids, and the following linear presentation $(\Phi, E)$ for $P$:
\begin{itemize}
\item $\Phi_0 = \{e\}, \Phi_2 = \{.\}$, all other $\Phi_i$s are empty;
\item $E$ contains the equations
\begin{enumerate}
\item $x_1.(x_2.x_3) = (x_1.x_2).x_3$
\item $e.x_1 = x_1$
\item $x_1.e = x_1$
\item $x_1.x_2 = x_2.x_1$
\end{enumerate}
\end{itemize}
This linear presentation gives rise to a symmetric-operadic presentation
$(\Phi, E)$, as described in Lemma \ref{lem:syntactic_semantic_equiv}.
Let $\phi : \Fsig \Phi \to P$  be the coequalizer in the diagram
\[
\NoCompileMatrices
\fork {\Fsig E} {} {} {\Fsig \Phi} \phi P
\]

We shall now prove that the algebras for $\Q$ are classical symmetric monoidal
categories.
More precisely, we shall show the following:
\begin{enumerate}
\item for a given category \C, the $\Q$-algebra structures on \C\ are in
one-to-one correspondence with the \smc\ structures on \C;
\item there exists an isomorphism (which we construct) between the category
$\WkPCat$ and the category of symmetric monoidal categories and weak functors;
\item the isomorphism in (2) respects the correspondence in (1).
\end{enumerate}

To fix notation, we recall the classical notions of symmetric monoidal category
and symmetric monoidal functor:

\begin{defn} 
\index{symmetric monoidal category}
A \defterm{symmetric monoidal category} is a 7-tuple $\Csmc$, where
\begin{itemize}
\item $\C$ is a category;
\item $\otimes : \C \times \C \to \C$ is a functor;
\item $I$ is an object of $\C$,
\item $\alpha: A \otimes (B \otimes C) \to (A \otimes B) \otimes C$ is natural
in $A, B, C \in \C$;
\item $\lambda : I \otimes A \to A$ and $\rho : A \otimes I \to A$ are natural in $A \in \C$;
\item $\tau : A \otimes B \to B \otimes A$ is natural in $A, B \in \C$,
\end{itemize}
$\alpha, \lambda, \rho, \tau$ are all invertible, and the following diagrams
commute:
\begin{equation}
\xymatrix{
& (A \otimes B) \otimes (C \otimes D) \ar[ddr]^\alpha \\
A \otimes (B \otimes (C \otimes D)) \ar[dd]_{1\otimes\alpha} \ar[ur]^\alpha \\
& & ((A \otimes B) \otimes C) \otimes D \\
A \otimes ((B \otimes C) \otimes D) \ar[dr]_\alpha \\
& (A \otimes (B \otimes C) ) \otimes D \ar[uur]_{\alpha \otimes 1}
}
\end{equation}
\begin{equation}
\xymatrix{
A \otimes (I \otimes C) \ar[rr]^\alpha \ar[ddr]_{1 \otimes \lambda}
& & (A \otimes I) \otimes C \ar[ddl]^{\rho \otimes 1} \\ \\
& A \otimes C
}
\end{equation}
\begin{equation}
\xymatrix{
A \otimes I \ar[rr]^\tau \ar[ddr]_\rho && I \otimes A\ar[ddl]^\lambda \\ \\
& A
}
\end{equation}
\begin{equation}
\xymatrix{
(A \otimes B) \otimes C \ar[r]^\tau \ar[d]_{\alpha^{-1}}
& C \otimes (A \otimes B) \ar[d]^\alpha \\
A \otimes (B \otimes C) \ar[d]_{1 \otimes \tau}
& (C \otimes A) \otimes B \ar[d]^{\tau \otimes 1} \\
A \otimes (C \otimes B) \ar[r]^\alpha
& (A \otimes C) \otimes B
}
\hskip 2em
\xymatrix{
A \otimes (B \otimes C) \ar[r]^\tau \ar[d]_\alpha
& (B \otimes C) \otimes A \ar[d]^{\alpha^{-1}} \\
(A \otimes B) \otimes C \ar[d]_{\tau \otimes 1}
& B \otimes (C \otimes A) \ar[d]^{1 \otimes \tau} \\
(B \otimes A) \otimes C \ar[r]^{\alpha^{-1}}
& B \otimes (A \otimes C)
}
\end{equation}
\begin{equation}
\xymatrix{
A \otimes B \ar[dr]_1 \ar[r]^{\tau_{A,B}} & B \otimes A \ar[d]^{\tau_{B,A}} \\
& A \otimes B.
}
\end{equation}
\end{defn}

\begin{defn}
Let $M = \Csmc$ and $N = \Csmcpprime$ be symmetric monoidal categories.
\index{symmetric monoidal functor}
A \defterm{lax symmetric monoidal functor} $F: M \to N$ consists of
\begin{itemize}
\item a functor $F : \C \to \C'$,
\item morphisms $F_2 : (F A) \otimes' (F B) \to F(A \otimes B)$, natural in
$A,B \in \C$,
\item a morphism $F_0 : I' \to F I$ in $\C'$,
\end{itemize}
such that the following diagrams commute:
\begin{equation}
\label{eqn:symm_mon_fctr_assoc}
\xymatrix{
FA \otimes' (FB \otimes' FC) \ar[r]^{\alpha'} \ar[d]_{1 \otimes' F_2}
& (FA \otimes' FB) \otimes' FC \ar[d]^{F_2 \otimes' 1} \\
FA \otimes' (F(B \otimes C)) \ar[d]_{F_2}
& F(A \otimes B) \otimes' FC \ar[d]^{F_2} \\
F(A \otimes (B \otimes C)) \ar[r]^{F\alpha}
& F((A \otimes B) \otimes C)
}
\end{equation}
\begin{equation}
\label{eqn:symm_mon_fctr_unit}
\xymatrix{
(FB) \otimes' I' \ar[r]^{\rho'} \ar[d]_{1 \otimes' F_0}FB
& F B \\
(FB) \otimes' (FI) \ar[r]^{F_2}
& F(B \otimes I) \ar[u]_{F\rho}
}
\hskip 2em
\xymatrix{
I' \otimes' (FB) \ar[r]^{\lambda'}
& FB \\
(FI) \otimes' (FB) \ar[r]^{F_2} \ar[u]^{F_0 \otimes' 1}
& F(I \otimes B) \ar[u]_{F\lambda}
}
\end{equation}
\begin{equation}
\label{eqn:symm_mon_fctr_tau}
\xymatrix{
(FA)\otimes'(FB) \ar[r]^{\tau'} \ar[d]_{F_2}
& (FB) \otimes' (FA) \ar[d]^{F_2} \\
F(A \otimes B) \ar[r]^{F\tau}
& F(B \otimes A).
}
\end{equation}
$F$ is said to be \defterm{weak} when $F_0, F_2$ are isomorphisms, and
\defterm{strict} when $F_0, F_2$ are identities.
\index{lax!symmetric monoidal functor}
\index{weak!symmetric monoidal functor}
\index{strict!symmetric monoidal functor}
\end{defn}

Recall also the coherence theorem for classical symmetric monoidal categories.
For any \nary\ permuted $\Phi$-tree $(\sigma \act t)$, let $(\sigma \act t)_M$
be the functor $M^n \to M$ obtained by replacing every $.$ in $t$ by $\otimes$
and every $e$ by $I$, and permuting the arguments according to $\sigma$, so
$(\sigma \act t)_M(A_1, \dots, A_n) = t_M(A_{\sigma 1}, \dots, A_{\sigma n})$
for all $A_1, \dots, A_n \in M$.
In particular, we do not make use of the symmetry maps on $M$ in constructing
these functors.
Then:

\begin{theorem}
\label{thm:smc_coherence}
(Mac Lane) \index{Mac Lane, Saunders}
\index{coherence!for symmetric monoidal categories}
In each weak symmetric monoidal category $M$ there is a function which assigns
to each pair $(\sigma \act t_1, \rho \act t_2)$ of permuted $\Phi$-trees of
the same arity $n$ a unique natural isomorphism
\[
\can_M (\sigma \act t_1, \rho \act t_2) : (\sigma \act t_1)_M
\to (\rho \act t_2)_M : M^n \to M
\]
called the \defterm{canonical map} from $\sigma \act t_1$ to $\rho \act t_2$,
in such a way that the identity of $M$ and all instances of $\alpha, \lambda,
\rho$ and $\tau$ are canonical, and the composite as well as the
$\otimes$-product of two canonical maps is canonical.
\end{theorem}
\begin{proof}
See \cite{catwork} XI.1.
\end{proof}

Finally, recall the coherence theorem for weak monoidal functors:

\begin{lemma}
\label{lem:mon_fctr_coherence}
\index{coherence!for monoidal functors}
Let $M, N$ be monoidal categories, and $F : M \to N$ be a weak monoidal
functor.
For every $n \in \natural$ and every strongly regular $\Phi$-tree $v$ of arity
$n$, there is a unique map $F_v : v_N(F A_1, \dots, F A_n) \to Fv_M(A_1, \dots,
A_n)$ natural in $A_1, \dots, A_n \in M$ and formed by taking composites and
tensors of $F_0$ and $F_2$, such that the diagram
\[
\commsquare{v_N(F A_1, \dots F A_n)}{F_{v}}
{Fv_M(A_1, \dots, A_n)} {\can_N}{F\can_M}
{w_N(F A_1, \dots F A_n)}
{F_{w}}{F(w)_M(A_1, \dots, A_n)}
\]
commutes for all $n \in \natural$, all $v, w \in (\Fp \Phi)_n$, and all $A_1,
\dots, A_n \in \M$.
\end{lemma}
\begin{proof}
See \cite{catwork}, p. 257.
\end{proof}

We may use this result to sketch a proof of a coherence theorem for weak
symmetric monoidal functors:

\begin{theorem}
\label{thm:smcf_coherence}
\index{coherence!for symmetric monoidal functors}
Let $M,N$ be \smcs, and $F:M \to N$ be a weak symmetric monoidal functor.
Let $\sigma\act v$ be an \nary\ permuted $\Phi$-tree.
Then there is a unique natural transformation
\[
\xymatrix{
	M^n \ar[r]^{F^n} \ar[d]_{(\sigma\act v)_M}
		\drtwocell\omit{^*{!(1,-1)\object{F_{\sigma\act v}}}}
	& N^n \ar[d]^{(\sigma\act v)_N} \\
	M \ar[r]_F & N
}
\]
formed by composing tensor products of $F_2$ and $F_0$, possibly with their
arguments permuted.
Furthermore, if $\rho\act w$ is another permuted $\Phi$-tree, then the diagram
\[
\commsquare{(\sigma\act v)_N(F A_1, \dots F A_n)}{F_{\sigma\act v}}
{F(\sigma \act v)_M(A_1, \dots, A_n)} {\can_N}{F\can_M}
{(\rho \act w)_N(F A_1, \dots F A_n)}
{F_{\rho \act w}}{F(\rho\act w)_M(A_1, \dots, A_n)}
\]
commutes.
\end{theorem}
\begin{proof}
\def\perm{\mbox{\upshape{perm}}}
Let $F_{\sigma\act v}(A_1, \dots, A_n) = F_v(A_{\sigma(1)}, \dots,
A_{\sigma(n)})$, and similarly on morphisms.
Then $F_{\sigma\act v}$ has the required type. 
We may decompose $\can_M(\sigma\act v, \rho\act w)$ as
$\perm_M(\sigma,\rho)\ \can_M(v, w)$, where $\perm_M(\sigma,\rho) :
F_{\sigma\act v} \to F_{\rho\act v}$ is a composite of $\tau$s.

Equation \ref{eqn:symm_mon_fctr_tau} and Lemma \ref{lem:mon_fctr_coherence}
together imply that the diagram
\[
\xymatrix{
(\sigma\act v)_N(F A_1, \dots F A_n)\ar[r]^{F_{\sigma\act v}} \ar[d]_{\perm_N}
& {F((\sigma \act v)_M(A_1, \dots, A_n))} \ar[d]^{F\perm_M} \\
(\rho\act v)_N(F A_1, \dots F A_n)\ar[r]^{F_{\rho\act v}} \ar[d]_{\can_N}
& {F((\rho \act v)_M(A_1, \dots, A_n))} \ar[d]^{F\can_M} \\
{(\rho \act w)_N(F A_1, \dots F A_n)} \ar[r]^{F_{\rho \act w}}
& {F((\rho\act w)_M(A_1, \dots, A_n))}
}
\]
commutes.
It remains to show that $F_{\sigma \act v}$ is unique with this property.

Suppose that $F_{\sigma\act v}$ is not unique for some $\sigma \act v$, and
that there exists some natural transformation $G: (\sigma \act v)_N(FA_1,
\dots, FA_n) \to F((\sigma\act v)_M (A_1, \dots, A_n))$, composed of tensor
products of components of $F_0$ and $F_2$, such that $G \neq F_{\sigma\act
v}$.
Suppose further that $\sigma\act v$ and $G$ have been chosen to be a minimal
counterexample, in the sense that of all such counterexamples, $\sigma$ may be
written as a product of the smallest number of transpositions.
If no transpositions are used, then we have a contradiction, because then
$\sigma = 1_n$, and Lemma \ref{lem:mon_fctr_coherence} tells us that $G = F_v$.
But suppose $\sigma = t_1t_2\dots t_m$ where each $t_i$ is a transposition:
then $t_1 \act G$ is a natural transformation $(\sigma\act v)_N (FA_{t_1
1}, \dots, FA_{t_1 n}) \to F((\sigma\act v)_M(A_{t_1 1}, \dots, A_{t_1 n}))$,
and thus a transformation $(t_1\sigma\act v)_N (FA\seq) \to F((t_1\sigma\act v)_M(A\seq))$.
But $t_1\sigma = t_2t_3\dots t_m$, and thus (by minimality of $\sigma$), it
must be the case that $t_1 \act G = F_{t_1\sigma\act v} = t_1\act F_{\sigma\act
v}$.
Hence $G = F_{\sigma \act v}$.
\end{proof}

We now proceed to relate the classical theory of \smcs\ to the more general
notion of categorification we developed in previous sections.

By Lemma \ref{lem:explicit_wkp}, if $\tau_1, \tau_2$ are \nary\ 1-cells in $\Q$
(in other words \nary\ permuted $\Phi$-trees), there is a (unique) 2-cell
$\tau_1 \to \tau_2$ in $\Q$ iff $\tau_1 \sim \tau_2$ under the congruence
generated by $E$.
By standard properties of commutative monoids, this relation holds iff $\tau_1$
and $\tau_2$ take the same number of arguments, so there is exactly one 2-cell
$\tau_1 \to \tau_2$ for every $n \in \natural$ and every pair $(\tau_1,
\tau_2)$ of \nary\ 1-cells in $\Q$.

Let $\SMCCat$ denote the category of \smcs\ and weak maps between them.
We shall define functors $S : \SMCCat \to \WkPCat$ and $R : \WkPCat \to
\SMCCat$, and show that they are inverses of each other.

Let $M = \Csmc$ be a \smc.
Let $SM$ be the weak $P$-category $(\hatmap) :\Q \to \End(\C)$ defined as
follows:
\begin{itemize}
\item On 1-cells of $\Q$, $(\hatmap)$ is determined by $\hat . = \otimes$ and
$\hat e = I$.
\item If $\delta : \tau_1 \to \tau_2$ is an \nary\ 2-cell in $\Q$ (i.e. a
morphism in the category $\Q_n$), let $\hat \delta$ be the canonical map $\hat
\tau_1 \to \hat \tau_2$.
\end{itemize}

\begin{lemma}
$SM$ is a well-defined $\Q$-algebra for all $M \in \SMCCat$.
\end{lemma}
\begin{proof}
The 1-cells of $\Q$ are the same as those of $\Fsig \Phi$; hence, $(\hatmap)$ is
entirely determined on 1-cells by a map of signatures $\Phi \to \Usig \End(\C)$,
which we have given.
On 2-cells, Theorem \ref{thm:smc_coherence} and the uniqueness property of
2-cells in $\Q$ tell us that if $\delta_1, \delta_2$ are 2-cells in $\Q$, then
$\widehat{\delta_1 . \delta_2} = \hat \delta_1 \otimes \hat \delta_2 = \hat
\delta_1 \hat . \hat \delta_2$, and $\widehat{\delta_1 \delta_2} =
\hat{\delta_1} \hat{\delta_2}$ wherever $\delta_1, \delta_2$ are composable.
Hence, $(\C, (\hatmap))$ is a well-defined $\Q$-algebra.
\end{proof}

Given symmetric monoidal categories $M$ and $N$, and a weak symmetric monoidal
functor $F : M \to N$, we would like to define a weak $P$-functor $SF = (F,
\psi): SM \to SN$.
Let $\psi_{\sigma\act v, A\seq} = F_{\sigma\act v}$ for all $n \in \natural$, all $\sigma\act v \in (\Fsig \Phi)_n$, and all $A_1, \dots, A_n \in M$.
By Theorem \ref{thm:smcf_coherence}, this is natural in $\sigma\act v$ and in
$A_1, \dots, A_n$.
The other axioms for a weak $P$-functor are all implied by the coherence
theorem (Theorem \ref{thm:smcf_coherence}).
This can be generalized: a lax symmetric monoidal functor $F$ determines a
lax $P$-functor $SF$, and a strict symmetric monoidal functor $F$
determines a strict $P$-functor $SF$.

Now, let \C\ be a $\Q$-algebra, with map $(\hatmap): \Q \to \End(\C)$.
We shall construct a \smc\ $R(\C, (\hatmap)) = \Csmc$.
Take
\begin{itemize}
\item $\otimes = \hat .$
\item $ I = \hat e$
\item $ \alpha = \hat \delta_1$, where $\delta_1 : -.(-.-) \to (-.-).-$ in
$\Q_3$,
\item $ \lambda = \hat \delta_2$, where $\delta_2 : e.- \to -$ in $\Q_1$,
\item $ \rho = \hat \delta_3$, where $\delta_3 : -.e \to e$ in $\Q_1$,
\item $ \tau = \hat \delta_4 $, where $\delta_4 : (-.-) \to (12)\cdot(-.-)$ in $\Q_2$.
\end{itemize}

\begin{lemma}
$R(\C, (\hatmap))$ is a symmetric monoidal category.
\end{lemma}
\begin{proof}
Because there is at most one 2-cell $\tau_1 \to \tau_2$ for any pair of 1-cells
$\tau_1, \tau_2$ in $\Q$, all diagrams involving these commute.
In particular, the axioms for a \smc\ are satisfied.
The 2-cells in $\End(\C)$ are natural transformations, so $\alpha, \lambda,
\rho$ and $\tau$ (as images of 2-cells in $\Q$ under the map $(\hatmap) : \Q \to
\End(\C)$) are natural transformations.
All 2-cells in $\Q$ are invertible, so $\alpha, \lambda, \rho$ and $\tau$ are
all natural isomorphisms.
Hence $(\C, \otimes, I, \alpha, \lambda, \rho, \tau)$ is a \smc.
\end{proof}

Let $(F, \psi) : (\C, (\hatmap)) \to (\C', (\checkmap))$ be a weak morphism of
$\Q$-algebras.
Then let $R(F,\psi) : R(\C, (\hatmap)) \to R(\C', (\checkmap))$ be the
following symmetric monoidal functor:
\begin{itemize}
\item the underlying functor is $F$,
\item $F_0$ is $\psi_{e \ocomp 1}: \check e \to F \hat e$,
\item $F_2$ is $\psi_{. \ocomp 1}: (\check . )F^2 \to F (\hat .)$.
\end{itemize}
The coherence diagrams (\ref{eqn:symm_mon_fctr_assoc}),
(\ref{eqn:symm_mon_fctr_unit}) and (\ref{eqn:symm_mon_fctr_tau}) all commute
by virtue of the coherence axioms for a weak morphism of $\Q$-algebras and the
naturality of $\psi$.
Hence $(F, F_0, F_2)$ is a symmetric monoidal functor.

\begin{lemma}
\label{lem:RS = 1}
Let \Csmc\ be a \smc. Then \[ RS\Csmc = \Csmc. \]
\end{lemma}
\begin{proof}
Let $RS\Csmc = \Csmcprime$.
Their underlying categories are equal, both being $\C$.
\[
\begin{array}{rcccl}
\otimes' &= &\hat . &= &\otimes \\
I' &= & \hat e &= & I \\
\alpha' & = &\hat \delta_1 & = & \alpha \mbox{, the unique canonical map of
the correct type}\\
\lambda' & = &\hat \delta_2 & = & \lambda \\
\rho' & = &\hat \delta_3 & = & \rho \\
\tau' & = &\hat \delta_4 & = & \tau \\
\end{array}
\]
\end{proof}

\begin{lemma}
\label{lem:SR = 1}
Let $(\C, (\hatmap))$ be a $\Q$-algebra, and let
$(\C', (\checkmap)) = SR(\C, (\hatmap))$.
Then $(\C, (\hatmap)) = (\C', (\checkmap))$.
\end{lemma}
\begin{proof}
Their underlying categories are the same.
As above, $(\checkmap)$ is determined on objects by the values it takes on $.$
and $e$: these are $\otimes = \hat .$ and $I = \hat e$ respectively.
So $(\checkmap) = (\hatmap)$ on objects.
If $\delta : \tau_1 \to \tau_2$, then $\check \delta$ is the unique canonical
map from $\check \tau_1 \to \check \tau_2$, which, by an easy induction, must be
$\hat \delta$.
So $(\checkmap) = (\hatmap)$, and hence $(\C, (\hatmap)) = SR(\C, (\hatmap))$.
\end{proof}

\begin{lemma}
\label{lem:RS = 1 on morphisms}
Let $M = \Csmc$ and $N = \Csmcpprime$ be \smcs, and let $(F, F_0, F_2)$ be a
weak symmetric monoidal functor $M \to N$.
Then $RS(F, F_0, F_2) = (F, F_0, F_2)$.
\end{lemma}
\begin{proof}
Let $(G, G_0, G_2) = RS(F, F_0, F_2)$. 
Then $G$ is the underlying functor of $S(F, F_0, F_2)$ which is $F$, and 
$G_0, G_2$ are both the canonical maps with the correct types given by
Theorem \ref{thm:smcf_coherence}: that is to say, they are $F_0$ and $F_2$
respectively.
\end{proof}

\begin{lemma}
\label{lem:SR = 1 on morphisms}
Let $(\C, (\hatmap))$ and $(\C', (\checkmap))$ be $\Q$-algebras, and let $(F,
\phi) : (\C, (\hatmap)) \to (\C', (\checkmap))$ be a weak morphism of
$\Q$-algebras.
Then $SR(F, \phi) = (F, \phi)$.
\end{lemma}
\begin{proof}
Let  $(G, \gamma) = SR(F, \phi)$.
Then $G$ is the underlying functor of $R(F, \phi)$, which is $F$.
Let $(F, F_0, F_2) = R(F, \phi)$.
Each component of $\gamma$ is then by definition the correct component of the
canonical map arising from $F_0, F_2$ in the process described in Theorem
\ref{thm:smcf_coherence}.
By the ``uniqueness'' part of the Theorem, this must be the corresponding
component of $\phi$.
Hence $\gamma = \phi$.
\end{proof}

\begin{theorem}
$S$ and $R$ form an isomorphism of categories $\SMCCat \cong \WkPCat$.
\end{theorem}
\begin{proof}
Lemmas \ref{lem:RS = 1} and \ref{lem:SR = 1} show that $R$ and $S$ are
bijective on objects; Lemmas \ref{lem:RS = 1 on morphisms} and \ref{lem:SR = 1
on morphisms} show that $R$ and $S$ are locally bijective on morphisms.
Hence, $R$ and $S$ are a pair of mutually inverse isomorphisms of categories.
\end{proof}

\section{Multicategories}

We can tell this whole story for (symmetric) multicategories as well as just
operads.
We sketch this development briefly here, although the remainder of the thesis
will continue to focus on the special case of operads.

\begin{defn} A (directed) \defterm{multigraph} consists of
\index{multigraph}
\begin{enumerate}
\item a set of \defterm{vertices} $V$,
\item for each $n \in \natural$ and each sequence $v_1, v_2, \dots,
v_n, w$ of vertices, a set $E(v_1, \dots, v_n; w)$ of \defterm{funnels from
$v_1, \dots, v_n$ to $w$}. 
\index{funnel}
\end{enumerate}
\end{defn}

\begin{defn} Let $M_1 = (V_1, E_1)$ and $M_2 = (V_2, E_2)$ be multigraphs.
\index{morphism!of multigraphs}
A \defterm{morphism of multigraphs} $f: M_1 \to M_2$ is
\begin{enumerate}
\item a function $f_V : V_1 \to V_2$,
\item for each finite sequence $v_1, v_2 \dots v_n, w$ of vertices in $M_1$, a
function
\[
f^{v_1,\dots,v_n}_w : E_1(v_1, \dots, v_n; w)
	\to E_2(f_V(v_1), \dots, f_V(v_n); f_V(w)).
\]
\end{enumerate}
\end{defn}

\begin{figure}[h]
\centerline{
	\epsfxsize=290pt
	\epsfbox{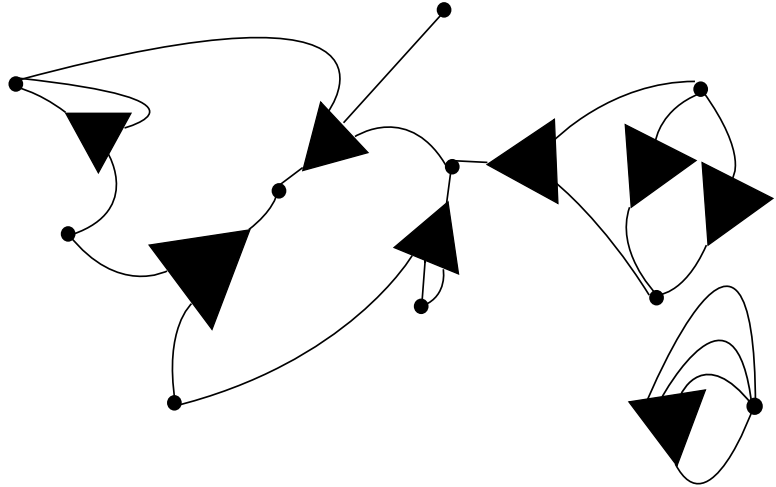}
}
\caption{A multigraph}
\label{fig:mgraph}
\end{figure}

We say that a funnel $f \in E(v_1, \dots, v_n; w)$ has \defterm{source} $v_1,
\dots, v_n$ and \defterm{target} $w$; we say that two funnels are
\defterm{parallel} if they have the same source and target.
The reason for the ``funnel'' terminology should be clear from Figure
\ref{fig:mgraph}.
We shall say that a multigraph has some property $P$ \defterm{locally} if every
$E(v_1, \dots, v_n; w)$ is $P$, and similarly a morphism $f$ of multicategories
is locally $P$ if every $f^{v_1,\dots,v_n}_w$ is $P$.

Multigraphs and their morphisms form a category which we shall call \Multigraph.
\index{\Multigraph}

In order to proceed with the rest of the construction, we will need to consider
subcategories of \Multicat, \Multigraph\ etc.

\begin{defn}
\index{$\Multigraph_X$}
Let $X$ be a set.
Then $\Multigraph_X$ is the subcategory of $\Multigraph$ whose objects are
multigraphs with vertex set $X$, and whose morphisms are identity-on-vertices
maps of multigraphs.
We define $\Multicat_X$ and $\SymmMulticat_X$ similarly.
\end{defn}

For each $X \in \Set$, there is a chain of adjunctions similar to that given in
Section \ref{sec:adj}:
\index{adjunctions!for multicategories}
\[
\xymatrix{
  \DMCat_X \ar@<1.2ex>[d]^\Ufpsig
	 \ar@/r1.4cm/[dd]^\Ufpp
	 \ar@/r3cm/[ddd]^\Ufp \\
  \SymmMulticat_X \ar@<1.2ex>[d]^\Usigp \ar@<1.2ex>[u]^\Ffpsig_{\dashv}
	  \ar@/r1.4cm/[dd]^\Usig \\
  \Multicat_X \ar@<1.2ex>[u]^\Fsigp_{\dashv} \ar@<1.2ex>[d]^\Up
	 \ar@/l1.4cm/[uu]^\Ffpp \\
  \Multigraph_X \ar@<1.2ex>[u]^\Fp_{\dashv} \ar@/l1.4cm/[uu]^\Fsig 
  	 \ar@/l3cm/[uuu]^\Ffp
}
\]

These adjunctions are monadic, by Lemma \ref{lem:frees_exist} and Lemma
\ref{lem:monadj}.
Note that $\Set^\natural$ can be regarded as $\Multigraph_1$: thus, the
adjunctions of Section \ref{sec:adj} are just the restrictions of the
adjunctions above to the one-vertex case.

We can consider multigraphs enriched in some category \V:

\begin{defn}
\index{multigraph!enriched}
Let \V\ be a category.
A \defterm{\V-multigraph} $M = (V,E)$ consists of
\begin{enumerate}
\item a set $V$ of \defterm{vertices},
\item for each $n \in \natural$ and each finite sequence $v_1, v_2 \dots v_n, w$
of vertices, an object of \V\ called $E(v_1, \dots, v_n; w)$ of
\defterm{funnels from $v_1, \dots, v_n$ to $w$}. 
\end{enumerate}
\end{defn}

\begin{defn} Let $M_1 = (V_1, E_1)$ and $M_2 = (V_2, E_2)$ be \V-multigraphs.
A \defterm{morphism of \V-multigraphs} $f: M_1 \to M_2$ is
\index{morphism!of enriched multigraphs}
\begin{enumerate}
\item a function $f_V : V_1 \to V_2$,
\item for each $n \in \natural$ and each finite sequence $v_1, v_2 \dots v_n,
w$ of vertices in $M_1$, an arrow $E_1(v_1, \dots, v_n; w)
	\to E_2(f_V(v_1), \dots, f_V(v_n); f_V(w))$ in \V.
\end{enumerate}
\end{defn}

The category of \V-multigraphs and their morphisms is called \VMultigraph.
The category whose objects are \V-multigraphs with vertex-set $X$ and whose
morphisms are identity-on-vertices maps is called $\VMultigraph_X$.
In particular, we shall consider multigraphs enriched in the category \Digraph\
of directed graphs.
An object of the category \DigraphMultigraph\ consists of
\begin{enumerate}
\item \defterm{vertices} (or \defterm{objects});
\item \defterm{funnels}, each of which has one object as its target, and a
sequence of objects as its source;
\item \defterm{edges}, which each have one funnel as a source and one as a
target: the source and target of a given edge must be parallel.
\end{enumerate}

\begin{figure}[h]
\centerline{
	\epsfxsize=290pt
	\epsfbox{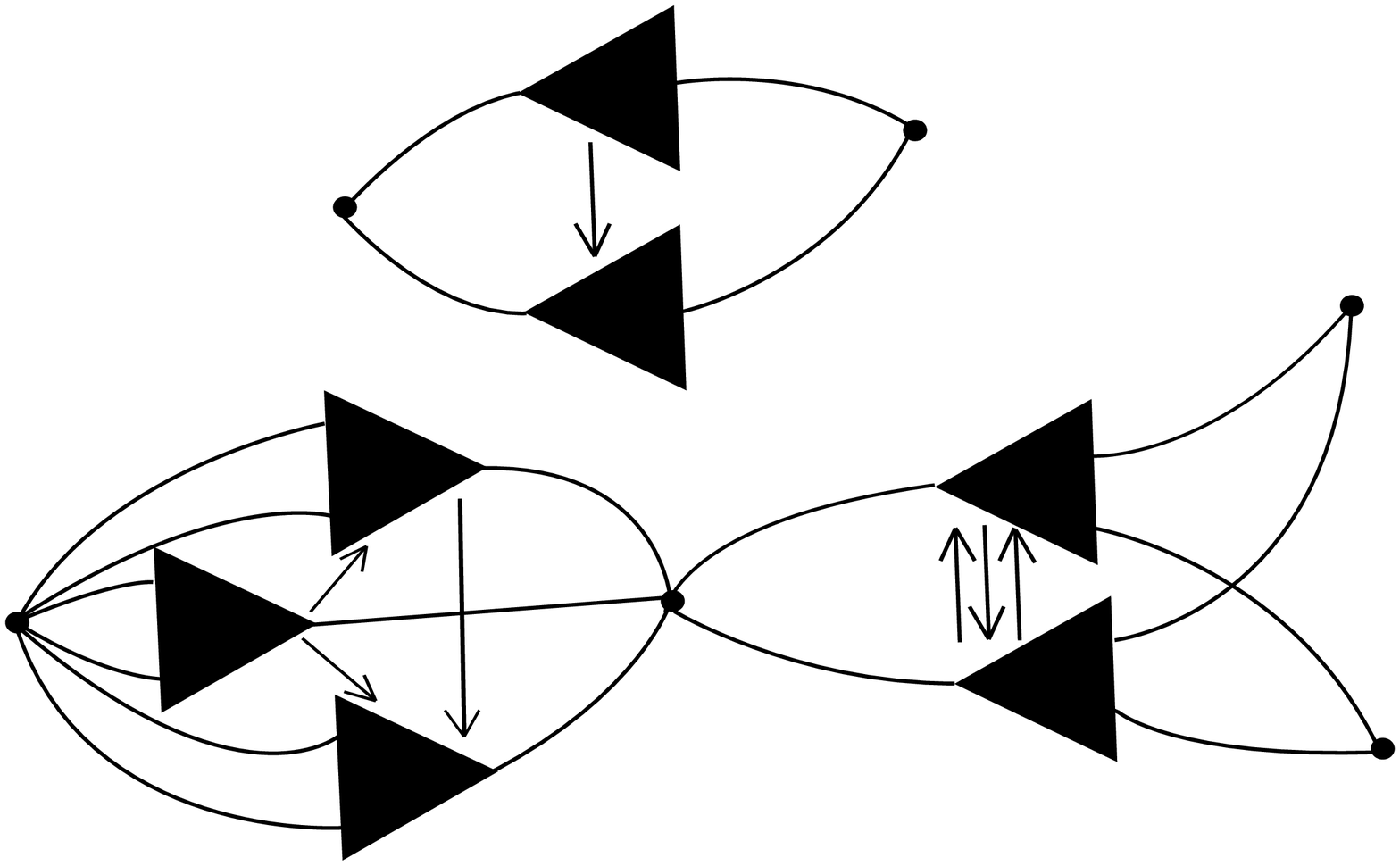}
}
\caption{A multigraph enriched in directed graphs}
\label{fig:gphmgph}
\end{figure}

The factorization system construction of Example \ref{ex:catopdFS} works in
this broader setting too.
Let $X$ be a set.
\index{factorization system!bijective on objects/full and faithful!on
$\DigraphMultigraph_X$}
The factorization system $(\E, \M)$ on \Digraph\ of Example \ref{ex:digraphFS}
gives rise to a factorization system $(\E', \M')$ on $\DigraphMultigraph_X$,
where \E\ consists of maps which are bijective on objects and funnels, and \M\
consists of maps which are locally full-and-faithful.
This lifts to a factorization system $(\E'', \M'')$ on $\CatMultigraph_X$ via
Lemma \ref{lem:monad}.
By the usual argument, there is a chain of monadic adjunctions:
\index{adjunctions!for \Cat-multicategories}
\[
\xymatrixrowsep{3pc}
\xymatrix{
  \DCMCat_X \ar@<1.2ex>[d]^\Ufpsig
	 \ar@/r1.7cm/[dd]^\Ufpp
	 \ar@/r3cm/[ddd]^\Ufp \\
  \CatSymmMulticat_X \ar@<1.2ex>[d]^\Usigp \ar@<1.2ex>[u]^\Ffpsig_{\dashv}
	  \ar@/r1.7cm/[dd]^\Usig \\
  \CatMulticat_X \ar@<1.2ex>[u]^\Fsigp_{\dashv} \ar@<1.2ex>[d]^\Up
	 \ar@/l1.7cm/[uu]^\Ffpp \\
  \CatMultigraph_X \ar@<1.2ex>[u]^\Fp_{\dashv} \ar@/l1.7cm/[uu]^\Fsig 
  	 \ar@/l3cm/[uuu]^\Ffp
}
\]
Since $\CatMulticat_X$ is monadic over $\CatMultigraph_X$, this in turn
lifts to a factorization system on $\CatMulticat_X$.
Similarly, we obtain a factorization system on $\CatSymmMulticat_X$.

A \defterm{generator} for a plain multicategory $M$ with object-set $X$ is a
multigraph $\Phi = (X,E)$ together with a regular epi $\Fp \Phi \to M$ in
$\Multicat_X$.
Similarly, a \defterm{generator} for a symmetric multicategory $M$ with
object-set $X$ is a multigraph $\Phi = (X,E)$ together with a regular epi
$\Fsig \Phi \to M$ in $\SymmMulticat_X$.

We can therefore extend Definition \ref{def:phi-wkpcat_plain} above, in the
obvious way. 
Let $D_*$ be the embedding of $\Multicat_X$\ into $\CatMulticat_X$ via the
(full and faithful) discrete category functor applied locally.

\begin{defn}
\index{weakening!of a multicategory}
Let $M$ be a plain multicategory with object-set $X$, and let $\phi: \Fp\Phi
\to M$ be a regular epi in $\Multicat_X$.
Then the \defterm{weakening of $M$ with respect to $\phi$} is the
unique-up-to-isomorphism \Cat-multicategory $\Wkwrt M \phi$ such that the
following diagram commutes:
\[
\xymatrix{
	D_* F\Phi \ar[rr]^{D_*\phi} \booar[dr]_b
	&& D_* M \\
	& \Wkwrt M \phi \lffar[ur]_f
}
\]
where $f$ is locally full and faithful, and $b$ is locally bijective on objects
(i.e., each map of sets of funnels in $b$ is a bijection).
The uniqueness of $\Wkwrt M \phi$ follows from properties of the factorization
system on $\CatMulticat_X$\ given above.
\end{defn}

\begin{defn}
\index{unbiased!weakening of a multicategory}
Let $M$ be a (symmetric) multicategory.
The \defterm{unbiased weakening} of $M$ is the weakening of $M$ with respect to
the counit $\epsilon$ of the adjunction $\Fp \dashv \Up$ (respectively, the
adjunction $\Fsig \dashv \Usig$).
\end{defn}

\begin{defn}
Let $M$ be a multicategory, and let $\phi: \Fp\Phi \to M$ be a regular epi in
$\Multicat_X$.
A \defterm{$\phi$-weak $M$-algebra} is an algebra for $\Wkwrt M \phi$.
An \defterm{unbiased weak $M$-algebra} is an algebra for $\Wk M$.
\end{defn}

We define weakenings of symmetric multicategories analogously.

\section{Examples}

\begin{example}
Let $M$ be the multicategory generated by
\[
\xymatrix{0 \ar[r]^f & 1 \ar[r]^g & 2}
\]
Then a weak algebra for $M$ in \Cat\ with respect to this generating set
consists of a diagram
\[
\xymatrix{\hat 0 \ar[r]^{\hat f} & \hat 1 \ar[r]^{\hat g} & \hat 2}
\]
in \Cat, whereas an unbiased weak $M$-algebra is a diagram
\[
\xymatrix{
\hat 0 \ar[dr]_{\hat f} \ar[rr]^{\widehat{gf}} \rrlowertwocell\omit{<3>\sim}
\ar@(ul,dl)[]_{\widehat{I_0}}
& & \hat 2 \ar@(ur,dr)[]^{\widehat{I_2}} \\
& \hat 1 \ar[ur]_{\hat g} \ar@(dl,dr)[]_{\widehat{I_1}}
}
\]
where $\widehat{I_0}, \widehat{I_1}$ and $\widehat{I_2}$ are trivial monads.
\end{example}

\begin{example}
\index{$M$-sets}
Consider the theory $T$ whose algebras are a monoid $M$ together with an
$M$-set $A$.
Then a weak $T$-algebra with respect to the standard presentation (a binary and
nullary operation on $M$, and an operation $M \times A \to A$) is a classical
monoidal category $\hat M$, a category $\hat A$, and a weak monoidal
functor $\hat M \to \End(\hat A)$.
An unbiased weak $T$-algebra is an unbiased monoidal category $\hat M$, a
category $\hat A$ equipped with a trivial monad $\hat I_A$, and an unbiased
monoidal functor $\hat M \to \End(\hat A)$ which commutes up to coherent
isomorphism with $\hat I_A$.
\end{example}

\begin{example}
\index{weak!$P$-functor}
Let $P$ be an operad, and let $\bar P$ be the multicategory from Section
\ref{sec:Pbar} whose algebras are pairs of $P$-algebras with a morphism between
them.
It seems clear that an unbiased weak $\bar P$-category is a pair of unbiased
weak $P$-categories and an unbiased weak $P$-functor between them; a rigorous
proof would first require a coherence theorem to be proven for weak maps of
$\Cat$-operad algebras, and currently no such theorem is known.
\end{example}

\section{Evaluation}
\label{sec:evaluation}

At the beginning of this chapter, we proposed three criteria that a successful
definition of categorification should satisfy: namely, it should be broad,
consistent with earlier work, and canonical.
The examples considered throughout the chapter show that our theory agrees with
the standard categorifications that are within its scope.
It is determined by the universal property given by the factorization system on
$\CatSymmOperad$: the only tunable parameter is the choice of generator of a
given theory, and in Chapter \ref{ch:coherence} we shall see that the weakening
of a given theory is independent (up to equivalence) of the generator used.
The main problem is the breadth of our theory: as presented, it is restricted
to linear theories, preventing us from categorifying the theories of groups,
rings, Lie algebras, and many other interesting nonlinear theories.
We shall now show what happens when we try to extend our theory to general
algebraic theories.

\begin{lemma}
There is a factorization system $(\E, \M)$ on \DCOCat\, where \E\ is the
collection of maps which are bijective on objects, and \M\ is the collection of
maps which are levelwise full and faithful.
\end{lemma}
\begin{proof}
The proof is exactly as for the proof of existence of the factorization systems
on \CatOperad\ and \CatSymmOperad\ given in Example \ref{ex:catopdFS}.
\end{proof}

\begin{theorem}
Let $P$ be the \dop\ whose algebras are commutative monoids, and $D_* : \DOpCat
\to \DCOCat$ be the levelwise ``discrete category'' functor.
Let $Q$ be the \dco\ given by the factorization
\[
\xymatrix{
D_* \Ffp\Ufp P \ar[rr]^{D_* \epsilon_P} \booar[dr] & & D_* P \\
& Q \lffar[ur]
}
\]
Then an algebra for $Q$ is an unbiased symmetric monoidal category $\C$ such
that, for all $A \in \C$, the component $\tau_{AA}$ of the symmetry map $\tau$
is the identity on $A \otimes A$.
\end{theorem}
\begin{proof}
\def\oneone{\left[1 \atop 1\right]}
We adopt the notation for elements of $P$ introduced in Example
\ref{ex:comm_monoid_fp}.
Let $f$ be the unique function $\fs 2 \to \fs 1$, and let $t : \fs 2 \to \fs 2$
be the permutation transposing 1 and 2.
Then $\epsilon(f \cdot \oneone) = [2] \in P_1$.
Let $(\C, (\hatmap))$ be a $Q$-algebra.
We shall write $\hat \oneone(A, B)$ as $A \otimes B$.
We may impose a symmetric monoidal category structure on $\C$, where the
symmetry map is the image under $(\hatmap)$ of the unique map $\delta : \oneone
\to t \cdot \oneone$ in $Q_1$.
All diagrams in $Q_1$ commute, so in particular, the following diagram
commutes:
\[
\xymatrix{
	[2] \ar[r] & f\cdot\oneone \\
	f\cdot\oneone \ar[ur]_{f\cdot\delta} \ar[u]
}
\]
The two unlabelled arrows are mutually inverse.
Applying $(\hatmap)$ to the entire diagram, and evaluating the resulting
functors at $A \in \C$, we see that the following diagram commutes:
\[
\xymatrix{
	\hat{[2]}(A) \ar[r] & A \otimes A \\
	A \otimes A \ar[ur]_{\tau_{A A}} \ar[u]
}
\]
and hence $\tau_{A A} = 1_{A \otimes A}$.
\end{proof}
This is not the case for most interesting symmetric monoidal categories.
Hence this definition of categorification would fail to be consistent with
earlier work.

\chapter{Coherence}
\label{ch:coherence}

\index{Kelly, Max}
\index{coherence!general case}
There are many ``coherence theorems'' in category theory, but in practice they usually fall into one of two classes:
\begin{enumerate}
\item
\label{all_commute}
``All diagrams commute'', or more precisely, that diagrams in a given class
commute if and only if some quantity is invariant.
\item
\label{weak=strict}
Every ``weak'' object is equivalent to an appropriate ``strict'' object.
\end{enumerate}
Since the diagrams of interest in theorems of type \ref{all_commute} will
usually commute trivially in a strict object, a coherence theorem of type
\ref{weak=strict} usually implies one of type \ref{all_commute}.
However, establishing the converse is usually harder.
In the previous chapter, our ``weak $P$-categories'' were defined explicitly in
terms of an infinite class of commuting diagrams (namely, those diagrams which
become identities under the application of the counit of the adjunctions
$\Fsig \dashv \Usig$ or $\Fp \dashv \Up$): it is therefore interesting to see
if we can prove a theorem of type \ref{weak=strict} about them.
We do this in Section \ref{sec:strictification}, and investigate an interesting property of the strictification functor in Section \ref{stsect}.

In Section \ref{sec:presindep}, we investigate how the operad defining weak
$P$-categories is affected by our choice of presentation for $P$.
While the independence result obtained is not a coherence theorem of the usual
form, it can be seen as a coherence theorem in a higher-dimensional sense: that
the process of categorification is itself coherent.

For other related work, see Power's paper \cite{power}.

\section{Strictification}
\label{sec:strictification}

Let $P$ be a plain operad, and $Q = \mbox{Wk}(P)$, with $\pi:Q \to D_* P$ the
levelwise full-and-faithful map in Theorem \ref{thm:WkP_characterization}.
We again adopt the ${}\seq$ notation from chain complexes and write, for
instance, $p\seq$ for a finite sequence of objects in $P$, and $p\seq\udot$ for
a double sequence.
Let $Q \kel A \toletter{h} A$ be a weak $P$-category.
We shall construct a strict $P$-category \st A and a weak $P$-functor
$(F,\phi):\st A \to A$, and show that it is an equivalence of weak
$P$-categories.
This ``strictification'' construction is closely related to that given in
\cite{j+s} for monoidal categories; however, it is more general, and since we
work for the moment with unbiased weak $P$-categories, our construction has
some additional pleasant properties.

In fact, \stfunc is functorial, and is left adjoint to the forgetful
functor $\StrPCat\linebreak[0] \to \WkPCat$ (see Section \ref{stsect}).
The theorem then says that the
unit of this adjunction is pseudo-invertible, and that the strict
$P$-categories and strict $P$-functors form a weakly coreflective
sub-2-category of \WkPCat.

If $P$ is a plain operad, let $\iota$ be the embedding
\begin{eqnarray*}
\iota : \Up D_* P & \to & \Up \WkP \\
\iota(p) & = & p \ocomp (|, \dots, |)
\end{eqnarray*}
Note that this is a morphism in $\Cat^\natural$, not in $\CatOperad$.

\begin{defn}
\index{strictification}
Let $P$, $Q$, $h$, $A$, $\iota$ be as above. The \defterm{strictification of
$A$}, written \st A, is given by the bijective-on-objects/full and
faithful factorization of $h(\iota \kel 1)$ in $\Cat$:
\[
\xymatrix{
	P\kel A \ar[r]^{\iota \kel 1} \booar[dr] & Q \kel A \ar[r]^h & A. \\
	& \st A \lffar[ur]
}
\]
\end{defn}

We shall show that \st A is a strict $P$-category.
We may describe \st A explicitly as follows:
\begin{itemize}
\item An object of \st A is an object of $P \kel A$.
\item If $p \in P_n$ and $a_1, \dots, a_n \in A$, an arrow $(p, \adot) \to (p',
\adot')$ in \st A is an arrow $h(p, \adot) \to h(p', \adot')$ in $A$.
We say that such an arrow is a \defterm{lifting} of $h(p, \adot)$.
Composition and identities are as in $A$.
\end{itemize}

We define an action $h'$ of $P$ on \st A as follows:
\begin{itemize}
\item On objects, $h'$ acts by $h'(q, (p, \adot)\udot)
= (\pi(p\cmp {p\udot}), \adot\udot)$ where $p \in P_n$ and $(p^i, \adot^i) \in
\st A$ for $n \in \natural$ and $i = 1, \dots n$.
\item Let $f_i :(p_i, a_i) \to (p_i', a_i')$ for $i = 1, \dots,
  n$. Then $h'(p, f\seq)$ is the composite
\[
\begin{array}{rcl}
h(p\ocomp(p\seq), \adot)& \toletter{\delta_{p\cmp {p\seq}}^{-1}}&
        h(p\cmp {p\seq}, \adot) = h(p, h(p_1, a_1), \dots , h(p_n, a_n))\\
    &\toletter{h(p, f\seq)}
        &h(p, h(p_1', a_1'), \dots, h(p_n', a_n'))
        = h(p\cmp {p\seq'}, \adot')\\
    &\toletter{\delta_{p\cmp {p\seq'}}}
        &h(p\ocomp(p\seq'), \adot').
\end{array}
\]
\end{itemize}

\begin{lemma}\st A is a strict $P$-category.\end{lemma}

\begin{proof}
It is clear that the action we have defined is strict and associative on
objects and that $1_P$ acts as a unit: we must show that the action on arrows
is associative.
Let $f_i^j:(p_i^j, a\ijseq) \to (q_i^j, b\ijseq)$, $\sigma \in Q_n$, and
$\tau_i \in Q_{k_i}$ for $j = 1, \dots, k_i$ and $i = 1,\dots,n$.
We wish to show that $h'(\sigma\ocomp(\tau\seq), f\dseq) = h'(\sigma,
h'(\tau_1, f_1\udot), \dots, h'(\tau_n, f_n\udot))$.

The LHS is
\[
\begin{array}{rcl}
h(\sigma\ocomp(\tau\seq)\ocomp(p\dseq), \adot\udot)&
	\toletter{\delta_{\sigma\ocomp(\tau_i)\cmp{p\dseq}}^{-1}}&
h(\sigma\ocomp(\tau\seq), h(p_1^1, a_{1\bullet}^1), \dots, h(p_n^{k_n},
a_{n\bullet}^{k_n}))\\
	&\toletter{h(\sigma\ocomp(\tau\seq), f\dseq)}
&h(\sigma\ocomp(\tau\seq), h(q_1^1, b_{1\bullet}^1), \dots, h(q_n^{k'_n},
b_{n\bullet}^{k'_n}))\\
	&\toletter{\delta_{\sigma\ocomp(\tau\seq)\cmp{q\dseq}}}
&h(\sigma\ocomp(\tau\seq)\ocomp(q\dseq), b\dseq).
\end{array}
\]

The RHS is 
\[
\begin{array}{rcl}
h(\sigma\ocomp(\tau\seq)\ocomp(p\dseq), \adot\udot)&
	\toletter{\delta_{\sigma\cmp{\tau_i\ocomp(p\dseq)}}^{-1}}&
h(\sigma, h(\tau_1\ocomp(p_1\udot), a_{1\bullet}\udot), \dots,
		h(\tau_n\ocomp(p_n\udot), a_{n\bullet}\udot)) \\
	&\toletter{h(\sigma, h'(\tau\seq, f\dseq))}
&h(\sigma, h(\tau_1\ocomp(q_1\udot), b_{1\bullet}\udot), \dots,
		h(\tau_n\ocomp(q_n\udot), b_{n\bullet}\udot))\\
	&\toletter{\delta_{\sigma\cmp{\tau_i\ocomp(p\dseq)}}}
&h(\sigma\ocomp(\tau\seq)\ocomp(q\dseq), b\dseq),
\end{array}
\]
where each $h'(\tau_i, f_i\udot)$ is
\[
\begin{array}{rcl}
h(\tau_i\ocomp(p_i\udot),a_i\udot))
	&\toletter{\delta_{\tau_i\cmp{p_i\udot}}^{-1}}
& h(\tau_i, h(p_i^1, a_{i\bullet}^1), \dots, h(p_i^{k_i}, a_{i\bullet}^{k_i}))\\
	&\toletter{h(\tau_i, f_i\udot)}
& h(\tau_i, h(q_i^1, b_{i\bullet}^1), \dots, h(q_i^{k_i}, b_{i\bullet}^{k_i}))\\
	&\toletter{\delta_{\tau_i\cmp{p_i\udot}}}
&h(\tau_i\ocomp(q_i\udot),b_i\udot).
\end{array}
\]

So the equation holds if the following diagram commutes:

\[
\xymatrix @+10pt {
& h(\sigma\ocomp(\tau\seq)\ocomp(p\dseq), \adot\udot)
	\ar[dl]_{\delta_{\sigma\ocomp(\tau_i)\cmp{p\dseq}}^{-1}}
	\ar[dr]^{\delta_{\sigma\cmp{\tau_i\ocomp(p\dseq)}}^{-1}}
	\ar[d]|{\delta_{\sigma\cmp{\tau\seq}\cmp{p\dseq}}^{-1}} \\
h(\sigma\ocomp(\tau\seq), h(p\dseq, a\dseq))
	\ar[d]|{h(\sigma\ocomp(\tau\seq), f\dseq)}
& h(\sigma, h(\tau\seq\ocomp(p\dseq),a\dseq))
	\ar[l]^{\delta_{\sigma\cmp{\tau\seq}}}
	\ar[d]|{h(\sigma, h(\tau_i, f_i\udot))}
	\ar @{} [dl]|*+[o][F-]{1}
	\ar @{} [dr]|*+[o][F-]{2}
& h(\sigma\ocomp(\tau\seq), h(p\dseq, a\dseq))
	\ar[l]^{h(\sigma, \delta_{\tau\seq\cmp{p\dseq}}^{-1})}
	\ar[d]|{h(\sigma, h'(\tau\seq, f\dseq))}\\
h(\sigma\ocomp(\tau\seq), h(q\dseq, b\dseq))
	\ar[dr]_{\delta_{\sigma\ocomp(\tau\seq)\cmp{q\dseq}}}
& h(\sigma, h(\tau\seq\ocomp(q\dseq), b\dseq))
	\ar[d]|{\delta_{\sigma\cmp{\tau\seq\ocomp(q\dseq)}}}
	\ar[l]^{\delta_{\sigma\cmp{\tau\seq}}}
	\ar[r]_{h(\sigma, \delta_{\tau_i\cmp{p_i\udot}})}
& h(\sigma, h(\tau\seq\ocomp(q\dseq), b\dseq))
	\ar[dl]^{\delta_{\sigma\cmp{\tau\seq\ocomp(q\dseq)}}}\\
&h(\sigma\ocomp(\tau\seq)\ocomp(q\dseq), b\dseq)\\
}
\]
The triangles all commute because all $\delta$s are images of arrows in $Q$,
and there is at most one 2-cell between any two 1-cells in $Q$. \xylabel 2
commutes by the definition of $h'(\tau_i, f_i\udot)$, and \xylabel 1
commutes by naturality of $\delta$.
\end{proof}

\begin{lemma}
\label{lem:equiv<=>Q-equiv}
Let $Q\kel A \toletter{h} A$ and $Q \kel B \toletter{h'} B$
be weak $P$-categories, $(F, \pi):A \to B$ be a weak $P$-functor, and $(F, G,
\eta, \epsilon)$ be an adjoint equivalence in \Cat.
Then $G$ naturally carries the structure of a weak $P$-functor, and $(F, G,
\eta, \epsilon)$ is an adjoint equivalence in \WkPCat.
\end{lemma}

\begin{proof}
We want a sequence $(\psi\seq)$ of natural transformations:
\[
\xymatrixrowsep{3pc}
\xymatrixcolsep{3pc}
\xymatrix{
	Q_i\times B^i
		\ar[d]_{1\times G^i}
		\ar[r]^{h'_i}
		\drtwocell\omit{^*{!(-1,-1.5)\object{\psi_i}}}
		& B \ar[d]^G \\
	Q_i\times A^i \ar[r]_{h_i}
		& A \\
}
\]
Let $\psi_i$ be given by
\def\epsiloncell#1{
	\ultwocell\omit{<-0.7>*{!(-2,-1.2){1\times \epsilon^{#1}}}}
}
\def\epsilonsumki{\epsiloncell{\sum k_i}}
\def\epsilonn{\epsiloncell{n}}
\[
\xymatrixrowsep{3.5pc}
\xymatrixcolsep{3pc}
\xymatrix{
	Q_i \times B^i
		\ar[d]_{1\times G^i}
		\ar[r]^{h'_i}
		\drtwocell\omit{^*{!(-1,-1.5)\object{\psi_i}}}
	& B
		\ar[d]^G \\
	Q_i \times A^i
		\ar[r]_{h_i}
	& A \\
}
\midequals
\xymatrixrowsep{1.5pc}
\xymatrixcolsep{1.2pc}
\xymatrix{
	Q_i \times B^i
		\ar[rr]^1
		\ar[dd]_{1\times G^i}
	&& Q_i \times B^i
		\ar[rr]^{h'_i}
	&& B
		\ar[dd]^G
		\\
	& {}
	\ultwocell\omit{<-0.7>*{!(-2,-1.8){1\times \epsilon^{i}}}}
	&& {}
		\lltwocell\omit{*{!(-5,-4)\object{{\pi^{-1}_i}}}}
		\\
	Q_i \times A^i
		\ar[urur]|{1 \times F^i}
		\ar[rr]_{h_i}
	&& A
		\ar[urur]^F
		\ar[rr]_1
	&& A
		\ultwocell\omit{{\eta}} \\
}
\]
We must check that $\psi$ satisfies (\ref{weakfunctordef}) and
(\ref{weakfunctordef2}) from Lemma \ref{lem:wkfunc_explicit}.
For (\ref{weakfunctordef}):
\begin{eqnarray*}
\midlabel{\mbox{LHS}} & \midequals
& \xymatrixrowsep{5pc}
\xymatrixcolsep{5pc}
\xymatrix{
	{}
		\ar[d]_{1\times G^{\sum k_i}}
		\ar[r]^{\prodkn{h'}}
		\drtwocell\omit{^*{!(-1,-1.5)\object{\prodkn{\psi}}}}
	& {}
		\ar[d]|{1\times G^n}
		\ar[r]^{h'_n}
		\drtwocell\omit{^*{!(-1,-1.5)\object{\psi_n}}}
	& {}
		\ar[d]^G
		\\
	{}
		\ar[r]_{\prodkn{h}}
	& {}
		\ar[r]_h
	& {} \\
} \\
& \midequals
& \xymatrixrowsep{2.3pc}
\xymatrixcolsep{2.3pc}
\xymatrix{
	{}
		\ar[rr]^1
		\ar[dd]_{1\times G^{\sum k_i}}
	&& {}
		\ar[rr]^{\prodkn{h'}}
	& {}
	& {}
		\ar[rr]^1
		\ar[dd]|{1\times G^n}
	&& {}
		\ar[rr]^{h'_n}
	& {}
	& {}
		\ar[dd]^G
		\\
	& {}
		\epsilonsumki
	&& {}
		\lltwocell\omit{*{!(-5,-4)\object{\prodkn{\pi^{-1}}}}}
	&& {}
		\epsilonn
	&& {}
		\lltwocell\omit{*{!(-5,-4)\object{\pi_n^{-1}}}}
	&& \\
	{}
		\ar[uurr]|{1 \times F^{\sum k_i}}
		\ar[rr]_{\prodkn{h}}
	&& {}
		\ar[uurr]|{1 \times F^n}
		\ar[rr]_1
	&& {}
		\ar[uurr]|{1 \times F^n}
		\ar[rr]_h
		\ultwocell\omit{*{!(1,0)\object{1 \times \eta^n}}}
	&& {}
		\ar[uurr]|F
		\ar[rr]_1
	&& {}
		\ultwocell\omit{{\eta}} \\
} \\
& \midequals
& \xymatrixrowsep{2.3pc}
\xymatrixcolsep{2.3pc}
\xymatrix{
	{}
		\ar[rr]^1
		\ar[dd]_{1\times G^{\sum k_i}}
	&& {}
		\ar[rr]^{\prodkn{h'}}
	& {}
	& {}
		\ar[rr]^1
	&& {}
		\ar[rr]^{h'_n}
	& {}
	& {}
		\ar[dd]^G
		\\
	& {}
		\epsilonsumki
	&& {}
		\lltwocell\omit{*{!(-5,-4)\object{\prodkn{\pi^{-1}}}}}
	&&&& {}
		\lltwocell\omit{*{!(-5,-4)\object{\pi_n^{-1}}}}
	& \\
	{}
		\ar[urur]|{1 \times F^{\sum k_i}}
		\ar[rr]_{\prodkn{h}}
	&& {}
		\ar[urur]|{1 \times F^n}
		\ar[rr]_1
	&& {}
		\ar[urur]|{1 \times F^n}
		\ar[rr]_h
		\uutwocell\omit{=}
	&& {}
		\ar[urur]|F
		\ar[rr]_1
	&& {}
		\ultwocell\omit{{\eta}} \\
} \\
& \midequals
& \xymatrixrowsep{2.3pc}
\xymatrixcolsep{2.3pc}
\xymatrix{
	{}
		\ar[rr]^1
		\ar[dd]_{1\times G^{\sum k_i}}
	&& {}
		\ar[rr]^{\prodkn{h'}}
	& {}
	& {}
		\ar[rr]^{h'_n}
	& {}
	& {}
		\ar[dd]^G
		\\
	& {}
		\epsilonsumki
	&& {}
		\lltwocell\omit{*{!(-5,-4)\object{\prodkn{\pi^{-1}}}}}
	&& {}
		\lltwocell\omit{*{!(-5,-4)\object{\pi_n^{-1}}}}
	&& \\
	{}
		\ar[urur]|{1 \times F^{\sum k_i}}
		\ar[rr]_{\prodkn{h}}
	&& {}
		\ar[urur]|{1 \times F^n}
		\ar[rr]_h
	&& {}
		\ar[urur]|F
		\ar[rr]_1
	&& {}
		\ultwocell\omit{{\eta}} \\
} \\
& \midequals
& \xymatrixrowsep{2.3pc}
\xymatrixcolsep{2.3pc}
\xymatrix{
	{}
		\ar[rr]^1
		\ar[dd]_{1\times G^{\sum k_i}}
	&& {}
		\ar[rr]^{h'_{\sum k_i}}
	& {}
	& {}
		\ar[dd]^G
		\\
	& {}
		\epsilonsumki
	&& {}
		\lltwocell\omit{*{!(-5,-4)\object{\pi_{\sum k_i}^{-1}}}}
	& \\
	{}
		\ar[urur]|{1 \times F^{\sum k_i}}
		\ar[rr]_{h_{\sum k_i}}
	&& {}
		\ar[urur]|F 
		\ar[rr]_1
	&& {}
		\ultwocell\omit{{\eta}} \\
} \\
& \midequals
& \xymatrixrowsep{5pc}
\xymatrixcolsep{5pc}
\xymatrix{
	{}
		\ar[d]_{1\times G^{\sum k_i}}
		\ar[r]^{h'_{\sum k_i}}
		\drtwocell\omit{^*{!(-1,-1.5)\object{\psi_{\sum k_i}}}}
	& {}
		\ar[d]^G
		\\
	{}
		\ar[r]_{h_{\sum k_i}}
	& {} \\
} \\
& = & \mbox{RHS.}
\end{eqnarray*}

For (\ref{weakfunctordef2}), consider the following diagram:
\[
\xymatrix{
	&&G b \ar[r]^{\delta_{1_Q}} \ar@/l3cm/[ddd]_{1} \ar[d]_{\eta G}
	  \ar @{} [dr]|*+[o][F-]{2} 
	  & h(1_P, G b) \ar@/r3cm/[ddd]^{\psi_1} \ar[d]^{\eta} \\
	\ar @{} [drr]|*+[o][F-]{1} 
	&&GFG b \ar[r]^{GF\delta_{1_Q}} \ar[d]_{1}
	  \ar @{} [dr]|*+[o][F-]{3} 
	  & GFh(1_P, Gb) \ar[d]^{\pi^{-1}_1}
	  \ar @{} [drr]|*+[o][F-]{5} 
	  \\
	&&GFG b \ar[r]_{G\delta'_{1_Q}} \ar[d]_{G\epsilon}
	  \ar @{} [dr]|*+[o][F-]{4} 
	  & G h'(1_P, FGb) \ar[d]^{Gh'(1_P, \epsilon)} &&\\
	&&Gb \ar[r]_{G\delta'_{1_Q}}
	  & Gh'(1_P, b)
}
\]
(\ref{weakfunctordef2}) is the outside of the diagram.
\xylabel{1} commutes by the triangle identities.
\xylabel{2} commutes by naturality of $\eta$.
\xylabel{3} commutes since $(F, \pi)$ is a $P$-functor.
\xylabel{4} commutes by naturality of $\delta$.
\xylabel{5} is the definition of $\psi$.
Hence the whole diagram commutes, and $(G, \psi)$ is a $P$-functor.

To see that $(F, G, \eta, \epsilon)$ is a $P$-equivalence, it is now enough to
show that $\eta$ and $\epsilon$ are $P$-transformations, since they satisfy the
triangle identities by hypothesis.

Write $(GF, \chi) = (G, \psi) \fcomp (F, \pi)$.
We wish to show that $\eta$ is a $P$-transformation $(1,1) \to (GF, \chi)$.
Each $\chi_{q, \adot}$ is the composite
\[
\xymatrix{
	h(q, GF\adot) \ar[r]^{\psi_{q, F\adot}}
	& G h(q, F\adot) \ar[r]^{G\pi_{q,\adot}}
	& GFh(q, \adot)
}
\]
Applying the definition of $\psi$, this is
\[
\centerline{
\xymatrix{
	h(q, GF\adot) \ar[r]^\eta
	& GF h(q, GF\adot) \ar[r]^{G\pi^{-1}}
	& G h(q, FGF\adot) \ar[r]^{Gh_q \epsilon F}
	& Gh(q, F \adot)   \ar[r]^{G\pi}
	& GFh(q, \adot)
}
}
\]

The axiom on $\eta$ is the outside of the diagram
\[
\centerline{
\xymatrix{
	h(q, \adot) \ar[rrrr]^1 \ar[dd]_{h(q, \eta)} \ar[dr]^\eta
	&&&& h(q, \adot) \ar[dd]^\eta \\
	& GFh(q, \adot) \ar[r]^{G\pi^{-1}} \ar[d]|{GFh(q, \eta)}
	\ar @{} [l]|*+[o][F-]{1} 
	\ar @{} [dr]|*+[o][F-]{2} 
	& Gh(q, F\adot) \ar[d]|{Gh(q, F\eta)} \ar[dr]^1
	\ar @{} [rr]|*+[o][F-]{3} & & \\ 
	h(q, GF\adot) \ar[r]^\eta
	& GFh(q, GF\adot) \ar[r]^{G\pi^{-1}}
	& Gh(q, FGF\adot) \ar[r]^{Gh(q, \epsilon F)}
	& Gh(q, F\adot) \ar[r]^{G\pi}
	& GFh(q, \adot)
}
}
\]
\xylabel{1} commutes by naturality of $\eta$, \xylabel{2} commutes by
naturality of $\pi^{-1}$, and \xylabel{3} commutes since $G\pi \fcomp G\pi^{-1}
= G(\pi \fcomp \pi^{-1}) = G1 = 1G$. The triangle commutes by the triangle
identities. So the whole diagram commutes, and $\eta$ is a $P$-transformation.
By Lemma \ref{inv <=> P-inv}, $\eta^{-1}$ is also a $P$-transformation.
Similarly, $\epsilon$ and $\epsilon^{-1}$ are $P$-transformations.

\end{proof}

The statement of the lemma is a fragment of the statement that \WkPCat\ is
2-monadic over \Cat.
Compare the fact that monadic functors reflect isos.

\begin{theorem}
\label{mainthm}
\index{coherence!for weak $P$-categories}
Let $Q\kel A \toletter{h} A$ be a weak $P$-category. Then $A$
is equivalent to $\st A$ in the 2-category \WkPCat.
\end{theorem}

\begin{proof}
Let $F: \st A \to A$ be given by $F(p, \adot) = h(p, \adot)$ and
identification of maps. This is certainly full and faithful, and it is
essentially surjective on objects because $\delta^{-1}_{1_Q}: h(1_P, a)
\to a$ is an isomorphism. By Lemma \ref{lem:equiv<=>Q-equiv}, it remains only
to show that $F$ is a weak $P$-functor.

We must find a sequence $(\phi_i : h_i (1 \times F^i) \to Fh')$ of natural
transformations satisfying equations (\ref{weakfunctordef}) and
(\ref{weakfunctordef2}) from Lemma \ref{lem:wkfunc_explicit}.
We can take $(\phi_i)_{q, (p\seq, a\dseq)} = (\delta_{q \cmp{p\seq}})_{a\dseq}$
for $q \in Q_n$ and $(p_1, a^1\seq), \dots, (p_n, a^n\seq) \in \st A$.
For (\ref{weakfunctordef}), we must show that

\[
\xymatrixrowsep{4pc}
\xymatrixcolsep{4pc}
\xymatrix{
	{}
		\ar[d]_{1\times F^{\sum k_i}}
		\ar[r]^{\prodkn{h'}}
		\drtwocell\omit{^*{!(-1,-1.5)\object{\prodkn{\phi}}}}
	& {}
		\ar[d]|{1\times F^n}
		\ar[r]^{h'_n}
		\drtwocell\omit{^*{!(-1,-1.5)\object{\phi_n}}}
	& {}
		\ar[d]^F
		\\
	{}
		\ar[r]_{\prodkn{h}}
	& {}
		\ar[r]_h
	& {} \\
}
\midequals
\xymatrixrowsep{4pc}
\xymatrixcolsep{4pc}
\xymatrix{
	{}
		\ar[d]_{1\times F^{\sum k_i}}
		\ar[r]^{h'_{\sum k_i}}
		\drtwocell\omit{^*{!(-1,-1.5)\object{\phi_{\sum k_i}}}}
	& {}
		\ar[d]^F
		\\
	{}
		\ar[r]_{h_{\sum k_i}}
	& {} \\
}
\]

All 2-cells in this equation are instances of $\delta$. Since there is at most
one 2-cell between two 1-cells in $Q$, the equation holds.

For (\ref{weakfunctordef2}) to hold, we must have 
\begin{equation}
\label{Fwkfctr}
\xymatrix{
	F(p,\adot) \ar[d]_1 \ar[r]^{\delta_{1_Q}}
	  & h(1_P, F(p,\adot)) \ar[d]^{\phi_{1_P}} \\
	F(p,\adot) \ar[r]_{F\delta'_{1_Q}}
	  & Fh'(1_P,(p, \adot))
}
\end{equation}
Since \st A is a strict monoidal category, $\delta' = 1$.
Apply this observation, and the definitions of $F$, $\phi$ and $h'$; then
(\ref{Fwkfctr}) becomes
\[
\xymatrix{
	h(p,\adot) \ar[d]_1 \ar[r]^{\delta_{1_Q}}
	  & h(1_P, h(p,\adot)) \ar[d]^{\delta_{1_P\cmp{p}}} \\
	h(p,\adot) \ar[r]_1
	  & h(p,\adot)
}
\]
Since there is at most one arrow between two 1-cells in $Q$, this diagram
commutes.
So $(F,\phi)$ is a weak $P$-functor, and hence (by Lemma
\ref{lem:equiv<=>Q-equiv}) an equivalence in \WkPCat.
\end{proof}

\begin{example}
\index{sets}
Consider the initial operad $0$, whose algebras are sets.
We saw in Example \ref{ex:trivial_theory} that unbiased weak $0$-categories
are categories equipped with a trivial monad.
By Theorem \ref{mainthm}, every unbiased weak $0$-category is equivalent via
weak $0$-functors to a category equipped with a monad which \emph{is} the
identity: in other words, a category.
\end{example}

\begin{example}
\index{monoids}
Consider the terminal operad $1$, whose algebras are monoids.
Theorem \ref{mainthm} tells us that every unbiased weak monoidal category is
monoidally equivalent to a strict monoidal category.
\end{example}

\section{Universal property of \stfunc}
\label{stsect}
\index{strictification!universal property}

Let $P$ be a plain operad.
\begin{theorem}
Let $U'$ be the forgetful functor $\StrPCat \to \WkPCat$ (considering both of
these as 1-categories). Then \stfunc is left adjoint to $U'$.
\end{theorem}
\begin{proof}
For each $(A,h) \in \WkPCat$, we construct an
initial object $A \toletter{(F', \psi)} \st A$ of the comma category
$(A\downarrow U')$, thus showing that \stfunc is functorial and that $\stfunc
\dashv U'$ (and that $(F', \psi)$ is the component of the unit at $A$).  Let
$(B, h'')$ be a strict $P$-category, and
$(G, \gamma): A \to U'B$ be a weak $P$-functor. We must show that there is a
unique strict $P$-functor $H$ making the following diagram commute:
\begin{equation}
\label{A->U}
\xymatrix{
	& A \ar[ddl]_{(F', \psi)} \ar[ddr]^{(G, \gamma)} \\ \\
	U'\,\st A \ar@{-->}[rr]^{(H, {\rm id})} & & U'B
}
\end{equation}
$(F', \psi)$ is given as follows:
\begin{itemize}
\item If $a \in A$, then ${F'}(a) = (1, a)$.
\item If $f: a \to a'$ in $A$ then ${F'}f$ is the lifting of $h(1,f)$ with
source $(1, a)$ and target $(1, a')$.
\item Each $\psi_{(p, \adot)}$ is the lifting of $(\delta_{1_Q})_{h(p, \adot)}:
h(p, \adot) \to h(1, h(p, \adot))$ to a morphism
$h'(p, F'(a)\seq) = (p, \adot) \to (1, h(p, \adot)) = F'(h(p, \adot))$.
\end{itemize}
For commutativity of (\ref{A->U}), we must have $H(1, a) = G(a)$, and
for strictness of $H$, we must have $H(p, \adot) = h''(p, H(1, a)\seq)$. These
two conditions completely determine $H$ on objects.

Now, take a morphism $f:(p, \adot) \to (p', \adot')$, which is a lifting of a
morphism $g: h(p,\adot) \to h(p',\adot')$ in $A$.
Then $Hf$ is a morphism $h''(p, G\adot) \to h''(p', G\adot')$: the obvious
thing for it to be is the composite
\[
\xymatrix{
	h''(p, G\adot) \ar[r]^\gamma 
	& G h''(p, \adot) \ar[r]^{Gg}
	& G h''(p', \adot') \ar[r]^{\gamma^{-1}}
	& h''(p', G\adot')
}
\]
and we shall show that this is in fact the only possibility.
Consider the composite
\[
\xymatrix{
	(1, h(p, \adot)) \ar[r]^-{\psi^{-1}}
	&(p, \adot) \ar[r]^{f}
	& (p', \adot) \ar[r]^-{\psi}
	& (1, h(p', \adot'))
}
\]
in \st A. Composition in \st A is given by composition in $A$, so this is
equal to the lifting of $\delta_{1_Q} \fcomp g \fcomp \delta_{1_Q}^{-1} = h(1,
g)$ to a morphism $(1, h(p, \adot)) \to (1, h(p',
\adot'))$, namely $F'g$. So $f = \psi^{-1} \fcomp F'g \fcomp \psi$, and
$Hf = H\psi^{-1} \fcomp HF'g \fcomp H\psi$. By commutativity of
(\ref{A->U}), $HF' = G$ and $H\psi = \gamma$, so $Hf = \gamma^{-1} \fcomp
Gg \fcomp \gamma$ as required.

This completely defines $H$. So we have constructed a unique
$H$ which makes (\ref{A->U}) commute and which is strict. Hence $(F',
\psi): A \to U'\,\st A$ is initial in $(A\downarrow U')$, and so $\stfunc
\dashv U'$.
\end{proof}

The $P$-functor $(F, \phi):\st A \to A$ constructed in Theorem \ref{mainthm} is
pseudo-inverse to $(F', \psi)$, which we have just shown to be the
$A$-component of the unit of the adjunction \stfunc $\dashv U'$. We can
therefore say that \StrPCat\ is a weakly coreflective sub-2-category of
\WkPCat. Note that the counit is \emph{not} pseudo-invertible, so this is not
a 2-equivalence.

\begin{example}
\index{sets}
Consider again the initial operad $0$, whose algebras are sets.
We saw in Example \ref{ex:trivial_theory} that unbiased weak $0$-categories
are categories equipped with a specified trivial monad.
Let $\Triv$ denote the category of such categories, with morphisms being
functors that preserve the trivial monad up to coherent isomorphism.
A strict unbiased $0$-category is a category equipped with a monad equal to
the identity monad, which is simply a category.
So $\Cat$ is a weakly coreflective sub-2-category of $\Triv$.
\end{example}

\section{Presentation-independence}
\label{sec:presindep}

We will now show that the weakening of a symmetric operad $P$ is essentially
independent of the generators chosen.
This generalizes Leinster's result (in \cite{hohc} section 3.2) that the
theory of weak monoidal categories is essentially unaffected by the choice of a
different presentation for the theory of monoids.

We will need the following lemma:
\begin{lemma}
\label{lem:fork}
In \CatSymmOperad, if $\fork P \alpha \beta Q \gamma R$ is a fork, and $\gamma$
is levelwise full and faithful, then $\alpha \cong \beta$.
\end{lemma}
\begin{proof}
We shall construct an invertible \Cat-$\Sigma$-operad transformation $\eta:\alpha \to
\beta$.
We form the $\eta_n$s as follows: for all $p \in P_n$, let $\gamma\alpha(p) =
\gamma\beta(p)$.
Since $\gamma$ is levelwise full, there exists an arrow $(\eta_n)_p: \alpha(p)
\to \beta(p)$ such that $\gamma_n((\eta_n)_p) = 1_{\gamma\alpha(p)}$.
Since $\gamma$ is levelwise full and faithful, this arrow is an isomorphism.
Each $\eta_n$ is natural because, for all $n \in \natural$ and $f: p \to q$ in
$P_n$,  the image under $\gamma$ of the naturality square
\[
\xymatrix{
	\alpha(p) \ar[r]^{(\eta_n)_p} \ar[d]_{\alpha(f)}
	& \beta(p) \ar[d]^{\beta(f)} \\
	\alpha(q) \ar[r]^{(\eta_n)_q}
	& \beta(q)
}
\]
is
\[
\xymatrix{
	\gamma\alpha(p) \ar[r]^1 \ar[d]_{\gamma\alpha(f)}
	& \gamma\beta(p) \ar[d]^{\gamma\beta(f)} \\
	\gamma\alpha(q) \ar[r]^1
	& \gamma\beta(q)
}
\]
which commutes since $\gamma\alpha = \gamma\beta$.
Since $\gamma$ is faithful, the naturality square commutes, and $\eta_n$ is
natural.
It remains to show that the collection $(\eta_n)_{n \in \natural}$ forms a
\Cat-$\Sigma$-operad transformation, in other words that the equations
{\def\nxseq#1{{#1}_n \times {#1}\seq}
\def\sumki{_{\sum{k_i}}}
\begin{eqnarray}
\xymatrixrowsep{5pc}
\xymatrixcolsep{5pc}
\xymatrix{
	\nxseq P \rtwocell^{\nxseq{\alpha}}_{\nxseq{\beta}}
		{*{!(-3.5,0)\object{\nxseq{\eta}}}}
		\ar[d]_\ocomp
		\drtwocell\omit{=}
	& \nxseq Q \ar[d]^\ocomp \\
	P\sumki \ar[r]_{\beta\sumki} & Q\sumki
}
&\midequals
&\xymatrixrowsep{5pc}
\xymatrixcolsep{5pc}
\xymatrix{
	\nxseq P \ar[r]^{\nxseq{\alpha}}
		\ar[d]_\ocomp
		\drtwocell\omit{=}
	& \nxseq Q \ar[d]^\ocomp \\
	P\sumki \rtwocell^{\alpha\sumki}_{\beta\sumki}
		{*{!(-3.5,0)\object{\eta\sumki}}}
	& Q\sumki
} \\
(\eta_1)_{1} &= &1 \\
\xymatrixrowsep{4pc}
\xymatrixcolsep{4pc}
\xymatrix{
	P_n \rtwocell^{\alpha_n}_{\beta_n}{*{!(-1,0)\object{\eta_n}}}
		\ar[d]_{\sigma\cdot-}
		\drtwocell\omit{=}
	& Q_n \ar[d]^{\sigma\cdot-} \\
	P_n \ar[r]_{\beta_n} & Q_n
}
&\midequals
&\xymatrixrowsep{4pc}
\xymatrixcolsep{4pc}
\xymatrix{
	P_n \ar[r]^{\alpha_n}
		\ar[d]_{\sigma\cdot-}
		\drtwocell\omit{=}
	& Q_n \ar[d]^{\sigma\cdot-} \\
	P_n \rtwocell^{\alpha_n}_{\beta_n}{*{!(-1,0)\object{\eta_n}}} & Q_n
}
\end{eqnarray}}
hold, for all $n, k_1 \dots k_n \in \natural$ and every $\sigma \in S_n$.
As above, it is enough to show that the images of both sides under $\gamma$
are equal, and this is trivially true by definition of $\eta$.
\end{proof}

Let $P$ be a symmetric operad.
\begin{theorem}
\label{thm:pres_indep}
\index{presentation!independence of weakenings}
Let $\Phi \in \Set^\natural$ and let $\phi : \Fsig \Phi \to P$ be a regular epi.
Then $\Wkwrt P \phi$ is equivalent as a symmetric \Cat-operad to \WkP.
\end{theorem}
\begin{proof}
Let $Q$ be the weakening of $P$ with respect to $\phi: \Fsig \Phi \to P$.
By the triangle identities, we have a commutative square
\[
\xymatrix{
	\Fsig \Phi \ar[r]^\phi \ar[d]_{\Fsig \trans\phi}
	& P \ar[d]^1 \\
	\Fsig \Usig P \ar[r]^{\epsilon_P}
	& P
}
\]
By functoriality of the factorization system, this gives rise to a unique map
$\chi: Q \to \WkP$ such that 
\[
\xymatrix{
	\Fsig \Phi \booar[r] \ar@/u0.7cm/[rr]^\phi \ar[d]_{\Fsig \trans\phi}
	& Q \lffar[r] \unar[d]^\chi
	& P \ar[d]^1 \\
	\Fsig \Usig P \booar[r] \ar@/d0.7cm/[rr]_{\epsilon_P}
	& \WkP \lffar[r]
	& P
}
\]
commutes.
We wish to find a pseudo-inverse to $\chi$.

Since \SymmOperad\ is monadic over $\Set^\natural$, a regular epi in
\SymmOperad\ is a levelwise surjection by Theorem
\ref{thm:reg epis split in Set^X^T}.
So we may choose a section $\psi_n$ of $\phi_n : (\Fsig \Phi)_n \to P_n$ for
all $n \in \natural$.
So we have a morphism $\psi : \Usig P \to \Usig \Fsig \Phi$ in $\Set^\natural$.
We wish to show that
\[
\xymatrix{
	\Fsig \Usig P \ar[r]^{\epsilon_P} \ar[d]_{\trans\psi}
	& P \ar[d]^1 \\
	\Fsig \Phi \ar[r]^\phi
	& P
}
\]
commutes.
This follows from a simple transpose argument:
\[
\newdir{((}{{}*!/-5pt/{}}
\xymatrixrowsep{0.15pc}
\xymatrix{
& \Fsig \Usig P \ar[r]^{\bar \psi} & \Fsig \Phi \ar[r]^\phi & P \\
\ar@{((-((}[rrrr] & & & & \\
& \Usig P \ar[r]^\psi & \Usig \Fsig \Phi \ar[r]^{\Usig \phi} & \Usig P & =
& \Usig P \ar[r]^1 & \Usig P \\
& & & & \ar@{-}[rrr] & & & \\
& & & & & \Fsig \Usig P \ar[r]^\epsilon & P.
}
\]
This induces a map
\[
\xymatrix{
	\Fsig \Usig P \booar[r] \ar@/u0.7cm/[rr]^{\epsilon_P}
		\ar[d]_{\trans\psi}
	& \WkP \lffar[r] \unar[d]^\omega
	& P \ar[d]^1 \\
	\Fsig \Phi \booar[r] \ar@/d0.7cm/[rr]_\phi
	& Q \lffar[r]
	& P
}
\]
We will show that $\omega$ is pseudo-inverse to $\chi$.
Now,
\[
\xymatrix{
	Q \lffar[r] \ar[d]_\omega
	& P \ar[d]^1 \\
	\WkP \lffar[r] \ar[d]_\chi
	& P \ar[d]^1 \\
	Q \lffar[r]
	& P
}
\]
commutes.
So $\xymatrix{Q\parallelars{1_Q}{\chi\omega} & Q \lffar[r] & P}$ is a fork.
By Lemma \ref{lem:fork}, $\chi\omega \cong 1_Q$, and similarly $\omega\chi
\cong 1_{\WkP}$.
So $Q \simeq \WkP$ as a symmetric \Cat-operad, as required.
\end{proof}

\begin{corollary}
\label{cor:unbiased_symm}
\index{unbiased!weakening of a plain operad!comparison to symmetric weakening}
Let $P$ be a plain operad.
Then $\Fsigp (\WkP) \simeq \Wk{\Fsigp P}$.
\end{corollary}
\begin{proof}
Let $\phi: \Fp\Up P \to P$ be the component at $P$ of the counit of the
adjunction $\Fp \dashv \Up$.
Let $\epsilon$ be the counit of the adjunction $\Fsigp \dashv \Usigp$.

By Theorem \ref{thm:plain_as_symmetric}, there is an isomorphism $\Fsigp
(\Wkwrt P \phi) \cong \Wkwrt{\Fsigp P}{\Fsigp \phi}$, and by Theorem
\ref{thm:pres_indep}, there is an equivalence $\Wkwrt{\Fsigp P}{\Fsigp \phi}
\simeq \Wk{\Fsigp P}$.
Hence $\Fsigp (\Wk P) \simeq \Wk{\Fsigp P}$.
\end{proof}

\begin{corollary}
Let $P$ be a plain operad.
Then the category \WkPCat\ is equivalent to the category
\cat{Wk-$\Fsigp P$-Cat}.
\end{corollary}

This tells us that the unbiased categorification of a strongly regular theory
is essentially unaffected by our treating it as a linear theory instead.

\begin{example}
\index{sets}
Considering again the trivial theory $0$, we see that $\Triv \simeq \Cat$.
\end{example}

This can be generalised to the multi-sorted situation:

\begin{lemma}
Let $X$ be a set, and $f$ be a regular epi in the category $\CatMulticat_X$ or
in the category $\CatSymmMulticat_X$.
Then $f$ is locally surjective on objects.
\end{lemma}
\begin{proof}
$\Multicat_X$ is monadic over $\Multigraph_X$ by Lemma
\ref{lem:frees_exist} and Theorem \ref{lem:monadj}, and an object of
$\Multigraph_X$ can be considered as an object of $\Set^Y$, where $Y
= X \times X^*$, and $X^*$ is the free monoid on $X$: for each $x \in X$, and
each sequence $x_1, \dots, x_n \in X^*$, there is a set of funnels $x_1, \dots,
x_n \to x$.
Hence, by \ref{thm:reg epis split in Set^X^T}, every regular epi in
$\Multicat_X$ is locally surjective.

The objects functor $O: \Cat \to \Set$ has both a left adjoint $D$
and a right adjoint $I$.
Hence $O$ and $I$ preserve products, and hence by Lemma \ref{lem:frees_exist}
they induce an adjunction
\[
\adjunction{\CatMulticat_X}{\Multicat_X}{O_*}{I_*}.
\]
Since $O_*$ is a left adjoint, it preserves colimits, and in particular regular
epis: hence, every regular epi in $\CatMulticat_X$ must be locally surjective
on objects.

The symmetric case is proved analogously.
\end{proof}

\begin{theorem}
Let $M$ be a (symmetric) multicategory, and $\phi : \Fp\Phi \to M$
(or in the symmetric case, $\phi : \Fsig\Phi \to M$) be a regular epi.
Then the weakening of $M$ with respect to $\phi$ is equivalent as a
\Cat-multicategory to $\Wk M$.
\end{theorem}
\begin{proof}
The proof is exactly as for Theorem \ref{thm:pres_indep}.
\end{proof}

\chapter{Other Approaches}
\label{ch:others}

\section{Pseudo-algebras for 2-monads}

We begin by recalling some standard notions of 2-monad theory.
\begin{defn}
\index{2-monad}
A \defterm{2-monad} is a monad object in the 2-category of 2-categories, in the
sense of \cite{street}; that is to say, a 2-category $\C$, a strict 2-functor
$T: \C \to \C$, and 2-transformations $\mu : T^2 \to T$ and $\eta : 1_\C \to T$
satisfying the usual monad laws:
\begin{equation}
\xymatrix{
	& T^3 \ar[dr]^{\mu T} \ar[dl]_{T \mu} \\
	T^2 \ar[dr]_{\mu} & & T^2 \ar[dl]^\mu \\
	& T
}
\end{equation}
\begin{equation}
\xymatrix{
	T \ar[r]^{T\eta} \ar[dr]_{1_T}
	& T^2 \ar[d]^\mu
	& T \ar[l]_{\eta T} \ar[dl]^{1_T} \\
	& T
}
\end{equation}
\end{defn}

As is common for ordinary 1-monads, we will usually refer to a 2-monad $(\C, T,
\mu, \eta)$ as simply $T$.

The usual notion of an algebra for a monad carries over simply to this case:

\begin{defn}
\index{algebra!for a 2-monad}
Let $(\C, T, \mu, \eta)$ be a 2-monad.
A \defterm{strict algebra} for $T$ is an object $A \in \C$ and a 1-cell $a : T
A \to A$ satisfying the following axioms:
\begin{equation}
\xymatrix{
& T^2 A \ar[dr]^{\mu} \ar[dl]_{T a} \\
	T A \ar[dr]_a & & T A \ar[dl]^a \\
	& A
}
\end{equation}
\begin{equation}
\xymatrix{
	A \ar[r]^{\eta} \ar[dr]_{1_A}
	& T A \ar[d]^a \\
	& A
}
\end{equation}
\end{defn}

For our purposes, it is more interesting to consider the well-known 
\defterm{pseudo-algebras} for a 2-monad.
These are algebras ``up to isomorphism'':

\begin{defn}
\index{pseudo-algebra!for a 2-monad}
Let $(\C, T, \mu, \eta)$ be a 2-monad.
A \defterm{pseudo-algebra} for $T$ is an object $A \in \C$, a 1-cell $a: TA
\to A$, and invertible 2-cells
\[
\xymatrixnocompile{
	& T^2 A \ar[dr]^{\mu} \ar[dl]_{T a}
	\ddtwocell\omit{^*{!(0, -0.5)\object{\alpha}}} \\
	T A \ar[dr]_a & & T A \ar[dl]^a \\
	& A
}
\hskip 2cm
\xymatrixnocompile{
	A \ar[r]^{\eta} \ar[dr]_{1_A}
	\drtwocell\omit{^<-1.7>*{!(0.5,0.5)\object{\beta}}}
	& T A \ar[d]^a \\
	& A
}
\]
satisfying the equations
\begin{equation}
\xymatrixnocompile{
& T^3 A \ar[dl]_{T^2 a} \ar[dr]^{\mu T} \ar[dd]^{T \mu}
& & & & T^3 A \ar[dl]_{T^2 a} \ar[dr]^{\mu T}  \ddtwocell\omit{=} \\
T^2 A \ar[dd]_{Ta}
& {} \drtwocell\omit{=} & T^2 A \ar[dd]^{\mu}
& & T^2 A \ar[dd]_{Ta} \ar[dr]_{\mu}
& & T^2 A \ar[dd]^{\mu} \ar[dl]^{Ta} \\
\urtwocell\omit{*{!(-1, -1.5)\object{T \alpha}}}
& T^2 A \ar[dl]_{T a} \ar[dr]_{\mu}
& {} & = 
& & T A \ar[dd]^a
& \\ 
T A \ar[dr]_a
& & T A \ar[dl]^a
& & T A \ar[dr]_a
& \ultwocell\omit{\alpha} \urtwocell\omit{\alpha}
& T A \ar[dl]^a \\
& A \uutwocell\omit{\alpha}
& & & & A
}
\end{equation}
\begin{equation}
\xymatrix{
& T A \ar@/u0.8cm/[ddl]_1 \ar@/u0.8cm/[ddr]^1 \ar[d]|{T \eta}
& & & TA \ar@/l0.5cm/[dd]_a \ar@/r0.5cm/[dd]^a \ddtwocell\omit{=}
& & & TA \ar[dl]_a \ar@/ur1cm/[ddrr]^1 \ar[dr]^{\eta T}  \ddtwocell\omit{=} \\
& T^2 A \ar[dl]_{T a} \ar[dr]^\mu 
& & =
& & =
& A \ar[dr]^\eta \ar@/dl1cm/[ddrr]_1
  \drlowertwocell\omit{^<1.7>*{!(-0.5,-0.5)\object{\beta}}}
& & T^2 A \ar[dr]^\mu \ar[dl]_{T a}\\
TA \ar[dr]_a
& & TA \ar[dl]^a
& & A
& & & TA \ar[dr]_a
& & TA \ar[dl]^a\\
& A \uutwocell\omit{\alpha}
& & & & & & & A \uutwocell\omit{\alpha}
}
\end{equation}
\end{defn}

\begin{defn}
\index{pseudo-algebra!for a 2-monad!pseudo-morphisms}
Let $(\C, T, \mu, \eta)$ be a 2-monad, and let $(A, a, \alpha_1, \alpha_2)$ and
$(B, b, \beta_1, \beta_2)$ be pseudo-algebras for $T$.
A \defterm{pseudo-morphism of pseudo-algebras} $(A, a, \alpha_1, \alpha_2)$ to
$(B, b, \beta_1, \beta_2)$ is a pair $(f, \phi)$, where $f : A \to B$ is a
1-cell in $\C$ and $\phi$ is an invertible 2-cell:
\[
\xymatrix{
TA \ar[r]^{Tf} \ar[d]_a \drtwocell\omit{\phi}
& TB \ar[d]^b \\
A \ar[r]_f & B
}
\]
satisfying the axioms
\begin{eqnarray}
\xymatrixrowsep{3pc}
\xymatrixcolsep{3pc}
\xymatrix{
	T^2 A \ar[r]^{T^2 f} \ar[d]_{\mu_A}\drtwocell\omit{=}
	& T^2 B \ar[d]^{\mu_B} \\
	T A \ar[r]^{Tf} \ar[d]_a \drtwocell\omit{\phi}
	& T B \ar[d]^b \\
	A \ar[r]_f
	& B
}
\xymatrix{ {} \\ = }
\xymatrix{
	T^2 A \ar[r]^{T^2 f} \ar[d]_{T a}
	\drtwocell\omit{*{!(-2,0)\object{T \phi}}}
	& T^2 B \ar[d]^{T b} \\
	T A \ar[r]^{Tf} \ar[d]_a \drtwocell\omit{\phi}
	& T B \ar[d]^b \\
	A \ar[r]_f
	& B
}
\\
\xymatrixrowsep{3pc}
\xymatrixcolsep{3pc}
\xymatrix{
	A \ar[r]^f \ar[d]_{\eta_A}
		\ddlowertwocell<-13>_{1_A}{*{!(0,1)\object{\alpha_2^{-1}}}}
		\drtwocell\omit{=}
	& B \ar[d]^{\eta_B}
		\dduppertwocell<13>^{1_B}{*{!(0,1)\object{\beta_2^{-1}}}} \\
	T A \ar[r]^{T f} \ar[d]_a \drtwocell\omit{\phi}
	& T B \ar[d]^b \\
	A \ar[r]_f
	& B
}
\xymatrix{ {} \\ = }
\xymatrix{
	\\
	A \rtwocell^f_f{*{!(-0.5,0)\object{1_f}}} & B
}
\end{eqnarray}
\end{defn}

This gives rise to a category $\PsAlg{T}$ for any 2-monad $T$.
\index{$\PsAlg{T}$}

Every cartesian monad $T$ on \Set\ gives rise to a 2-monad $\bar T$ on \Cat\ in
an obvious way, and (as we saw in Theorem \ref{thm:opds<=>cart. monad}) every
plain operad $P$ gives rise to a cartesian monad $T_P$ on \Set.
So an alternative definition of ``weak $P$-category'' might be ``pseudo-algebra
for $\bar T_P$''.

In order to explore the connections between this idea and the notion of weak
$P$-category given in previous chapters, we shall need some theorems from
\cite{bkp} and related papers.

\begin{theorem} (Blackwell, Kelly, Power)
\label{thm:ps-classifier}
Let $T$ be a 2-monad with rank on a cocomplete 2-category $K$, let $\Alg T\str$
be the 2-category of strict $T$-algebras and strict morphisms, and $\Alg T\wk$
be the 2-category of strict $T$-algebras and weak morphisms.
Then the inclusion $J : \Alg T\str \to \Alg T\wk$ has a left adjoint $L$.
Thus every strict $T$-algebra $A$ has a \defterm{pseudo-morphism classifier} $p
: A \to A'$ (where $A' = JLA$), such that for all $B \in K$, and every
pseudo-morphism $f: A \to B$, we may express $f$ uniquely as the composite of
$p$ and a strict morphism:
\[
\xymatrix{
	A \ar[r]^p \ar[dr]_f & A' \unar[d]^{J\bar f} \\ & B
}
\]
\end{theorem}
\begin{proof}
See \cite{bkp}, Theorem 3.13.
\end{proof}

\begin{theorem} (Blackwell, Kelly, Power)
\label{thm:f* adj}
Let $f:S \to T$ be a strict map between 2-monads with rank on a cocomplete
2-category $K$.
Then the induced map $f^* : \Alg T\str \to \Alg S\str$ has a left adjoint, and 
the induced map $f^* : \Alg T\wk \to \Alg S\wk$ has a left biadjoint.
\end{theorem}
\begin{proof}
See \cite{bkp}, Theorem 5.12.
\end{proof}

\begin{corollary}
\label{cor:ps-alg adj}
Composing this left adjoint with the left adjoint of Theorem
\ref{thm:ps-classifier} gives us an adjunction
\[
\NoCompileMatrices
\adjunction{\Alg S\wk}{\Alg T\str} F U
\]
\end{corollary}

\begin{theorem} (Power, Lack)
\label{thm:ps-classifier as boff}
Let $T$ be a 2-monad with rank on a cocomplete 2-category $K$ of the form
$\Cat^X$ for some set $X$, and let $T$ preserve pointwise
bijectivity-on-objects.
Let $(A, a)$ be a strict $T$-algebra.
Then the pseudo-morphism classifier $A'$ for $A$ may be found by factorizing the
structure map $a : TA \to A$ as a pointwise bijective-on-objects map followed
by a locally full-and-faithful map:
\[
\xymatrix{
	TA \ar[rr]^a \booar[dr] & & A \\ & A' \lffar[ur]
}
\]
\end{theorem}
\begin{proof}
The construction is given in Power's paper \cite{power}, and the universal
property of the algebra constructed is proved in Lack's paper \cite{lack}.
\end{proof}
This argument is due to Steve Lack (private communication).

\begin{theorem}
\label{thm:ps-algs are weak algs}
Let $P$ be a plain operad.
Let $T_P$ be the monad induced by $P$ on \Set.
Then a pseudo-algebra for $\bar T_P$ is a weak $P$-category in the sense of
Definition \ref{def:weakening}.
Furthermore, there is an isomorphism of categories $\PsAlg{\bar T_P} \cong
\Algwk{P}$.
\end{theorem}
\begin{proof}
\CatOperad\ is monadic over $\Cat^\natural$ via one of the special monads of
Theorem \ref{thm:ps-classifier as boff}, and hence, for every plain
$\Cat$-operad $P$, the pseudo-morphism classifier of $P$ is none other than
$\WkP$.
Hence, if $A$ is a category, then a strict map of \Cat-operads $\WkP \to
\End(A)$ is precisely a weak map $P \to \End(A)$, or equivalently a
$\bar T_P$-pseudo-algebra structure on $A$.
\end{proof}

We may also use these ideas to provide a simple proof of the strictification
result in Theorem \ref{mainthm}.
The map $P \to \WkP$ given by Theorem \ref{thm:ps-classifier} is pseudo, but it
has a strict retraction $q: \WkP \to P$.
This is equivalent to a strict map of monads $T_\WkP \to T_P$.
By Corollary \ref{cor:ps-alg adj}, this induces a 2-functor $\Alg P\str \to \Alg
\WkP\wk$ with a left adjoint.
This functor is simply the inclusion of the 2-category of strict
$P$-categories, strict $P$-functors and $P$-transformations into the 2-category
of weak $P$-(categories, functors, transformations), and its left adjoint is
the functor $\stfunc$ constructed in Section \ref{sec:strictification}.
The fact that any weak $P$-category $A$ is equivalent to $\st A$ is a
consequence of the fact that any pseudo $P$-algebra is equivalent to a strict
one, and this holds by the General Coherence Result of Power.

However, pseudo-algebras are less useful in the case of linear theories.
Since the monads arising from symmetric operads are not in general cartesian,
we may not perform the construction given above.
We may, however, use the existence of colimits in \Cat, and consider the
2-monad
\[
A \mapsto \int^{n \in \bbB} P_n \times A^n
\]
for any symmetric operad $P$.
If $P$ is the free symmetric operad on a plain operad $P'$, this 2-monad is
equal to $\bar T_{P'}$.
Yet this coend construction also leads to problems.

Let $T$ be the ``free commutative monoid'' monad on \Set, and $S$ be the ``free
monoid'' monad on \Set.
Since these both arise from symmetric operads, we may lift them to 2-monads
$T', S'$ on \Cat\ as described above.
$T'$ is the free commutative monoid monad on \Cat, which is to say the free
strict symmetric monoidal category 2-monad; similarly, $S'$ is the free strict
monoidal category 2-monad.
For each category $A$, there is a functor $\pi_A: S'A \to T'A$ which is full
and surjective-on-objects; hence, if $(A, a, \alpha_1,\alpha_2)$ is a
pseudo-algebra for $T'$, we obtain an $S'$-pseudo-algebra structure by
precomposing with
$\pi_A$:
\[
\xymatrix{
	S' A \ear[d]^{\pi_A} \\
	T' A \ar[d]^a \\
	A
}
\hskip 1cm
\xymatrix{
	& S'^2 A \ear[d]^{(\pi * \pi)_A} \\
	& T'^2 A \ar[dl]_{T'a} \ar[dr]^{\mu} \ddtwocell\omit{^\alpha_1} \\
	T'A \ar[dr]_a & & T'A \ar[dl]^a \\
	& A
}
\hskip 1cm
\xymatrix{
	A \ar[d]_{1_A} \ar[r]^\eta & S'A \ear[d]^{\pi_A} \\
	A \ar[r]^\eta \ar[dr]_{1_A}
	\druppertwocell\omit{^<-2>*{!(1,0)\object{\alpha_2}}}
	& T'A \ar[d]^a \\
	& A
}
\]
The $S'$-pseudo-algebra structure so obtained is uniquely determined.
Since $\pi_A$ is full and surjective-on-objects, it is epic, and so every
pseudo-algebra for $S'$ is a pseudo-algebra for $T'$ in at most one way.
Hence we may view all pseudo-algebras for $T'$ as pseudo-algebras for $S'$
(that is, as monoidal categories) with extra properties.
But there exist monoidal categories with several choices of symmetric structure on them.
For instance, consider the category of graded Abelian groups, with tensor
product $(A\otimes B)_n = \bigoplus_{i + j = n} A_i \otimes B_j$.
As well as the obvious symmetry, there is another given by $\tau_{AB}(a\otimes
b) = (-1)^{ij} b\otimes a$, where $a \in A_i, b \in B_j$.

We can say at least something about the extra properties that pseudo-algebras
for $T'$ must have:

\begin{theorem}
\index{commutative monoids}
A pseudo-algebra for $T'$ is a symmetric monoidal category $A$ in which $x
\otimes y = y \otimes x$ for all $x, y \in A$.
\end{theorem}
\begin{proof}
Recall our construction of the \dop\ whose algebras are commutative monoids in
Example \ref{ex:comm_monoid_fp}.
From this, we may deduce that if $A$ is a set, then an element of $T_P A$ is a
function $A \to \natural$ assigning each element of $A$ its multiplicity: in
other words, a multiset of elements of $A$.
Let $(A, a, \alpha, \beta)$ be a pseudo-algebra for $T'$ in \Cat.
Then we have a binary tensor product:
\[
x\otimes y := a(x^1y^1)
\]
where $x^1y^1$ is the function $A \to \natural$ sending $x$ and $y$ to 1 and
all other objects of $A$ to 0.
The tensor is defined analogously on morphisms.
The components of $\alpha$ and $\beta$ give us associator, symmetry and unit
maps, and it can be shown that they satisfy the axioms for a monoidal category.
However, since the function $x^1y^1$ is equal to the function $y^1x^1$ for all
$x, y \in A$, it must be the case that $x \otimes y = y \otimes x$.
\end{proof}

Since not all symmetric monoidal categories satisfy this condition, it is
apparent that a na\"\i ve approach to categorification based on pseudo-algebras
is doomed to fail, and that more sophistication is required.
In fact, I conjecture that a stronger condition holds: that the symmetry maps
are all identities.

In the specific case of symmetric monoidal categories, we may remedy the
situation as follows. \index{symmetric monoidal category}
Let $T$ be the ``free symmetric strict monoidal category'' 2-monad.
Then pseudo-algebras for $T$ are precisely symmetric monoidal categories.

\section{Laplaza sets}

This notion was introduced by T. Fiore, P. Hu and I. Kriz in
\cite{fhk_laplaza}, as a generalization of Laplaza's categorification of rigs
in \cite{laplaza}.
It was introduced as an attempt to correct an error in the earlier definition
of categorification proposed in \cite{fiore}; the error in question is
essentially that discussed in Section \ref{sec:evaluation} above.
\index{Fiore, Tom}\index{Hu, Po}\index{Kriz, Igor}\index{Laplaza, Miguel}

\begin{defn}
\index{Laplaza set}
Let $T$ be a \dop.
A \defterm{Laplaza set} for $T$ is a subsignature of $\Ufp T$.
\end{defn}

Concretely, a Laplaza set $S$ for $T$ is a sequence $S_0 \subset T_0, S_1
\subset T_1, \dots$ of subsets of $T_0, T_1, \dots$.

\begin{defn}
\index{pseudo-algebra!for a \dop\ with Laplaza set}
Let $T$ be a \dop, and $S$ be a Laplaza set for $T$.
A \defterm{$(T,S)$-pseudo algebra} is
\begin{itemize}
\item a category \C\,
\item for each $\phi \in T_n$, a functor $\hat \phi : \C^n \to \C$,
\item \defterm{coherence morphisms} witnessing all equations that are true in
$T$,
\end{itemize}
such that, if
\begin{itemize}
\item $s_1, s_2, t_1$ and $t_2$ are elements of $(\Ffp\Ufp T)_n$,
\item $\delta_1: \hat s_1 \to \hat t_1$ and $\delta_2: \hat s_2 \to \hat t_2$
are coherence morphisms,
\item $\epsilon(s_1) = \epsilon(s_2) \in S$ and $\epsilon(t_1) = \epsilon(t_2)
\in S$,
\end{itemize}
then $\delta_1 = \delta_2$.
\end{defn}

This definition can be recast in terms of strict algebras for a \dco.

By judicious choice of Laplaza set, one can recover the classical notion of
symmetric monoidal category and Laplaza's categorification of the theory of
rigs.

\section{Non-algebraic definitions}
\index{Rosick\'y, Ji\v r\'i}
\index{Leinster, Tom}

Various definitions have appeared that are inspired by the notions of homotopy
monoids etc. in topology.
In \cite{hty_opd}, Leinster proposes a definition of a ``homotopy $P$-algebra
in $M$'' for any plain operad $P$ and any monoidal category $M$; his shorter
paper \cite{hty_mon} explores this definition in the case $P = 1$.
Related (but more general) is Rosick\'y's work described in \cite{rosicky}.

These definitions stand roughly in relation to ours as do the ``non-algebraic''
definitions of $n$-category in relation to the ``algebraic'' ones: see
\cite{cheng+lauda}.

\bibliographystyle{halpha}
\bibliography{thesis}
\printindex
\addcontentsline{toc}{chapter}{Index}

\end{document}